\documentclass[11pt]{amsart}

\usepackage[active]{srcltx}
\usepackage[utf8]{inputenc}
\usepackage[T1]{fontenc}
\usepackage[hcentering,vcentering,hscale=0.7,vscale=0.8]{geometry}
\usepackage{hyperref}
\usepackage{array}
\usepackage{enumitem}

\usepackage{tikz}
\usepgfmodule{plot}
\usetikzlibrary{hobby}

\usepackage{amsmath, amsthm, amssymb, amscd}
\newtheorem{theorem}{Theorem}[section]
\newtheorem{corollary}[theorem]{Corollary}
\newtheorem{lemma}[theorem]{Lemma}
\newtheorem{proposition}[theorem]{Proposition}
\newtheorem{maintheorem}{Theorem}

\theoremstyle{remark}
\newtheorem{definition}[theorem]{Definition}
\newtheorem{remark}[theorem]{Remark}

\newtheorem{example}[theorem]{Example}
\newtheorem{question}[theorem]{Question}

\usepackage{mathabx} %\righttoleftarrow
\newcommand{\A}{\ensuremath{\mathbb{A}}}

\newcommand{\C}{\ensuremath{\mathbb{C}}}
\newcommand{\D}{\ensuremath{\mathbb{D}}}
\newcommand{\E}{\ensuremath{\mathbb{E}}}
\renewcommand{\H}{\ensuremath{\mathbb{H}}}
\newcommand{\N}{\ensuremath{\mathbb{N}}}

\newcommand{\R}{\ensuremath{\mathbb{R}}}
\renewcommand{\SS}{\ensuremath{\mathbb{S}}}
\newcommand{\T}{\ensuremath{\mathbb{T}}}

\newcommand{\Z}{\ensuremath{\mathbb{Z}}}

\newcommand{\JJ}{\mathcal{J}}
\newcommand{\RR}{\mathcal{R}}
\newcommand{\ra}{\varrho}
\newcommand{\Dd}{\mathbb{D}}
\newcommand{\Sig}{\Sigma}
\newcommand{\eps}{\varepsilon}
\newcommand{\Rem}{\mathcal{R}}
\newcommand{\ddd}{\mathrm{d}}
\newcommand{\Wh}{\widehat{W}}
\DeclareMathOperator{\Ree}{Re}
\DeclareMathOperator{\Imm}{Im}
\DeclareMathOperator{\tr}{tr}
\DeclareMathOperator{\Area}{Area}
\DeclareMathOperator{\sgn}{sgn}
\DeclareMathOperator{\dist}{dist}

\title[On the dynamics of two vortices]{On the dynamics of two vortices\\ in a domain}
\author[J.~Fejoz]{Jacques \textsc{Fejoz}}
\address{CEREMADE, Université Paris Dauphine -- PSL, Paris, France \\
  and LTE, Observatoire de Paris -- PSL, Paris, France}
\email{jacques.fejoz@dauphine.fr}
\thanks{Corresponding author: Jacques Fejoz, \texttt{jacques.fejoz@dauphine.fr}.}

\author[Wang Q.]{\textsc{Wang} Qun}
\address{Department of Mathematical and Computational Sciences,
  University of Toronto Mississauga, Mississauga, Ontario, Canada}
\email{qun.wang@utoronto.ca}

\begin{document}

\begin{abstract}
  We study the dynamics of two point vortices in a bounded, finitely connected planar domain $D$. The global structure of the energy levels is governed by two properties: the opposite of the Hamiltonian is strictly plurisubharmonic whatever the signs of the intensities, whereas the energy levels are compact only when the intensities have opposite signs. In the opposite-sign problem, every regular energy level is a closed, Stein fillable contact $3$-manifold, and almost every energy level carries a periodic orbit. In the same-sign problem, the energy levels have one end per boundary component of $D$, at quadruple collisions of the vortices with their images across $\partial D$; yet the flow is complete, and these collisions are not even approached asymptotically. We then analyse three perturbative regimes --- near interior binary collisions, near the boundary of $D$, and near quadruple collisions, which we desingularize by a quasi-homogeneous blow-up --- and obtain in each of them periodic orbits and invariant tori carrying quasi-periodic motions. A dipole (opposite intensities) behaves differently: at leading order, it moves along straight lines and reflects on $\partial D$ like a billiard ball, up to an explicit shift along the boundary. Finally, in the same-sign problem at high energy, the energy levels are of contact type and, when the Robin function has a single critical point (for instance in strictly convex domains), they carry a Birkhoff annulus of section, whose return map is a twist map.
\end{abstract}

\maketitle

\tableofcontents

\section{Introduction}
\label{sec:intro}

Point vortices are singular solutions of Euler's equations for a planar perfect incompressible fluid, whose vorticity is concentrated at finitely many points. Since the works of Helmholtz, Kirchhoff and Routh, one knows that the positions of $N$ vortices obey a Hamiltonian system of ordinary differential equations; in a domain with boundary, the Hamiltonian is built from the hydrodynamic Green function of the domain, as shown by Lin~\cite{Lin:1941,Lin:1941:2}. We refer to \cite{Flucher:1999, Marchioro:1994, Newton:2001} and references therein for background. This article is devoted to the dynamics of $N=2$ vortices in a bounded domain, which is the first case where the dynamics is not expected to be integrable in general; non-integrability is proved in~\cite{Fejoz:2026} for generic deformations of the disc and equal, or nearly equal, intensities.

Let $D\subset\C$ be a bounded, finitely connected domain with $C^\infty$ boundary $\partial D=\partial_1D\sqcup\dots\sqcup\partial_mD$. Let $g(z,\zeta)=-\log|\zeta-z|+h(z,\zeta)$ be its hydrodynamic Green function, $r(z)=h(z,z)$ its Robin function and $q=e^r$ its harmonic radius, which vanishes on $\partial D$ (precise definitions are given in section~\ref{sec:introduction}). Two vortices $z_1,z_2\in D$ of intensities $\gamma_1,\gamma_2\in\R_*$ move according to the Hamiltonian
\[
  H(z_1,z_2)=\gamma_1\gamma_2\,g(z_1,z_2)+\tfrac12\bigl(\gamma_1^2\,r(z_1)+\gamma_2^2\,r(z_2)\bigr)
\]
on the configuration space $C_2(D)=\{(z_1,z_2)\in D^2,\ z_1\neq z_2\}$, endowed with the symplectic form $\omega=\frac1{2i}\sum_j\gamma_j\,\ddd\bar z_j\wedge\ddd z_j$. In section~\ref{sec:introduction} we recall why this is so, for any number of vortices (theorem~\ref{thm:H}); our proof relies on a dynamical argument, as an alternative to existing presentations. For a single vortex the Hamiltonian is $\frac{\gamma^2}2r$ and the dynamics is integrable; section~\ref{sec:1-v} reviews it, and relates the curl of the one-vortex flow to the Bergman kernel of $D$ and to Suita's conjecture: when $D$ is simply connected this curl equals $4/q^2$, hence is a first integral of the flow, and this identity characterizes simply connected domains (corollary~\ref{cor:suita} and proposition~\ref{prop:converse-suita}).

We write $\Sig_E=H^{-1}(E)$ for the energy levels, and
\[
  \Gamma=\gamma_1+\gamma_2,\qquad\kappa=\frac{\Gamma^2}{2\gamma_1\gamma_2}
\]
(so that $\kappa\geq2$ if $\gamma_1\gamma_2>0$).

\medskip\noindent\textbf{Global structure.} Two properties of different natures govern the energy levels. The first one is a symplectic convexity property, due to Flucher and Gustafsson: $-H$ is strictly plurisubharmonic, whatever the signs of the intensities (theorem~\ref{thm:Hpsh}). The second one is compactness, which does depend on the signs.

\begin{maintheorem}[Global structure]
  \label{thm:A}
  Let $\gamma_1,\gamma_2\in\R_*$.
  \begin{enumerate}
  \item If $\gamma_1\gamma_2<0$, the energy levels are compact. Every regular energy level is a closed, tight, Stein fillable contact $3$-manifold, and almost every energy level carries a periodic orbit.
  \item If $\gamma_1\gamma_2>0$, the energy levels are not compact. Each of them has exactly $m$ ends, diffeomorphic to $\T^2\times(0,1)$, one per boundary component of $D$, along which the two vortices tend to a common point of $\partial D$. Yet the flow is complete, and each end contains invariant tori which separate the energy level and whose relative measure tends to $1$; consequently, every solution stays in a compact subset of $C_2(D)$.
  \end{enumerate}
\end{maintheorem}

\begin{figure}[ht]
  \centering
  \includegraphics[width=.86\textwidth]{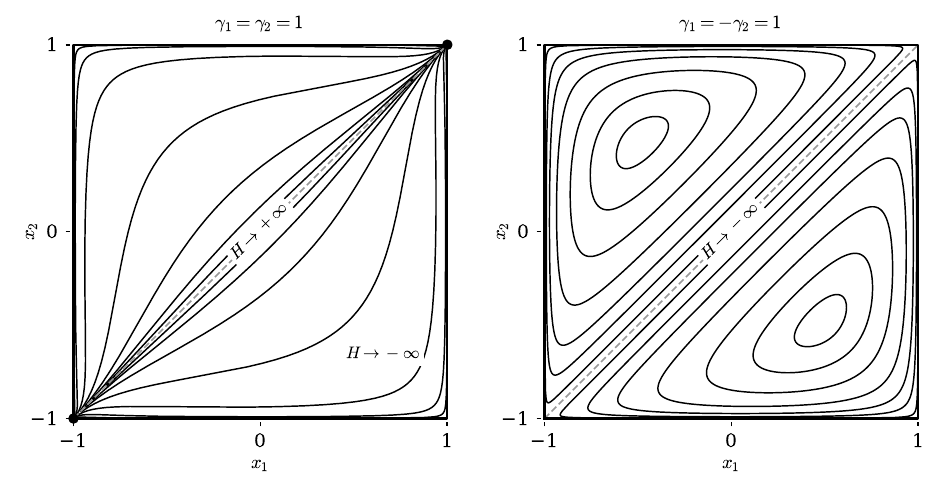}
  \caption{Theorem~\ref{thm:A} in the unit disc: traces of the energy levels on the slice $\{z_1=x_1,\,z_2=x_2\}$, $x_1,x_2\in(-1,1)$, of the configuration space (a slice which is not invariant by the flow). Left, same signs ($\kappa=2$; $E=-3,-2,\dots,3$): every energy level reaches the two corners $x_1=x_2=\pm1$, which are quadruple collisions, and is tangent there to the diagonal of binary collisions (dashed), the distance between the vortices being of the order of the square of their distance to the boundary. Right, opposite signs: the energy levels are compact.}
  \label{fig:levels}
\end{figure}

See theorems~\ref{thm:compact} and~\ref{thm:HZ} and proposition~\ref{prop:fillable} for the first item, and proposition~\ref{prop:ends} and theorems~\ref{thm:complete} and~\ref{thm:kam} for the second one. In the first item, the contact structure is the $CR$ structure of the level; the vortex flow is not its Reeb flow (remark~\ref{rem:notreeb}), which is why one only gets periodic orbits on almost every level. The ends of the second item are \emph{quadruple collisions}: both vortices collide on the boundary, hence with their two images across $\partial D$. Figure~\ref{fig:levels} shows a planar slice of the energy levels in both cases.

\medskip\noindent\textbf{Perturbative regimes.} The boundary of $C_2(D)$ consists of binary collisions and of configurations where at least one vortex lies on $\partial D$. We analyse the dynamics in the three corresponding regimes: near interior binary collisions (section~\ref{sec:collision}), near $\partial D$ away from collisions (section~\ref{sec:near-domain-boundary}), and near quadruple collisions (section~\ref{sec:quadruple}). We state the results in the first regime only, the other two being similar in spirit. Assume that $\Gamma\neq0$ and let
\[
  \eps=e^{-E/(\gamma_1\gamma_2)} ,
\]
so that the vortices are close to each other, away from $\partial D$, when $\eps$ is small. Let $\Theta=q^{2\kappa}$. An annulus $\mathcal U\Subset D$ foliated by closed regular level curves of $r$ is called \emph{non-isochronous} if the second derivative of $\Theta$ with respect to the area enclosed by these curves does not vanish, i.e.\ if the frequency of the Hamiltonian flow of $\Theta$ varies monotonically, with a non-vanishing derivative, across $\mathcal U$; this always holds near a non-degenerate maximum of $r$ (lemma~\ref{lem:auto}).

\begin{maintheorem}[Near interior binary collisions]
  \label{thm:B}
  Assume $\Gamma\neq0$ and let $K\Subset D$ be compact. For $\eps$ small enough:
  \begin{enumerate}
  \item To every non-degenerate critical point $Z_*$ of $r$ there corresponds a periodic orbit in $\Sig_E$, depending real analytically on $\eps$, along which the two vortices, at distance $\eps\,q(Z_*)^\kappa(1+O(\eps^2))$ from each other, rotate about a point $O(\eps^4)$-close to $Z_*$ with period $2\pi\eps^2q(Z_*)^{2\kappa}/|\Gamma|\,(1+O(\eps^2))$. It is elliptic if $Z_*$ is a local maximum of $r$, and hyperbolic if $Z_*$ is a saddle. In particular there is at least one elliptic family if the maximum of $r$ is non-degenerate (e.g.\ if $D$ is strictly convex), and at least $m-1$ hyperbolic ones if $r$ is a Morse function.
  \item Over every non-isochronous annulus $\mathcal U\Subset D$, the energy level $\Sig_E$ contains real analytic invariant Lagrangian tori carrying quasi-periodic motions, whose relative measure is $1-O(e^{-c/\eps^2})$ for some $c>0$. They separate $\Sig_E$, so that the solutions lying between two of them neither collide nor reach $\partial D$; between two consecutive tori there are periodic orbits and Mather sets of every intermediate rotation number.
  \item Along every solution whose centre of vorticity stays in $K$ over the time interval $[0,t]$, the distance between the vortices varies by at most $O(\eps^3)$, as long as $|t|\leq e^{c'/\eps^2}$, for some $c'>0$.
  \end{enumerate}
\end{maintheorem}

See theorems~\ref{thm:periodic} and~\ref{thm:qp}, proposition~\ref{prop:po} and corollary~\ref{cor:adiabatic}. These results do not depend on the signs of the intensities. Here is a description of the corresponding motions. As $\eps\to0$, the two vortices form a tight pair, of separation $\rho=|z_1-z_2|\simeq\eps\,q(Z)^{\kappa}$, where $Z=(\gamma_1z_1+\gamma_2z_2)/\Gamma$ is their centre of vorticity, and two motions are combined:
\begin{itemize}
\item a fast spin: the pair rotates about $Z$ like an isolated pair in the plane, in the direction given by the sign of $\Gamma$, at the angular velocity $\Gamma/\rho^2\simeq\Gamma\,\eps^{-2}q(Z)^{-2\kappa}$; the vortices $z_1$ and $z_2$ describe circles of radii $|\gamma_2/\Gamma|\,\rho$ and $|\gamma_1/\Gamma|\,\rho$ about $Z$ (on the same side of $Z$ when $\gamma_1\gamma_2<0$), slightly deformed by the strain of the images, so that $\rho$ oscillates at twice the spin frequency with relative amplitude $O(\eps^2)$;
\item a slow drift: the centre of vorticity moves like a single vortex of intensity $\Gamma$, that is, at the velocity $\frac\Gamma2\nabla^\perp r$, which does not depend on $\eps$, along the level curves of the Robin function $r$.
\end{itemize}
A lap of $Z$ around a closed level curve $\mathcal O$ of $r$ thus takes a time of order $1$, during which the pair makes
\[
  \frac{1}{\pi\,\eps^{2}q_{\mathcal O}^{2\kappa}}\oint_{\mathcal O}\frac{\ddd s}{|\nabla r|}\,\bigl(1+O(\eps^2)\bigr)
\]
turns, $q_{\mathcal O}$ being the value of the harmonic radius on $\mathcal O$. Since $q$ is constant on $\mathcal O$, so are, along a lap, the separation and the spin rate; but they vary from one level curve to the next, and this variation is the twist which makes the two frequencies non-resonant for most level curves.

\begin{figure}[ht]
  \centering
  \includegraphics[width=.8\textwidth]{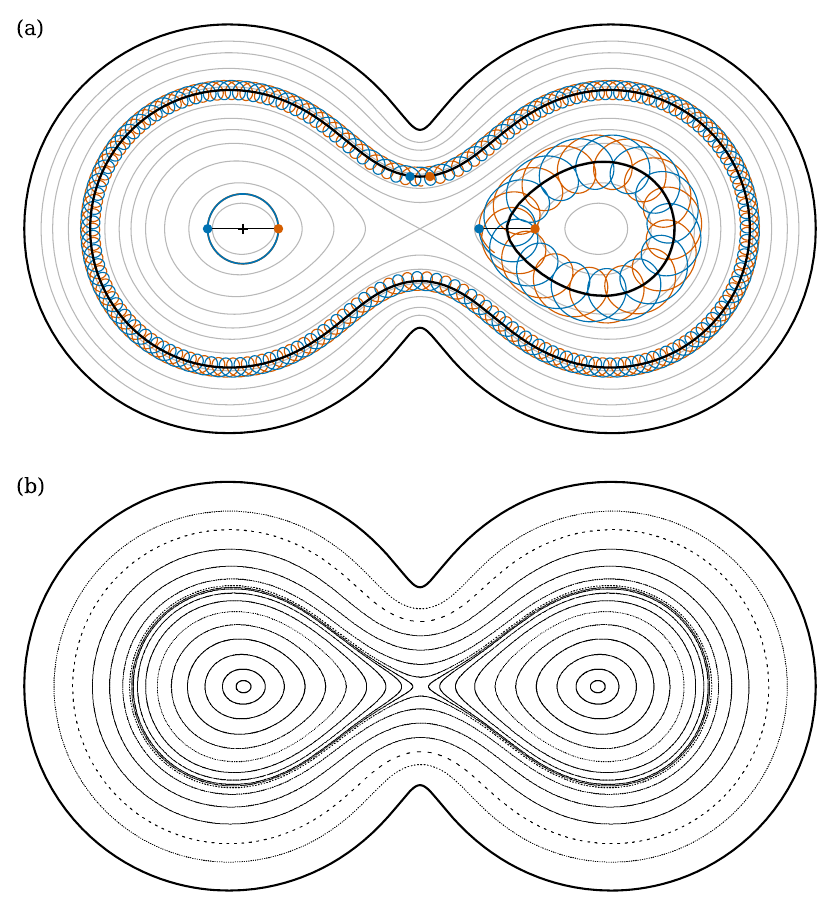}
  \caption{Theorem~\ref{thm:B} in the image $D$ of the unit disc by $w\mapsto w/(1-0.6\,w^2)$, whose Robin function has two maxima and a saddle, for $\gamma_1=\gamma_2=1$ and $\eps=0.25$ (a value which is not very small, chosen so that the spin of the pair is visible). (a)~Three solutions: the two vortices (in colours; the dots are their initial positions) and their centre of vorticity (in black), over one lap of the latter, against the level curves of $r$ (in grey). On the left, the elliptic periodic orbit associated with a maximum of $r$ (cross). (b)~Iterates of the return map of the energy level to the section $\{\arg(z_1-z_2)=0\}$, i.e.\ positions of the centre of vorticity each time the pair has completed a turn: the invariant curves are close to the level curves of $r$, around two elliptic fixed points and one hyperbolic fixed point.}
  \label{fig:pair}
\end{figure}

The periodic orbits of theorem~\ref{thm:B} are the motions in which the centre of vorticity does not move: the pair spins in place, close to an equilibrium $Z_*$ of a single vortex of intensity $\Gamma$, and it is the rotation which stabilizes the pair against the local strain created by the images. The invariant tori carry the generic motion, quasi-periodic with the two frequencies above: the centre of vorticity follows a level curve of $r$ up to $O(\eps^2)$, the pair spins incommensurably faster, and it remains bound forever. See figure~\ref{fig:pair}.

The periodic orbits and Mather sets of the second item describe what happens in the thin layers left between invariant tori. Record, at each turn of the pair (say each time $z_1-z_2$ points in a given direction), the position of $Z$ along its circuit; the rotation number $\varpi$ is the mean fraction of a lap per turn, of order $\eps^2$. For $\varpi=p/k$ rational, the periodic orbits are resonant motions: after exactly $k$ turns of the pair, its centre of vorticity has made exactly $p$ laps and the configuration repeats; those given by Aubry--Mather theory are moreover well ordered: the successive positions of $Z$ along its circuit come in the same cyclic order as under the rigid rotation by $p/k$. For $\varpi$ irrational, the Mather set is either an invariant torus, possibly non-smooth, or a \emph{cantorus}, the remnant of a destroyed torus. On a cantorus, the motion still has two well-defined frequencies but it is not quasi-periodic: the successive positions of $Z$ along its circuit are of the form $\Phi(\phi_0+n\varpi)$, where $\Phi$ is a non-decreasing function with jumps (Aubry's hull function), so that, whenever the pair has a given orientation, its centre of vorticity lies on a Cantor subset of its circuit. Since the complement of the tori has a measure which is exponentially small with respect to $1/\eps^2$, so is the thickness of the layers containing these motions: geometrically, the motions of the Mather sets are indistinguishable from those of the neighbouring tori, and differ from them only by the combinatorics of the succession of spins and laps.

\medskip The same mechanism --- a single fast angle, whose action, read on the energy level, generates the return map to a section (proposition~\ref{prop:universal}) --- is at work \emph{near the boundary of $D$}, at very negative energies. Three cases occur there, according to the position of the vortices (section~\ref{sec:near-domain-boundary}).
\begin{itemize}
\item One vortex is close to a boundary component $\partial_jD$ and the other one is not (regime~I). The first vortex skims along $\partial_jD$ in a thin layer, of width $d_1$ exponentially small with respect to $|E|$, at speed $\simeq|\gamma_1|/(2d_1)$. Meanwhile the second vortex moves, in a time of order $1$, along a level curve of the effective potential $V_j=\frac{\gamma_2^2}2r+\gamma_1\gamma_2C_j$, where $C_j$ is the limit of $g(\cdot,z_1)$ as $z_1$ tends to any point of $\partial_jD$: it feels the first vortex as if the latter had been absorbed by the wall. The width $d_1$ is constant along each level curve of $V_j$, and its variation from one level curve to the next provides the twist. The non-degenerate critical points of $V_j$ give rise to periodic orbits, and non-isochronous annuli of level curves of $V_j$ to invariant tori (theorems~\ref{thm:poI} and~\ref{thm:kamI}).
\item The two vortices are close to two different boundary components (regime~III). Each of them skims along its own component, at distance $d_i$ and speed $\simeq|\gamma_i|/(2d_i)$. The two frequencies are of the same order, their ratio varies monotonically along the energy level whatever the signs of the intensities, and the two vortices, which stay away from each other, interact only at second order: this is a classical, two-frequency nearly integrable system, most of whose motions are quasi-periodic (theorem~\ref{thm:III}).
\item The two vortices are close to the same boundary component (regime~II). They drift along it at different speeds, hence overtake each other and undergo a strong interaction once per relative period. This regime is not perturbative; we only describe its local structure, which is reminiscent of the co-orbital motions of celestial mechanics.
\end{itemize}

\emph{Near quadruple collisions}, in the same-sign problem, the motions are the analogue of those of theorem~\ref{thm:B}, at any energy, the small parameter now being the harmonic radius $q$ of the centre of vorticity, which is of the order of the distance of the pair to the wall. A tight pair, of separation $\rho\simeq\eps\,q^{\kappa}\ll q$, spins at angular velocity $\Gamma/\rho^2$, while its centre of vorticity drifts along $\partial_jD$ like a single vortex of intensity $\Gamma$, at a speed of order $\Gamma/q$. Both frequencies blow up as $q\to0$, the spin much faster than the drift: the number of turns per lap is of order $\eps^{-2}q^{1-2\kappa}$, with $2\kappa-1\geq3$. We desingularize this regime by the quasi-homogeneous blow-up $z_1-z_2=q(Z)^\kappa\zeta$, for which the logarithmic divergences of $H$ cancel exactly (proposition~\ref{prop:cancel}). In contrast with the collision manifolds of the Newtonian $N$-body problem, the collision manifold thus obtained is normally degenerate; it is nevertheless accumulated by the invariant tori of theorem~\ref{thm:A} (theorem~\ref{thm:kam}).

\begin{figure}[ht]
  \centering
  \includegraphics[width=\textwidth]{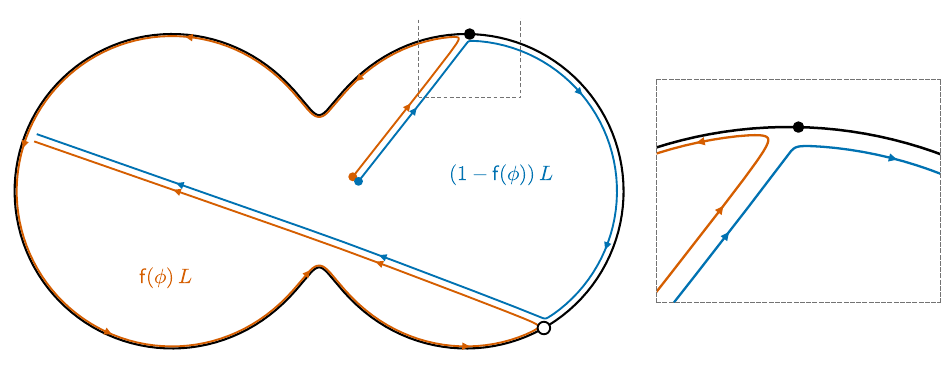}
  \caption{A dipole ($\gamma_1=-\gamma_2=1$, $\eps=0.06$) in the domain of figure~\ref{fig:pair}, over one wall event. Starting from the two coloured dots, the dipole moves along a nearly straight line to the boundary (black dot), where it splits. The two vortices then slide along $\partial D$ in opposite directions, each of them at a constant distance from the boundary, and at a speed which is inversely proportional to this distance. These two distances differ, as the close-up of the dashed window shows (right): here the vortex of positive intensity (in orange) is about three times closer to the boundary than the other one, hence about three times faster. When the two vortices meet again, having covered together the whole perimeter $L$ of the domain, the former has thus travelled the fraction $\mathsf f(\phi)\simeq0.76$ of it, and the latter the fraction $1-\mathsf f(\phi)$. A dipole re-forms there, close to the point predicted by theorem~\ref{thm:limitmap} (white dot), and it leaves the boundary at the angle $\phi$ at which it arrived.}
  \label{fig:dipole}
\end{figure}

\emph{A dipole} ($\Gamma=0$, say $\gamma_1=-\gamma_2=1$) behaves differently (section~\ref{sec:dipole}): its fast motion is a translation, not a rotation. The separation $\rho$ of the two vortices is almost constant by conservation of energy, and the dipole moves at speed $\simeq1/\rho$ in the direction orthogonal to its axis, along nearly straight lines. When it reaches the wall it splits: the two vortices slide along $\partial_jD$ in opposite directions, at different heights, and a dipole re-forms where they meet again. At leading order, the dipole thus reflects on $\partial D$ like a billiard ball, up to an explicit shift along the boundary, which depends only on the angle of incidence (theorem~\ref{thm:limitmap}); we do not estimate the error terms in this last regime, but the limit law agrees well with numerical solutions (figure~\ref{fig:dipole}).

\medskip\noindent\textbf{Contact type and Birkhoff annulus.} Let us come back to the same-sign problem. Whether the vortex flow on an energy level is a Reeb flow, i.e.\ whether the level is of contact type for $\omega$ itself, and whether it has a global surface of section, are natural questions (questions~\ref{qu:contacttype} and~\ref{qu:Birkhoff}), which we answer at high energy.

\begin{maintheorem}[High energy]
  \label{thm:C}
  Assume $\gamma_1,\gamma_2>0$.
  \begin{enumerate}
  \item Every energy level is of contact type near its ends, and there is $E_c$ such that $\Sig_E$ is of contact type for every $E\geq E_c$. In nearly round domains, this holds for every $E$.
  \item Assume that $D$ is simply connected and that $r$ has a unique critical point, which is non-degenerate (for instance, that $D$ is strictly convex). For $E$ large enough, the elliptic periodic orbit $P_E$ of theorem~\ref{thm:B} bounds an annulus $\mathcal A_E\subset\Sig_E\setminus P_E$, whose other end lies at the quadruple collisions, and which every solution other than $P_E$ crosses transversally, infinitely many times. The return map of $\mathcal A_E$ is an area-preserving twist map, whose rotation number goes from a finite value at $P_E$, of order $\eps^{-2}$, to infinity at the collision end; it has periodic orbits of every larger rational rotation number.
  \end{enumerate}
\end{maintheorem}

See proposition~\ref{prop:contact-highE}, remark~\ref{rem:nearlyround}, theorem~\ref{thm:annulus} and corollary~\ref{cor:annulus-PB}. In the blown-up energy level, the second boundary component of $\mathcal A_E$ is a periodic orbit of the collision manifold. Remark~\ref{rem:lowE-hard} explains why the low energy regime is harder.

\medskip\noindent\textbf{Notations.} Each regime has its own natural coordinates and its own small parameter. Rather than aiming at globally unambiguous notations, we have kept them local: the names of the variables (such as $Z$, $u$, $I$, $\theta$, $\JJ$, $\varphi$, $A$, $\psi$) and of the small parameters ($\eps$, $\delta$, $\lambda$) are defined anew in each section, and their meaning may change from one section to another, according to the region of the phase space under study. For instance, $\eps=e^{-E/(\gamma_1\gamma_2)}$ in sections~\ref{sec:collision} and~\ref{sec:highE}, whereas $\eps$ is a ratio of frequencies in section~\ref{sec:quadruple} and $\eps=e^{E}$ in section~\ref{sec:dipole}. Only the objects attached to the domain and to the problem ($D$, $m$, $g$, $h$, $r$, $q$, $\gamma_j$, $\Gamma$, $\kappa$, $H$, $\omega$, $\Sig_E$) keep the same meaning throughout.

\medskip\noindent\textbf{Organization.} Section~\ref{sec:introduction} derives the Hamiltonian equations of $N$ vortices and recalls the properties of the Green function. Section~\ref{sec:1-v} deals with one vortex. Section~\ref{sec:first-properties} establishes the global properties of the $2$-vortex problem. Sections~\ref{sec:collision}, \ref{sec:near-domain-boundary} and~\ref{sec:quadruple} are devoted to the three perturbative regimes, section~\ref{sec:highE} to contact type and to the Birkhoff annulus at high energy, and section~\ref{sec:dipole} to the dipole.

\section{Point vortices in a domain}
\label{sec:introduction}

\begin{figure}[h]
  \centering
  \includegraphics[width=5cm]{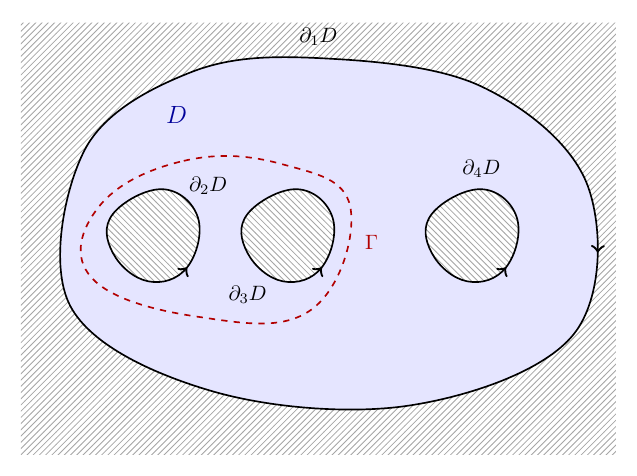}
  \caption{Example of a $4$-connected domain $D$ with a curve $\Gamma$ (a possible level curve of the Robin function) which is not homotopic to a boundary component}
  \label{fig:domain4}
\end{figure}

Let $D \subset \C$ be a $C^\infty$-smooth, bounded and finitely connected domain. Write $\partial D=\partial_1D \sqcup \dots \sqcup \partial_mD$, where $m$ is the degree of connectivity and $\partial_1 D$ is the outer boundary component; the Euler characteristic of $D$ is $\chi(D)=2-m$ (see an example in figure~\ref{fig:domain4}). We denote by $L_j$ the length of $\partial_jD$.

According to a theorem of Koebe~\cite[Chap. VI]{Goluzin:1969}, any bounded domain which is finitely connected and such that the connected components of the complement are not reduced to points, is conformally equivalent to such a $D$. So, smoothness is not a key hypothesis. Boundedness is not a key hypothesis either (the point at infinity on the Riemann sphere playing no special role), but will avoid discussing boundary conditions at infinity.

\medskip Consider a fluid moving in $D$, with time-dependent velocity vector field $u = u(t,\zeta) = u_t(\zeta)$. Using the Euclidean scalar product of $\C$, we define the associated velocity $1$-form $v = u^\flat$. Suppose that
\begin{itemize}
\item the fluid cannot pass through the boundary:
  \[\tag{B} v(t,z) \cdot n(z) = 0 \quad (t \in \R, \; z \in \partial D),\]
  where $n(z)$ is a vector normal to $\partial D$ at $z$
\item the fluid is incompressible:
  \[\tag{I} \delta v=0,\]
  where $\delta = *d*$ is the codifferential (in the spatial variables) and $*$ the Hodge operator 
\item the fluid is perfect and satisfies Euler's equation
  \[\tag{E} \partial_t v + L_{u} v = d\left( \frac{\|v\|^2}{2} -p \right),\]
  where $p$ is the pressure, and the Lie derivative of the velocity $1$-form $v$ in the direction of $u$ is related to the covariant derivative via Koszul's formula: $\nabla_{u} v = L_{u} v- d\left( \frac{\|v\|^2}{2} \right)$.
\end{itemize}

\begin{lemma}[Helmholtz's theorem]
  The circulations
  \[c_k = \oint_{\partial_k D} v\]
  do not depend on time ($k=1,...,m$).
\end{lemma}

\begin{proof}
  Using equations (E), (B) (which entails that $u$ is tangent to $\partial_k D$) and Stokes formula, we see that
  \[\partial_t c_k = \oint_{\partial_k D} \partial_t v = \oint_{\partial_k D} -i_u \, dv =0.\]
\end{proof}

Introduce the vortex function $w = *dv$. Euler's equation (E) for some pressure $p$ is equivalent to the left hand side being exact, i.e. to the equation obtained by taking the exterior derivative:
\[\tag{dE} \partial_t w + dw \cdot u = 0\]
(meaning that the vorticity $w$ is transported into itself along the flow of $u$), provided the cohomological obstructions vanish:
\[\tag{O} \oint_{\partial_k D} \left( \partial_t v + L_u v \right) = 0 \quad (k=1,...,m).\]

\begin{remark}
  According to the Helmholtz-Hodge decomposition, $v$ uniquely decomposes into the sum of exact, coexact and harmonic parts. Because of incompressibility (I), the exact part vanishes. Due to the non-linearity in (E), in general the evolutions of the two remaining parts are coupled. In the vortex subproblem, as defined below, a miracle occurs and the harmonic part stays null.
\end{remark}

\medskip Following Koebe~\cite{Koebe:1914} and Lin~\cite{Lin:1941,Lin:1941:2}, we call \emph{hydrodynamic Green function} of $D$ of periods\footnote{also called circulations or fluxes} $c_1,...,c_m \in \R$ a function
\[g : (D \times \bar D) \setminus \Delta \to \R, \quad (z,\zeta) \mapsto g(z,\zeta) = g_z(\zeta)\]
(where $\Delta$ is the diagonal of $D^2$) of class $C^2$ such that, for all $z \in D$ (the ``source'' point),
\[
  \tag{G}
  \begin{cases}
    -\Delta_\zeta (g_z) = 2\pi \delta_z\\
    j_{\partial D}^* d_\zeta (g_z) = 0\\
    \oint_{\partial_k D} *d_\zeta (g_z) = c_k &(\forall k)\\
    \oint_{\partial D} g_z (*d_\zeta (g_z)) = 0.
  \end{cases}
\]
% $$\Delta (\log|z - w|) = 2\pi \delta_w$$
where $\Delta_\zeta$ and $d_\zeta$ are the Laplacian and the real derivative with respect to $\zeta$, $\delta_z$ is the Dirac distribution supported by $\{z\}$ and $j_{\partial D} : \partial D \to \bar D$ is the inclusion map. The sign choice in the first equation of (G) is useless here, except for being consistent with the usual convention, e.g. \cite{Flucher:1999}.

\begin{proposition}
  The Green function $g$ has the following properties:
  \begin{enumerate}
  \item (Existence) $g$ exists iff $\sum c_k = -2\pi$.
  \item (Uniqueness) $g$ is unique.
  \item (Symmetry) $g$ is symmetric:
    \[g(z,\zeta) = g(\zeta,z) \quad (\forall (z, \zeta) \in D^2 \setminus \Delta)\]
  \item (Regularity) $g$ is real analytic over $D^2 \setminus \Delta$ and has the same regularity on $D \times \partial D$ (in particular, $C^\infty$ or analytic regularity) as $\partial D$.
  \item (Conformal invariance) If $\varphi : D \to D'$ is a conformal map,
    \[g_{D}(z,\zeta) = g_{D'}(\varphi(z),\varphi(\zeta)).\]
    In particular, the logarithmic pole of the Green's function is conformally invariant.
  \end{enumerate}
\end{proposition}

There are different variants for Green functions and the terminology itself is not standardized. So we give references or hints for the proofs of these properties. 

\begin{proof}[Hints of proof]
  For (1), (2) (3), see \cite[Chap. 4]{Schiffer:1954} or \cite[15.2]{Flucher:1999}. Both these references start with the Dirichlet Green function $g_0$ (vanishing identically on all boundary components, see definition~\ref{def:Robin}). Then, in order to achieve arbitrary constants on boundary components and prescribed periods, some linear combination of the domain's harmonic measures (under the name of differentials of the third kind with prescribed imaginary periods, for the second reference) is added; note that this linear combination is $C^\infty$-smooth on $\bar D$. 

  Note that the last equation of (G) precisely normalizes the gauge so as to get a symmetric function. 

  (4) Regularity in $D$ is a consequence of elliptic regularity~\cite{Gilbarg:2001}. Boundary regularity then follows from the Schwarz reflection principle, as discussed for example in~\cite[Section 6.3]{Ahlfors:1979}.

  (5) The book \cite{Schiffer:1954} shows the conformal invariance of all such domain functionals. For the sake of clarity, we give a direct proof here. In order to show that $g(z,\zeta) = g_{D'}(\varphi(z),\varphi(\zeta))$ acts as the hydrodynamic Green function on $D$, we check the defining properties of (G) one by one:

  -- Laplace's equation: let  $z \in D$ and $z' = \varphi(z) \in D'$. We have
  \[-\Delta_{\zeta'}\, g_{D'} (z',\cdot) = 2\pi\delta_{z'}.\]

  Since $\varphi$ is conformal, using $\Delta = 4 \partial_{\bar\zeta}\partial_\zeta$ we see that, for a $C^2$ real function $f$ on $D'$, 
  \[\Delta (f \circ \varphi) = |\varphi'|^2 \, (\Delta f) \circ \varphi.\]
  The case when $f(\zeta) = g_{D'}(z',\zeta)$ yields
  \[\Delta_\zeta\, g_z = |\varphi'|^2 \, \big(\Delta g_{D'}(z',\cdot)\big) \circ \varphi = - 2\pi |\varphi'|^2 \, \delta_{z'}\circ \varphi.\]
  The map $\varphi$ is a local diffeomorphism near $z$, with real Jacobian $|\varphi'(z)|^2$ at $z$. The change-of-variables rule for the delta distribution therefore gives
  \[\delta_{z'}(\varphi(\zeta)) = \delta_{\varphi(z)}(\varphi(\zeta))
    = \frac{\delta_z(\zeta)}{|\varphi'(z)|^2},\]
  whence
  \[-\, \Delta_\zeta\, g_z = 2\pi\frac{|\varphi'(z)|^2}{|\varphi'(z)|^2}\, \delta_z
    = 2\pi \delta_z.\]

  -- Boundary constants: if $g_{D'}$ is constant on each boundary component of $D'$, then so is $g$ on each $\partial_k D$.

  -- Circulations: the circulation integral is invariant under conformal mappings, hence $\oint_{\partial_k D} *d_\zeta g_z = c_k$.

  - Symmetry being equivalent to the gauge normalization, $g_{D'}$ satisfies the gauge normalization. 

  By uniqueness, $g$ must be the hydrodynamic Green function for $D$.
\end{proof}

\begin{definition}
  \label{def:Robin}
  The \emph{regular part} $h$ of $g$ is the function \[h(z,\zeta) = g(z,\zeta) + \log |\zeta-z|;\]
  it is biharmonic, and thus real analytic over $D^2$ (as a consequence of Cauchy's estimates for harmonic functions). The \emph{Robin function} $r$ is the restriction of $h$ to the diagonal:
  \[r(z) = h(z,z), \quad z \in D,\]
  and the (hydrodynamic) \emph{harmonic radius} $q$ is
  \[q(z) = e^{r(z)}, \quad z \in D.\]
  The \emph{Dirichlet Green function} $g_0$ of $D$ is the function satisfying the first equation of (G) and vanishing on $\partial D$; we denote by $h_0$, $r_0$ and $q_0 = e^{r_0}$ its regular part, its Robin function and its harmonic radius, defined as above.
\end{definition}

When $m=1$, $g_0$ satisfies (G) with $c_1 = -2\pi$, so that, by uniqueness, $g = g_0$, $h = h_0$, $r = r_0$ and $q = q_0$. When $m\geq 2$, the periods $\oint_{\partial_kD} *d_\zeta\, g_0(z,\cdot) = -2\pi\omega_k(z)$, where $\omega_k$ is the harmonic measure of $\partial_kD$, depend on $z$, so that $g_0$ is not a hydrodynamic Green function. The index $0$ will always refer to the Dirichlet normalization; the objects $g_0$, $h_0$, $r_0$, $q_0$ only appear in auxiliary arguments and in the discussion of Suita's conjecture at the end of section~\ref{sec:1-v}.

\begin{corollary}
  \label{cor:conformal-r}
  Functions $h$ and $r$ are real analytic on $D^2$ and $D$ respectively. Moreover, if $\varphi : D \to D'$ is a conformal map,
  \[
    \begin{cases}
      h_D(z, w) = h_{D'}(\varphi(z), \varphi(w)) - \log \left| \frac{\varphi(z) - \varphi(w)}{z - w} \right|\\
      r_D(z) = r_{D'}(\varphi(z)) - \log|\varphi'(z)|.
    \end{cases}\]
\end{corollary}

There are few explicit examples. 

\begin{example}[Disc]
  \label{ex:disc}
  For the unit disc $\D = \{|z| < 1\}$ (with circulation $c_1 = -2\pi$ along the circle),
  $$g(z, \zeta) = \log \left| \frac{1-z\bar{\zeta}}{\zeta-z} \right|, \quad h(z, \zeta) = \log |1 - z\bar{\zeta}|, \quad 
  r(\zeta) = \log (1 - |\zeta|^2).$$
\end{example}

The next example is not bounded. Seen on the Riemann sphere, $\E$ is simply connected, so its Green function is that of the disc read through $z \mapsto 1/z$; the formula below corresponds to the other natural normalization, in which the circulation around the disc vanishes and the whole of the circulation $-2\pi$ is carried by the point at infinity. The two differ by the harmonic function $-\log|z\zeta|$.

\begin{example}[Exterior disc]
  For the exterior disc $\E = \{|z|>1\}$,
  $$g(z, \zeta) = \log \left| \frac{1-z\bar{\zeta}}{\zeta-z} \right| - \log |z\zeta|,$$
  while
  $$h(z, \zeta) = \log |1 - z\bar{\zeta}| - \log|z \zeta|, \quad r(\zeta) = \log \frac{|\zeta|^2 - 1}{|\zeta|^2}.$$
\end{example}

We now generalize the two prior examples to the annulus. Considering the lift of the Green function $g$ to the double cover of an annulus, one sees that $g$ may be expressed in terms of Jacobi's theta functions; for the precise formula below, one may map the annulus to an infinite rectangular strip using the logarithm, use the method of images (creating an infinite chain of sources and sinks), and take the exponential back to the annulus \cite{Nehari:1952,Crowdy:2020}.

\begin{example}[Annulus]
  \label{ex:annulus}
  Consider the annulus
  \[\A(\varrho,R) = \{ \varrho < |z| < R \},\]
  with periods $c_1$ on the outer circle and $c_2$ on the inner one, $c_1+c_2=-2\pi$. Writing $\mathsf q = \varrho/R$ for the modulus of the annulus and $b$ for a real constant (which we do not care for, as it is fixed by the gauge normalization),
  \[g(z, \zeta) = - \log \left| \frac{\vartheta_1\left(\frac{1}{2\pi i}\log\frac{z}{\zeta}, \mathsf q\right)}{\vartheta_1\left(\frac{1}{2\pi i}\log\frac{z\bar{\zeta}}{R^2}, \mathsf q\right)} \right| - \frac{c_2}{2\pi} \bigl( \log |z| + \log|\zeta| \bigr) + b,\]
  and, subtracting the pole,
  $$r(\zeta) = - \log \left| \frac{\vartheta_1'(0, \mathsf q)}{2\pi \,\vartheta_1\left(\frac{1}{2\pi i}\log\frac{|\zeta|^2}{R^2}, \mathsf q\right)} \right| + \left( 1 - \frac{c_2}{\pi} \right)\log|\zeta| + b.$$
  Letting $\mathsf q \to 0$ with $R=1$ and $c_2=0$, the Jacobi theta function degenerates into $\vartheta_1(v,\mathsf q) \sim 2 \mathsf q^{1/4} \sin(\pi v)$ and one recovers example~\ref{ex:disc}, as it must be.
\end{example}

Conformal invariance allows to significantly extend these examples. In particular, from the examples~\ref{ex:disc} and~\ref{ex:annulus} we get expressions of the hydrodynamic Green functions for all simply or doubly connected domains, in terms of Riemann's conformal representation.

\medskip We can now define vortex flows. Consider a generalized function on $D$ of the form 
\[\psi  = \sum_1^N \gamma_j g_{z_j(t)},\]
where the moving points $z_j(t)$ are to be determined, assuming they do not collide ($z_j(t) \neq z_k(t)$ if $j \neq k$). Let
\[\tag{V} u = 4 \Im \left( \partial_{\bar \zeta} \psi \, \partial_\zeta \right)\]
be the associated Hamiltonian vector field with respect to the standard symplectic form $\omega_0= \frac{1}{2i} \, d\bar\zeta \wedge d\zeta$ (with the sign convention $d\psi = + i_u \omega_0$). The corresponding velocity form is $v = -*d\psi$. Due to the definition of the Green function $g$, the vorticity function is
\[w = *dv = - \Delta \psi = 2\pi \sum \gamma_j \delta_{z_j(t)};\]
thus the moving points $z_j(t)$ are called \emph{vortices} and the coefficients $\gamma_j \in \R_*$ are \emph{vortex intensities} (or vortices, by metonymy); the circulation carried by the $j$-th vortex is $2\pi\gamma_j$. 

Recall that our fluid dynamical problem is defined by equations (B), (I) and (E), and that the latter is equivalent to (dE) together with vanishing obstructions (O). Let $\Delta = \{z \in \C^N, \; \exists j,k \; z_j=z_k\}$ be the fat diagonal of $\C^N$, and
\[C_N(D) = D^N \setminus \Delta\]
be the configuration space of $N$ pairwise distinct points in $D$. 

\begin{theorem}[Helmholtz, Kirchhoff, Routh, Lin~\cite{Lin:1941}]
  \label{thm:H}%
  Equations (B), (I) and (E) are satisfied iff the joint motion $z(t) = (z_1(t),...,z_N(t))$ of vortices in $C_N(D)$ is an integral curve of the Hamiltonian vector field $X_H$ associated with the Hamiltonian
  \[\tag{H} H(z) = \sum_{j < k} \gamma_j \gamma_k\, g(z_j, z_k) + \frac{1}{2} \sum_j \gamma_j^2\, r(z_j)\]
  and the symplectic form
  \[\omega = \frac{1}{2i} \sum_j \gamma_j\, d\bar z_j \wedge dz_j,  \qquad i_{X_H} \omega = dH.\]
\end{theorem}

In real notations, the expression of $\omega$ entails that $x_j = \Re z_j$ and $y_j = \Im z_j$ are conjugate. 

We will explain after the statement of lemma~\ref{lm:dE} in which (eventually straightforward) sense equation (E) is understood. 

\begin{lemma}
  \label{lm:BIO}
  Equations (B), (I) and (O) are satisfied by any vector field $u$ of the form (V).
\end{lemma}

\begin{proof} [Proof of the lemma]
  The boundary condition (B) is satisfied because $\psi$ is constant on each $\partial_j D$.

  The incompressibility condition (I) is satisfied ($v$ being purely coexact):
  \[\delta v = - \delta *d\psi = -*d^2\psi = 0.\]

  Condition (O) is satisfied. Indeed, using Cartan's formula, we have
  \[\oint_{\partial_k D} \left( \partial_t v + L_u v \right) = \partial_t \oint_{\partial_k D} v + \oint_{\partial_k D} (i_u \,dv + d\, i_u v),\]
  where $\oint i_u \, dv = 0$ because of (B), $\oint d\, i_uv=0$ because of Stokes' formula, and, by definition of $g$,
  \[\oint_{\partial_k D} v = -\sum_j \gamma_j \oint_{\partial_k D} *dg_{z_j(t)} = - c_k \sum_j \gamma_j ,\]
  which is constant.
\end{proof}

For any $j=1,...,N$, the vector field $u$ splits into
\begin{equation}
  \label{eq:uj}
  u = u_j + U_j, \quad u_j(\zeta) = - 2\gamma_j \, \Im \left( \frac{\zeta-z_j(t)}{|\zeta-z_j(t)|^2} \, \partial_\zeta \right)
\end{equation}
where $U_j$ is a real analytic vector field in the neighborhood of $z_j(t)$.

\begin{lemma}
  \label{lm:dE}
  Equation (dE) is satisfied iff
  \[\dot z_j = U_j(z_j) \quad (j=1,...,N).\]
\end{lemma}

This statement calls for a clarification, because the term $L_u w = dw \cdot u$ in (dE) is the product of two distributions. In the neighborhood of $z_j(t)$, only the term $dw \cdot u_j$ demands a closer look. In the proof below, we will show that 
\begin{itemize}
\item $u_j$ has a well-defined, continuous flow $\phi_{s}$
\item $\frac{1}{s}(\phi_{s}^*w - w)$ has a null limit as $s \to 0$, so that $dw \cdot u_j = 0$ in the most usual sense.
\end{itemize}

\begin{proof}[Proof of the lemma]
  The first term of the vorticity equation $(dE)$ is
  \[\partial_t w = - \sum \gamma_j d \delta_{z_j(t)} \dot z_j(t).\]
  On the other hand, the advection term is $dw \cdot u$. Away from vortices, both these terms vanish, so (dE) holds. 

  We have assumed that vortices do not collide, so there remains to check the equation in the neighborhood of each separate vortex. Let $j \in \{1,...,N\}$. The Hamiltonian splits into   
  \[\psi = \psi_j + \Psi_j, \quad \psi_j(t,\zeta) = - \gamma_j \log |\zeta-z_j(t)|,\]
  where $\Psi_j$ is real analytic in the neighborhood $\{(t,z_j(t))\}$. The vector field correspondingly splits as in equation~\eqref{eq:uj}. Setting
  \[Z = \zeta - z_j(t)\]
  (where the physical time $t$ is frozen while we investigate the Eulerian flow),
  we see that
  \[u_j = - 2\gamma_j \, \Im \left( \frac{Z}{|Z|^2} \partial_Z \right)\]
  is purely rotational. Its flow
  \[\phi_s(Z) = \exp \left( \frac{i\gamma_j}{|Z|^2} \, s\right) Z\]
  extends by continuity to $\phi_s(0) = 0$. We see that, for any test function $\varphi$,
  \[\phi_s^* \delta_0 \cdot \varphi = \delta_0 \cdot (\varphi \circ \phi_{-s}) = \varphi(\phi_{-s}(0)) = \varphi(0) = \delta_0 \cdot \varphi.\] 
  So,
  $$dw \cdot u_j = \lim_{s \to 0} \frac{1}{s}(\phi_{s}^*w - w) = 0.$$
  Eventually, in the neighborhood of $z_j(t)$, 
  \[\partial_t w + dw \cdot u = \gamma_j d\delta_{z_j(t)} \cdot \left( - \dot z_j(t) + U_j \right).\]
  So (dE) holds on $D$ iff, for all $j$,
  \[\dot z_j = U_j.\]  
\end{proof}

\begin{proof}[End of the proof of the theorem]
  According to lemma~\ref{lm:dE}, (dE) holds on $D$ iff, for each $j$, the vortex $z_j$ moves with the velocity $U_j(z_j)$, where $U_j$ is the Hamiltonian vector field
  \[U_j(\zeta) = 4 \Im \left( \partial_{\bar\zeta} \Psi_j(\zeta) \, \partial_\zeta \right)\]
  of the smooth Hamiltonian
  \[\Psi_j(\zeta) = \psi(\zeta) - \psi_j(\zeta) = \gamma_j\, h(z_j(t),\zeta) + \sum_{k \neq j} \gamma_k\, g(z_k(t),\zeta).\]
  At this stage, each vortex obeys its own equation, with its own vector field $U_j$. The trick is that all these equations are the components of a single Hamiltonian system on $C_N(D)$. Indeed, on the diagonal, $\partial_{\bar z} \bigl( h(z,z) \bigr) = 2\, \partial_{\bar\zeta} h(z,\zeta)|_{\zeta=z}$, by symmetry of $h$; hence, differentiating (H),
  \[\partial_{\bar z_j} H = \gamma_j \left( \gamma_j\, \partial_{\bar\zeta} h(z_j,\zeta) + \sum_{k \neq j} \gamma_k\, \partial_{\bar\zeta} g(z_k,\zeta) \right)_{\zeta = z_j} = \gamma_j\, \partial_{\bar\zeta} \Psi_j(z_j),\]
  so that
  \[\dot z_j = U_j(z_j) = \frac{4}{\gamma_j} \Im \left( \partial_{\bar z_j} H \, \partial_{z_j}\right),\]
  which is precisely the equation $i_{X_H}\omega = dH$ read in the $j$-th factor.
\end{proof}

We will describe the topology of $C_N(D)$ in the simply and doubly connected cases, with $N=2$ vortices (examples~\ref{ex:C2(D)} and~\ref{ex:C2(A)} in section~\ref{sec:first-properties}). See~\cite{Kassel:2008} for more. 

\begin{theorem}[Flucher--Gustafsson]
  \label{thm:Hpsh}
  The opposite $-H$ of the Kirchhoff--Routh Hamiltonian is strictly plurisubharmonic on $C_N(D)$, whatever the intensities $\gamma_1,...,\gamma_N \in \R_*$. Equivalently, every regular level set of $H$ is a strictly pseudoconvex hypersurface of $C_N(D)$, with $\{H>E\}$ on the pseudoconvex side of $\{H=E\}$ ($E \in \R$). 
\end{theorem}

For the sake of completeness, we reproduce the argument of Flucher and Gustafsson (see \cite{Gustafsson:2019}), which relates the complex Hessian of $H$ to a Bergman kernel of $D$. Let $\mathcal{O}_{L^2}(D)$ be the Hilbert space of square integrable holomorphic functions on $D$, $K_0$ its reproducing kernel (the \emph{Bergman kernel} of $D$), and let
\[\widetilde{\mathcal{O}}_{L^2}(D) = \left\{ f \in \mathcal{O}_{L^2}(D), \; \oint_{\partial_k D} f(z)\, dz = 0 \; (\forall k) \right\}\]
be the closed subspace of those $f$ which admit a single-valued primitive; let $K$ be its reproducing kernel, the so-called \emph{reduced} Bergman kernel. Constants belong to $\widetilde{\mathcal{O}}_{L^2}(D)$, so $K(z,z)>0$; $K(z,z)\leq K_0(z,z)$; and $K = K_0$ when $D$ is simply connected. It is $K$ which is attached to the \emph{hydrodynamic} Green function $g$, whereas $K_0$ is attached to the Dirichlet Green function $g_0$ (see section~\ref{sec:1-v}): this makes no difference for the statement, since all we use of $K$ is that it is a reproducing kernel.

What makes the statement noteworthy is that no hypothesis at all is needed on the signs of the intensities: they enter the Hessian through a congruence, which does not see signs. This is in sharp contrast with the compactness properties of section~\ref{sec:first-properties}, which are governed precisely by the relative signs of $\gamma_j$'s.

\begin{proof}
  Let
  \[A_{j\bar{k}} = \frac{\partial^2 H}{\partial z_j \partial \bar{z}_k}\]
  be the coefficients of the complex Hessian matrix of $H$.

  Off the diagonal ($j \neq k$), the derivative acts on the regular part $h$ of the Green function only, the logarithmic term $-\log |z_j-z_k|$ being pluriharmonic away from $\Delta$. By Schiffer's formula (see also lemma~\ref{lm:superharmonic} below), the mixed derivative of $h$ is, up to a factor, the reduced Bergman kernel:
  \[\frac{\partial^2 g(z_j, z_k)}{\partial z_j \partial \bar{z}_k} = \frac{\partial^2 h(z_j, z_k)}{\partial z_j \partial \bar{z}_k}= -\frac{\pi}{2}\, K(z_j, z_k).\]

  On the diagonal ($j=k$), the derivative acts on the Robin function, and
  \[\frac{\partial^2 r}{\partial z \partial\bar z} = \frac14 \Delta r = - \pi K(z,z).\]

  Both contributions have the same shape, so that the whole matrix is
  \[A_{j\bar{k}} = -\frac{\pi}{2} \, \gamma_j \gamma_k\, K(z_j, z_k), \qquad \mbox{i.e.} \qquad A = -\frac{\pi}{2}\, \gamma \, \mathcal{K} \, \gamma,\]
  where $\gamma = \mathrm{diag}(\gamma_1,...,\gamma_N)$ and $\mathcal{K} = \bigl( K(z_j,z_k) \bigr)_{j,k}$. Being the Gram matrix of a reproducing kernel at $N$ pairwise distinct points, $\mathcal{K}$ is Hermitian positive definite. Since $\gamma$ is real, diagonal and invertible, $\gamma \mathcal{K} \gamma$ is congruent to $\mathcal{K}$, hence positive definite as well, regardless of the signs of the $\gamma_j$. So $A$ is negative definite.
\end{proof}

\begin{remark}[Warning on signs] With our convention $g = -\log|\zeta-z|+h$, the Green function tends to $+\infty$ at its pole and the Robin function tends to $-\infty$ at $\partial D$; the Robin function is superharmonic (lemma~\ref{lm:superharmonic}), so already the $1$-vortex Hamiltonian $\frac{\gamma^2}{2}r$ is plurisuperharmonic. It is with the opposite sign convention for $g$, frequent in the literature, that the Kirchhoff--Routh function itself is plurisubharmonic.
\end{remark}

Theorem~\ref{thm:Hpsh} acquires its full strength when combined with the following observation, which we will use in section~\ref{sec:first-properties}.

\begin{lemma}
  \label{lm:Stein}
  The configuration space $C_N(D)$ is a Stein manifold, regardless of the smoothness of $D$.
\end{lemma}

\begin{proof}
  Since $D \subset \C$ is open, it is a domain of holomorphy, hence Stein, and so is the product $D^N$. The configuration space $C_N(D) = D^N \setminus \Delta$ is obtained by removing the fat diagonal, which is the zero locus of the globally defined holomorphic function $\prod_{j<k}(z_j-z_k)$, hence a complex analytic hypersurface; the complement of an analytic hypersurface in a Stein manifold is Stein.
\end{proof}

\section{On the $1$-vortex problem}
\label{sec:1-v}

We will hereafter keep the same hypotheses on the domain $D$ as in section~\ref{sec:introduction}. In this section  we review some properties of the $1$-vortex problem, for which the Hamiltonian is merely the Robin function of $D$. The dynamics is integrable.

The vector field is smooth inside $D$. Let $\delta(z) = \text{dist}(z, \partial D)$ be the distance to the boundary. In lemma~\ref{lm:asymptotics}, key asymptotics of $h$ and $r$ near the boundary are given.

\begin{lemma}
  \label{lm:asymptotics}
  The regular part of the Green function and the Robin function split into
  \[\begin{cases}
    h(z,\zeta) &= \log |\zeta-z^*| + \tilde h(z,\zeta)\\
    r(z) &= \log (2\delta(z)) + \tilde r(z),
    \end{cases}\]
  where $z^*$ is a geometric reflection of $z$ across $\partial D$ and the remainders $\tilde{h}$ and $\tilde r$ are $C^\infty$-smooth across collisions and up to the boundary of $D$. 
\end{lemma}

This kind of asymptotics appears in Flucher's and others' works (see for example \cite{Flucher:1999,Bethuel:2017}).

\begin{proof}
  The difference between the Dirichlet and the hydrodynamic Green functions is smooth up to the boundary, since it consists in a linear combination of harmonic measures. So in order to prove the lemma, we may assume $g = g_0$, $h = h_0$ and $r = r_0$.

  Let $z_0 \in \partial D$. Because $\partial D$ is $C^\infty$, Kellogg's theorem ensures we can find a neighborhood $U$ of $z_0$ and a conformal mapping $f: U \cap  D \to \mathbb{H}$ (the upper half-plane) such that: $f$ extends to a $C^\infty$ diffeomorphism up to the boundary segment $U \cap \partial D$ and $f$ maps $U \cap \partial D$ locally onto the real axis $\R$. Because $f$ is a diffeomorphism up to the boundary, we can define a local $C^\infty$ ``Schwarz reflection'' $z^*$ for any $z \in U \cap  D$ across $\partial D$ by pulling back the standard complex conjugation from $\mathbb{H}$:$$z^* = f^{-1}(\overline{f(z)})$$
  Notice that $z^*$ lies outside $D$ (or on the boundary if $z \in \partial D$), and the map $z \mapsto z^*$ is a smooth involution.

  The standard Green functions for the whole and upper half plane are
  \[g_\C (z,\zeta) = -\log |\zeta-z| \quad \mbox{and} \quad g_{\mathbb{H}}(z,\zeta) = g_\C(z,\zeta) - g_\C(\bar z,\zeta) = \log \left|\frac{\zeta-\bar z}{\zeta- z} \right|.\]
The function $g_{\H}(f(z),f(\zeta))$ is an approximation of $g$ in $U \cap D$ because it has the required logarithmic singularity as $\zeta \to z$, and it vanishes on the boundary ($z = z^*$). So, for $z, \zeta \in U \cap  D$, the function 
  $$\hat g(z, \zeta) = g(z, \zeta) - g_{\mathbb{H}}(f(z), f(\zeta))$$
  is biharmonic and vanishes when $z \in U \cap \partial D$. Thus, by the Schwarz reflection principle, $\hat g(z, \zeta)$ is  $C^\infty$ smooth up to the boundary. Now, substitute the definition of $g_{\mathbb{H}}$ to isolate $h(z, w)$:
  $$g(z, \zeta) = \log \left| \frac{f(\zeta)-\overline{f(z)}}{f(\zeta)-f(z)} \right| + \hat g(z, \zeta).$$
  Since $\overline{f(z)} = f(z^*)$ and 
  $h(z, \zeta) = g(z, \zeta) + \log|\zeta-z|$,
  \[h(z, \zeta) = \log|f(\zeta)-f(z^*)| - \log\left|\frac{f(\zeta)-f(z)}{\zeta-z}\right| + \hat g(z, \zeta).\]
    Factoring out $|\zeta-z^*|$ from the first term, we get
    \[h(z, \zeta) = \log|\zeta-z^*| + \log\left|\frac{f(\zeta)-f(z^*)}{\zeta-z^*}\right| - \log\left|\frac{f(\zeta)-f(z)}{\zeta-z}\right| + \hat g(z, \zeta).\]
    Because $f$ is a smooth diffeomorphism up to the boundary, the second and third terms are smooth. Thus, we obtain the wanted asymptotics for the regular part, with
    \[\tilde h(z,\zeta) = \log\left|\frac{f(\zeta)-f(z^*)}{\zeta-z^*}\right| - \log\left|\frac{f(\zeta)-f(z)}{\zeta-z}\right| + \hat g(z, \zeta).\]

    To find the asymptotic behavior of the Robin function, we now take the limit as $\zeta \to z$:
    $$r(z) = h(z, z) = \log|z-z^*| + \tilde h(z, z).$$
    We must relate the distance to the reflection point $|z-z^*|$ to the geometric distance to the boundary $\delta(z)$. By definition of the reflection, $f(z) - f(z^*) = f(z) - \overline{f(z)} = 2i \operatorname{Im} f(z)$. We can express $|z-z^*|$ as:$$|z-z^*| = \frac{|f(z)-f(z^*)|}{\left| \frac{f(z)-f(z^*)}{z-z^*} \right|} = \frac{2\operatorname{Im} f(z)}{|\sigma(z)|}, \quad \sigma(z) := \frac{f(z)-f(z^*)}{z-z^*},$$ where $\sigma$ is smooth and evaluates to $f'(z)$ on $\partial D$. We now look at the ratio between $\operatorname{Im} f(z)$ and $\delta(z)$. Both are smooth equations for $\partial D$ (they vanish on the boundary and have positive inward normal derivatives). By Taylor's theorem, their ratio $\frac{\operatorname{Im} f(z)}{\delta(z)}$ is smooth, positive, and on the boundary it evaluates to $|f'(z)|$. Therefore, we can rewrite $|z-z^*|$ exactly as
    $$|z-z^*| = 2\delta(z) \, \varpi(z), \qquad \varpi(z) := \frac{\operatorname{Im} f(z)}{\delta(z)\, |\sigma(z)|},$$ where $\varpi$ is a strictly positive $C^\infty$ function up to the boundary and, crucially, $\varpi = 1$ on $\partial D$. Taking logarithms,$$\log|z-z^*| = \log(2\delta(z)) + \log \varpi(z).$$
    Substituting this back into the expression for $r(z)$:
    $$r(z) = \log(2\delta(z)) + \log \varpi(z) + \tilde h(z, z)$$
    The function
    \[\tilde r(z) = \log \varpi(z) + \tilde h(z, z)\]
    is $C^\infty$ up to $\partial D$. For $g = g_0$ one moreover has $\tilde h(z,z) = \log|f'(z)|-\log|f'(z)|+0 = 0$ on the boundary, and $\log\varpi=0$ there, so that $\tilde r$ vanishes on $\partial D$ and $r_0(z)-\log(2\delta(z)) = O(\delta(z))$. For the hydrodynamic Green function one has to add back the linear combination of harmonic measures, and $\tilde r$ tends to the constant $\ra_j$ of lemma~\ref{lem:C} on $\partial_jD$ instead of to $0$.
\end{proof}

\begin{remark}
  The proof above extends when the boundary is $C^{2,\alpha}$. More precisely, if it is $C^{k+1,\alpha}$ with $k \geq 1$, the reflection map is $C^{k,\alpha}$, and the remainder function $\tilde r$ has the same regularity $C^{k,\alpha}$. This can also be seen directly in the simply or doubly connected cases. Indeed, by Kellogg's Theorem, if $\partial D \in C^{k+1,\alpha}$, the conformal map is $C^{k+1,\alpha}(\bar{D})$. But the Robin function depends on the conformal radius, which requires evaluating the derivative of the conformal map. Thus, regularity drops to $C^{k, \alpha}$.
\end{remark}

Of higher interest is the following high regularity result.

\begin{lemma}
  If the boundary is real analytic, both $\tilde h$ and $\tilde r$ are real analytic too up to the boundary.
\end{lemma}

\begin{proof}
  Because the boundary is a real-analytic curve, the Schwarz reflection principle guarantees that the local conformal map $f: U \cap D \to \mathbb{H}$ can be holomorphically, hence real analytically, continued across $\partial D$. Consequently, the local reflection $z^* = f^{-1}(\overline{f(z)})$ too is real-analytic in a neighborhood of the boundary.

  Since $f$ is holomorphic and its derivative does not vanish on the boundary, the difference quotient $Q(z, \zeta) = \frac{f(\zeta)-f(z)}{\zeta-z}$ is a non-vanishing, holomorphic function of both $z$ and $\zeta$ in a neighborhood of the boundary. Taking the modulus and then the logarithm of a non-vanishing holomorphic function yields a real-analytic function. The same applies to $\frac{f(\zeta)-f(z^*)}{\zeta-z^*}$.

  The function $\hat g(\cdot,\zeta)$ is harmonic in $D$ and vanishes identically on $\partial D$. By the Schwarz reflection principle for harmonic functions, it extends as a harmonic --- therefore real-analytic --- function across the real-analytic boundary $\partial D$.

  For any real-analytic curve, the distance function $\delta(z)$ is  real-analytic in a sufficiently small tubular neighborhood of the curve. Because both $\operatorname{Im} f(z)$ and $\delta(z)$ are real-analytic and vanish to exactly first order on $\partial D$, their ratio $\frac{\operatorname{Im} f(z)}{\delta(z)}$ is uniquely extended as a strictly positive, real-analytic function up to and across the boundary. Because $\tilde h(z, \zeta)$ and $\tilde r(z)$ are constructed entirely from additions, subtractions, and logarithms of these non-vanishing, real-analytic components, they inherit the exact same regularity. Thus, they are real-analytic up to the boundary.
\end{proof}

\begin{proposition}
  \label{prop:Robin}
  The Robin function $r$ has the following properties:
  \begin{enumerate}
  \item Its critical level sets form a set of zero Lebesgue measure in $D$. 
  \item Each component $\partial_k D$ has a tubular neighborhood $\mathcal{V}(\partial_k D)$ with no critical point.
  \item In $\mathcal{V}(\partial_k D)$, level sets of $r$ are homotopic to $\partial_k D$.
  \end{enumerate}
\end{proposition}

\begin{proof}
  \begin{enumerate}
  \item The first item is a classical property of non constant real analytic functions, of which we recall the proof. Let $C = \{z \in D, \; dr_z = 0\}$ be the critical set of $r$. It follows from Łojasiewicz's theorem on the stratification and the triangulation of semi-analytic sets (see~\cite{Benedetti:1990}) that $C$ consists of isolated points and analytic curves, with locally finitely many connected components. Since $r$ is locally constant on these components, $r(C \cap K)$ must be finite for any compact subset $K$ of $D$. Since $D$ is open, it is $\sigma$-compact and can be written as the union of countably many compact subsets. Therefore, the set $r(C)$ of critical values of $r$ is (at most) countable.

    The principle of analytic continuation shows that $r$ cannot be constant on an open set of $D$. So, if $v$ is a critical value of $r$, $r^{-1}(v)$ must have empty interior, its Hausdorff dimension is $\leq 1$, and its Lebesgue measure is null. The first statement follows. 
  \item From lemma~\ref{lm:asymptotics},
    \[d r(z) = \frac{1}{\delta(z)}\, d \delta(z) + d\tilde r(z).\]
    Since $d\delta(z)$ and $d\tilde{r}$ are $O(1)$, $d r(z)$ cannot vanish if $\delta(z)$ is small. Since $\partial_k D$ is compact, there exists $\epsilon>0$ such that the tubular neighborhood
    \[\mathcal{V}(\partial_k D) = \{z \in D, \, \mathrm{dist}(z,\partial_k D) < \epsilon\}\]
    has no critical point. 
  \item The last claim is a consequence of Morse theory applied to the harmonic radius 
    \[q(z) = e^{r(z)} = 2\delta(z) \, \exp \left( \tilde r(z) \right)\]
    in $\mathcal{V}(\partial_k D)$.
  \end{enumerate}
\end{proof}

Considering the Hamiltonian vector field of $r$, lemma~\ref{lm:asymptotics} and proposition~\ref{prop:Robin} entail the following dynamical result.

\begin{corollary}
  \label{cor:one}
  Periodic orbits for one vortex in $D$ form a dense subset $D_{per}$ of $D$ of full Lebesgue measure, containing the neighborhood $\sqcup_k \mathcal{V}(\partial_k D)$ of the boundary. Close to $\partial D$, the period tends to $0$. 
\end{corollary}

The following classical lemma and remark aim at better understanding the homotopy classes of the level curves or $r$.

\begin{lemma}
  \label{lm:superharmonic}
  The Robin function $r$ is superharmonic.
\end{lemma}

We refer to~\cite{Ahlfors:1979} (and the Schiffer formula therein) or \cite{Blocki:2013} for details. 

\begin{proof}[Idea of proof]
  Using that the regular part of the Green function, $h$, is harmonic in each variable, and that $r(z) = h(z,z)$, the two pure terms drop out of the chain rule and we are left with the mixed ones, which are equal by symmetry of $h$:
  \[\Delta r = 4 \frac{\partial^2 }{\partial z \partial\bar z} \bigl( h(z,z) \bigr) = 8 \left. \frac{\partial^2 h}{\partial \zeta\partial\bar z} \right|_{\zeta=z}.\]
  The reduced Bergman kernel (as introduced above) may be defined by Schiffer's formula:
  \[K(z,\zeta) = - \frac{2}{\pi} \frac{\partial^2 h}{\partial z \partial \bar \zeta}(z,\zeta).\]
  So
  \begin{equation}
    \label{eq:laplacian-r}
    \Delta r = -4\pi K(z,z).
  \end{equation}
  Because of the Bergman equality
  \[K(z,z) = \sup_{f \in \widetilde{\mathcal{O}}_{L^2}(D), \, \|f\|_{L^2(D)} \leq 1} |f(z)|^2 >0\]
  (since $D$ is bounded, the supremum is taken over a set of functions containing non-zero constants, which do have a single-valued primitive),
  \[\Delta r <0.\]
\end{proof}

\begin{remark}
  Assume that the Robin function is Morse, which is true for a generic domain~\cite{Bartsch:2018}. Proposition~\ref{prop:Robin} shows that close to $\partial D$ (where $r \to - \infty$) components of level sets $\Lambda_{s} = \{r = s\}$ are homotopic to components of $\partial D$. As the value of $s$ increases, the homotopy class of $\Lambda_{s}$ does not change, unless $s$ crosses a critical value of $r$. In dimension $2$, a priori there are three possible cases:
  \begin{itemize}
  \item If we passed through a local minimum (null index), in the neighborhood of the corresponding critical point, a new point and then a new closed curve would appear as a connected component of $\Lambda_s$. But because of lemma~\ref{lm:superharmonic}, $r$ does not attain any local minimum inside $D$. % (The Laplacian is the sum of the eigenvalues of the Hessian, so the two eigenvalues cannot both be positive.)
  \item When we pass through a saddle point (index $1$), two connected components get close to each other and recombine.
  \item When we pass through a local maximum (index $2$), a connected component of $\Lambda_s$ shrinks and disappears. 
  \end{itemize}
  All connected components of level sets of $r$ are circles. They cannot appear ex nihilo inside $D$ (no critical point of $0$ index). So, any component of a regular level curve $\Lambda_s$ necessarily comes from an already existing one corresponding to smaller values of $s$.

  Recall that $m$ is the connectivity degree of $D$. When $m=1$ or $2$, any connected component of a level set of $r$ must be homotopic to a component of $\partial D$. When $m=3$ (pair of pants): if, for example, the level curves arising from $\partial_2D$ and $\partial_3D$ merge, they transform into a single circle enclosing both $\partial_2D$ and $\partial_3D$, which is therefore homotopic to the outer component $\partial_1D$. When $m \geq 4$, a new phenomenon occurs. Indeed, if two circles, say around $\partial_2D$ and $\partial_3D$, merge first, the merged circle encloses $\partial_2D$ and $\partial_3D$ while leaving $\partial_1D$ and $\partial_4D$ outside, and it is homotopic to no $\partial_kD$. So, for $m \geq 4$, some level sets of the Robin function are not homotopic to any connected component of $\partial D$.
\end{remark}

\medskip

We complete this section by relating the one-vortex flow to Suita's conjecture. Suita compared two ways of measuring, in a conformally invariant manner, how far a point $z\in D$ lies from $\partial D$: by potential theory, through the logarithmic capacity of $\partial D$ seen from $z$, and by complex analysis, through the largest value at $z$ of a square integrable holomorphic function of norm $1$. Suitably normalized, both quantities coincide on simply connected domains, and the conjecture asserts that, relatively to the capacity, the Bergman kernel is never smaller than in the simply connected case.

The conjecture concerns the Dirichlet normalization. So let $D$ be of arbitrary connectivity $m$, and consider the Dirichlet Green function $g_0$, its Robin function $r_0$ and its harmonic radius $q_0 = e^{r_0}$ (definition~\ref{def:Robin}), as well as the Bergman kernel $K_0$ of $D$ (section~\ref{sec:introduction}). The logarithmic capacity of $\partial D$ with respect to $z$ is
\[
  c_\beta(z) = \exp\Bigl(-\lim_{\zeta\to z}\bigl(g_0(z,\zeta)+\log|\zeta-z|\bigr)\Bigr) = e^{-r_0(z)} = \frac{1}{q_0(z)},
\]
and the metric $c_\beta(z)\,|dz| = |dz|/q_0(z)$ is called the \emph{capacity metric} of $D$; it is conformally invariant, by corollary~\ref{cor:conformal-r} applied to $g_0$. On the other hand,
\[K_0(z,z) = \sup\bigl\{|f(z)|^2 : f\in\mathcal O_{L^2}(D),\ \|f\|_{L^2(D)}\le1\bigr\},\]
and $K_0(z,z)^{1/2}\,|dz|$ is conformally invariant too. On the unit disc, $q_0(z) = 1-|z|^2$ (example~\ref{ex:disc}) and $\pi K_0(z,z) = (1-|z|^2)^{-2}$, so that $\pi K_0 = c_\beta^2$; by conformal invariance, this identity holds on every simply connected domain.

\begin{theorem}[Suita's conjecture~\cite{Suita:1972}; B\l{}ocki~\cite{Blocki:2013}, Guan--Zhou~\cite{Guan:2015}]
  \label{thm:suita}
  For every $z\in D$,
  \[\pi K_0(z,z) \geq c_\beta(z)^2 = \frac{1}{q_0(z)^2},\]
  with equality at some point of $D$ if and only if $m=1$.
\end{theorem}

B\l{}ocki proved the inequality for bounded domains of $\C$, and Guan and Zhou proved it on open Riemann surfaces, with equality at one point only if the surface is conformally equivalent to the disc minus a closed polar set, which for our domains means $m=1$; both proofs rely on the Ohsawa--Takegoshi $L^2$ extension theorem with an optimal constant. The inequality has a geometric meaning. Schiffer's formula, in its original form for the Dirichlet Green function, reads $K_0(z,\zeta) = -\frac{2}{\pi}\frac{\partial^2 h_0}{\partial z\,\partial\bar\zeta}(z,\zeta)$, so that, as in the proof of lemma~\ref{lm:superharmonic}, $\Delta r_0 = -4\pi K_0(z,z)$. Writing the capacity metric as $e^{-2r_0}|dz|^2$, its Gaussian curvature is
\[
  \kappa_0 = e^{2r_0}\,\Delta r_0 = -4\pi\, q_0^2\, K_0(z,z),
\]
and theorem~\ref{thm:suita} states that $\kappa_0\leq-4$, which is the value of the curvature on the disc, where the capacity metric is the Poincar\'e metric $|dz|/(1-|z|^2)$. This is the form in which Suita stated his conjecture.

\medskip

The vortex, however, does not see the Dirichlet normalization: its Hamiltonian is the hydrodynamic Robin function $r$, and $\Delta r = -4\pi K(z,z)$ involves the reduced kernel $K$ (lemma~\ref{lm:superharmonic}). The Hamiltonian vector field of $r$ for $\omega_0 = dx\wedge dy$ is $X_r = \partial_y r\,\partial_x - \partial_x r\,\partial_y$, whose curl is
\[\operatorname{curl} X_r = -\Delta r = 4\pi\,K(z,z) > 0 .\]
The curl depends neither on the additive normalization of $r$ nor, when $m\ge2$, on the circulations $c_1,\dots,c_m$. Indeed, two hydrodynamic Green functions of $D$ differ by a function of the separated form $\sum_k a_k\bigl(\omega_k(z)+\omega_k(\zeta)\bigr)+b$, $\omega_k$ being the harmonic measure of $\partial_kD$, which is harmonic on the diagonal.

When $D$ is simply connected, the two normalizations agree: $g=g_0$, hence $r=r_0$, $q=q_0$ and $K=K_0$. Theorem~\ref{thm:suita} is then an identity, and it expresses the curl of the one-vortex flow in terms of its Hamiltonian.

\begin{corollary}\label{cor:suita}
  Let $D$ be simply connected, and let $\phi : D\to\mathbb D$ be a
  Riemann map. Then $r = \log\bigl(1-|\phi|^2\bigr) - \log|\phi'|$, and
  \[
    \operatorname{curl} X_r(z) = 4\pi K(z,z) = \frac{4}{q(z)^2}
    = \frac{4\,|\phi'(z)|^2}{\bigl(1-|\phi(z)|^2\bigr)^2}
    \qquad (z\in D).
  \]
  Equivalently:
  \begin{enumerate}
    \item the Robin function solves Liouville's equation
      $\Delta r + 4e^{-2r} = 0$;
    \item the metric $(\operatorname{curl} X_r)\,|dz|^2$ is the Poincaré metric of $D$, of curvature $-1$;
    \item $\log (\operatorname{curl} X_r) + 2r \equiv \log 4$, so that $X_{\log \operatorname{curl} X_r} = -2X_r$. In particular the curl is a first integral of the one-vortex flow, and the vortex moves along the level curves of the hyperbolic density of $D$.
  \end{enumerate}
\end{corollary}

\begin{proof}
  Using example~\ref{ex:disc} and conformal invariance of the Green function, we see that $q = \bigl(1-|\phi|^2\bigr)/|\phi'|$. The Bergman kernel $K=K_0$ transforms as
  \[K(z,\zeta) = \frac{\phi'(z)\,\overline{\phi'(\zeta)}} {\pi\bigl(1-\phi(z)\overline{\phi(\zeta)}\bigr)^2}, \] which gives the displayed identity on the diagonal. Then $\Delta r = - \operatorname{curl} X_r = -4e^{-2r}$, which is~(1). Next, $\Delta\log \operatorname{curl} X_r = -2\Delta r = 2 \operatorname{curl} X_r$, and $(\operatorname{curl} X_r)\,|dz|^2$ is the pull-back by $\phi$ of $4|dw|^2/(1-|w|^2)^2$, which is~(2). Finally~(3) is a restatement of $\operatorname{curl} X_r = 4e^{-2r}$.
\end{proof}

\begin{remark}\label{rem:suita-m}
  When $m\geq2$, theorem~\ref{thm:suita} says nothing about the one-vortex flow. It bounds $\operatorname{curl}X_{r_0} = 4\pi K_0$ from below, but $r_0$ is the Hamiltonian of no vortex flow, since $g_0$ is not a hydrodynamic Green function; and the curl of the actual flow is smaller, $K(z,z)\le K_0(z,z)$, the supremum defining $K(z,z)$ being taken over a subspace of codimension $m-1$. Nor can $q_0$ be replaced by the hydrodynamic harmonic radius $q$: changing the circulations adds to $r$ the non-constant harmonic function $2\sum_k a_k\omega_k$, so that $q^2\operatorname{curl}X_r$ is not even determined by $D$. Proposition~\ref{prop:converse-suita} below shows that the identity of corollary~\ref{cor:suita}, written in terms of the flow only, characterizes simply connected domains.
\end{remark}

\medskip

Conversely, the local dynamics of one vortex detects whether $D$ is simply connected.

\begin{proposition}\label{prop:converse-suita}
  Let $D$ be as in section~\ref{sec:introduction}, of connectivity $m$, and let $u = \operatorname{curl} X_r$. The following conditions are equivalent:
  \begin{enumerate}
    \item[(i)] $m = 1$;
    \item[(ii)] $\Delta\log u = 2u$ on some non-empty open subset of $D$,
      i.e.\ the metric $u\,|dz|^2$ has curvature $-1$ there;
    \item[(iii)] $\log u + 2r$ is constant on some non-empty open subset
      of $D$.
  \end{enumerate}
\end{proposition}

Since $u$ is determined by the vector field $X_r$, and condition~(iii)
only involves $dr = \iota_{X_r}\omega_0$, either condition is a
property of the germ of the one-vortex flow at an arbitrary point of
$D$.

\begin{proof}
  (i)~$\Rightarrow$~(ii), (iii) is corollary~\ref{cor:suita}, and
  (iii)~$\Rightarrow$~(ii) follows from $\Delta\log u = -2\Delta r = 2u$.
  We prove (ii)~$\Rightarrow$~(i).

  \emph{Analyticity.} The reduced Bergman kernel
  $K(z,\zeta)$ is holomorphic in $(z,\bar\zeta)$ on
  $D\times D$, so $u$ is real analytic and positive on $D$. The
  function $\Delta\log u - 2u$ is real analytic and vanishes on an
  open set, hence on the connected set $D$.

  \emph{Completeness.} By lemma~\ref{lm:asymptotics},
  $r = \log(2\delta) + \tilde r$ with $\tilde r$ smooth up to
  $\partial D$. Near $\partial D$ the distance $\delta$ is smooth with
  $|\nabla\delta| = 1$, so
  $\Delta\log\delta = \Delta\delta/\delta - 1/\delta^2$ and
  $u = \delta^{-2}\bigl(1+O(\delta)\bigr)$. Hence
  $\sqrt u\ge 1/(2\delta)$ on a collar of $\partial D$. Since $\delta$
  is $1$-Lipschitz, every path tending to $\partial D$ has infinite
  $\sqrt u\,|dz|$-length, and the metric $u\,|dz|^2$ is complete.

  \emph{Lift to the universal cover.} Let
  $\pi:\mathbb D\to D$ be a universal covering ($D$ is bounded, hence
  hyperbolic), and $\tilde u := (u\circ\pi)\,|\pi'|^2$, so that
  $\tilde u\,|dw|^2 = \pi^*\bigl(u\,|dz|^2\bigr)$ is complete, of
  curvature $-1$, on the simply connected disc. By Liouville's
  representation theorem,
  $\tilde u = 4|f'|^2/(1-|f|^2)^2$ for some locally injective
  holomorphic $f:\mathbb D\to\mathbb D$ \cite{Liouville:1853,Kraus:2013}.
  Then $f$ is a local isometry from a complete surface to the Poincar\'e disc,
  hence a Riemannian covering, hence a biholomorphism since $\mathbb D$ is simply
  connected. Therefore $f\in\operatorname{Aut}(\mathbb D)$ and
  $\tilde u(w) = 4/(1-|w|^2)^2$.

  \emph{Polarization.} The function
  \[
    G(w,v) := 4\pi\,K\bigl(\pi(w),\pi(v)\bigr)\,
    \pi'(w)\,\overline{\pi'(v)}
  \]
  is holomorphic in $(w,\bar v)$ on $\mathbb D\times\mathbb D$, and so
  is $4/(1-w\bar v)^2$. Both functions agree on the diagonal, where
  they equal $\tilde u(w)$. A function holomorphic in $(w,\bar v)$
  which vanishes for $v=w$ vanishes identically, so
  $G(w,v) = 4/(1-w\bar v)^2$ on $\mathbb D\times\mathbb D$.
  It is here that we use that $K$ is single-valued on
  $D\times D$.

  \emph{Deck transformations.} Let $\gamma$ be a deck transformation:
  $\pi\circ\gamma = \pi$, hence
  $\pi'\bigl(\gamma(w)\bigr)\,\gamma'(w) = \pi'(w)$. Then
  \[
    G\bigl(w,\gamma(w)\bigr)
    = \frac{\tilde u(w)}{\overline{\gamma'(w)}}
    = \frac{4}{(1-|w|^2)^2\,\overline{\gamma'(w)}} ,
  \]
  while the formula for $G$ gives
  $G\bigl(w,\gamma(w)\bigr) = 4/\bigl(1-w\,\overline{\gamma(w)}\bigr)^2$.
  At $w=0$ we get $\gamma'(0) = 1$. Writing
  $\gamma(w) = e^{i\theta}(w-a)/(1-\bar a w)$, we have
  $\gamma'(0) = e^{i\theta}(1-|a|^2)$, so $a=0$, $\theta = 0$ and
  $\gamma = \mathrm{id}$. The deck group is trivial, $\pi$ is a
  biholomorphism, and $D$ is simply connected.
\end{proof}

\begin{remark}[Two vortices]\label{rem:converse-two}
  The same information can be read from the two-vortex problem near
  any configuration. By theorem~\ref{thm:Hpsh}, the complex Hessian of
  $H$ is $A_{j\bar k} = -\frac\pi2\gamma_j\gamma_k\,
  K(z_j,z_k)$. The coupling coefficient
  \[
    c(z_1,z_2) := \frac{|A_{1\bar2}|^2}{A_{1\bar1}A_{2\bar2}}
    = \frac{|K(z_1,z_2)|^2}
           {K(z_1,z_1)\,K(z_2,z_2)}
  \]
  depends neither on $\gamma_1,\gamma_2$, nor on the additive constant
  of $H$, nor on the circulations, and it is conformally invariant.
  Let $d_u$ be the distance of the metric $u\,|dz|^2$. If $m=1$, then
  $c = \cosh^{-4}(d_u/2)$ (compute at $z_1=0$ on the disc, where
  $c = (1-|z_2|^2)^2$ and $|z_2| = \tanh(d_u/2)$).

  Conversely, as $z_2\to z_1$,
  \begin{align*}
    \log c &= -\partial_z\partial_{\bar z}\log K(z_1,z_1)\,
      |z_2-z_1|^2 + O(|z_2-z_1|^3),\\
    -4\log\cosh\frac{d_u}2 &= -2\pi K(z_1,z_1)\,|z_2-z_1|^2
      + O(|z_2-z_1|^3).
  \end{align*}
  So the identity $c = \cosh^{-4}(d_u/2)$, already at second order
  near the diagonal of an open set of configurations, is equivalent to
  condition~(ii) of proposition~\ref{prop:converse-suita}, and hence
  to $m=1$.
\end{remark}

We will see in section~\ref{sec:quadruple} how some of these considerations generalize to two (or more) same-sign vortices.

\section{First properties of the $2$-vortex problem}
\label{sec:first-properties}

In the rest of the paper, we consider the case of $N=2$ vortices, in a domain $D$ as before. Its  Hamiltonian
\begin{equation}
  \label{eq:H2}
  H = \gamma_1\gamma_2 \left( -\log |z_2-z_1| + h(z_1,z_2) \right) + \frac{1}{2} \left( \gamma_1^2 r(z_1) + \gamma_2^2 r(z_2) \right),
\end{equation}
is defined on the configuration space $C_2(D) = D^2 \setminus \Delta$, with symplectic form
\[\omega = \frac{1}{2i} \left( \gamma_1 d \bar z_1 \wedge dz_1 + \gamma_2 d\bar z_2 \wedge dz_2 \right).\]

\begin{example}[Configuration space on the disc]
  \label{ex:C2(D)}
  Suppose $D$ is the unit disc $\D$, thus diffeomorphic to $\C$. The diffeomorphism
  \[C_2(\C) \to \C \times \C^\times, \quad (z_1,z_2) \mapsto \left( \frac{z_1+z_2}{2}, z_1-z_2 \right)\]
  shows that $C_2(\D)$ is diffeomorphic to $\R^3 \times \SS^1$. 
\end{example}

\begin{example}[Configuration space on the annulus]
  \label{ex:C2(A)}
  Suppose $D$ is the annulus, thus diffeomorphic to $\C^\times$. The map
  \[C_2(\C^\times) \to \C^\times \times (\C \setminus \{0,1\}),\quad (z_1, z_2) \mapsto \left( z_1, \frac{z_2}{z_1} \right),\]
  defined by scaling the second point by the inverse of the first, is a smooth, invertible bijection, making it a global diffeomorphism. Topologically, this implies the configuration space of the annulus is diffeomorphic to $(S^1 \times \mathbb{R}) \times P$, where $P$ is the interior of a pair of pants (a $2$-sphere with three punctures).
\end{example}

The space $C_2(D)$ is non-compact for ``two'' reasons: its boundary is the union of
\begin{itemize}
\item the boundary regime
  \[(\partial D \times \bar D) \cup (\bar D \times \partial D),\]
  where at least one vortex lies on the domain boundary, and
\item the collision regime
  \[\Delta_{\bar D} = \{(z,z), \; z \in \bar D\},\]
  representing binary collisions.
\end{itemize}
The two regimes overlap at ``quadruple collisions'' of the two vortices and of their respective images across $\partial D$, $\{(z,z), \; z = z^*\}$ (see lemma~\ref{lm:asymptotics} for a definition of images).

We write
\[\Sig_E = H^{-1}(E) \subset C_2(D)\]
for the energy level of energy $E \in \R$; this notation will be used throughout the paper. The first and most consequential difference between the same-sign and the opposite-sign problems is the following.

\begin{theorem}[Compactness dichotomy]
  \label{thm:compact}
  Level sets of $H$ are
  \begin{itemize}
  \item non-compact in the same-sign problem, their closure meeting $\partial C_2(D)$ exactly at quadruple collisions;
  \item compact in the opposite-sign problem.
  \end{itemize}
\end{theorem}

So the flow is complete in the opposite-sign problem. We will see in section~\ref{sec:quadruple} that it is complete in the same-sign problem too, although for a much less immediate reason.

\begin{proof}
  The boundary of $C_2(D)$ consists of the configurations where one or both vortices lie on $\partial D$, and of those where the two vortices collide (possibly on $\partial D$). From the expression~\eqref{eq:H2} of the Hamiltonian, the logarithmic term $-\gamma_1\gamma_2 \log|z_2-z_1|$ blows up at a collision, and the only terms that can compensate for it are the remaining three; each of them is bounded as long as both vortices stay away from $\partial D$. So the closure of an energy level set can contain a collision only if both vortices approach the boundary at the same time. Whether this does happen is precisely what the sign of $\gamma_1\gamma_2$ decides. We inject the asymptotics of lemma~\ref{lm:asymptotics}.

  First assume that $\gamma_1\gamma_2>0$, and, for example, that both intensities are $>0$. Let $E \in \R$. We exhibit a path in $\Sig_E$ tending to $(\partial D)^2 \cap \Delta$.

  Locally, approximate the boundary $\partial D$ by the real axis and the domain $D$ by the upper half-plane. For points $z_j = x_j + iy_j$, the distance to the boundary is $y_j$.  The regular part of the Green function and the Robin function satisfy
  \[h(z_1, z_2) \sim \log|z_1 - \bar{z}_2| \quad \mbox{and} \quad 
  r(z) \sim \log|z - \bar{z}| = \log(2y).\]
  Let's assume the vortices approach the boundary at the same height ($y_1 = y_2 = y \to 0$) while colliding horizontally with a separation $\delta = |x_1 - x_2| \to 0$. The distance between $z_1$ and the "image" of $z_2$ is $|z_1 - \bar{z}_2| = \sqrt{\delta^2 + (2y)^2}$. Substituting these into the Hamiltonian gives
  $$H \sim  -\gamma_1 \gamma_2 \log \delta + \gamma_1 \gamma_2 \log\sqrt{\delta^2 + 4y^2} + \frac{1}{2}(\gamma_1^2 + \gamma_2^2)\log(2y).$$
  The first term goes to $+\infty$ while the second and third go to $-\infty$. 
  Consider a regime where the vortices collide faster than they approach the boundary: $\delta \ll y$. The distance to the image point is dominated by $y$:
  $$\log\sqrt{\delta^2 + 4y^2} \sim \log(2y)$$
  Substituting this approximation back into the Hamiltonian pools all the $-\infty$ terms together:
  $$H \sim -\gamma_1 \gamma_2 \log \delta + \gamma_1 \gamma_2 \log(2y) + \frac{1}{2}(\gamma_1^2 + \gamma_2^2)\log(2y).$$
  So we may factorize the logarithmic term:
  $$H \sim -\gamma_1 \gamma_2 \log \delta + \frac{1}{2}(\gamma_1 + \gamma_2)^2 \log(2y)$$
  To maintain a constant energy $H \sim E$, we balance the two terms:
  $$\gamma_1 \gamma_2 \log \delta \sim \frac{1}{2}(\gamma_1 + \gamma_2)^2 \log(2y).$$
  Solving for the separation $\delta$ as a function of the boundary distance $y$ gives:
  \[\delta \sim y^\kappa, \quad \mbox{with} \quad \kappa = \frac{(\gamma_1 + \gamma_2)^2}{2\gamma_1 \gamma_2} \geq 2.\]
  Since $\kappa > 1$, our initial assumption that $\delta \ll y$ is validated. Using the implicit function theorem, we see that there is a path with the same asymptotics, drawn on $\Sig_E$, which is therefore non-compact. The exponent $\kappa$ is not an accident: it is the quasi-homogeneous weight which governs the whole of section~\ref{sec:quadruple}, and the computation above is the shadow of the exact cancellation of proposition~\ref{prop:cancel}.

  \medskip Second, assume that $\gamma_1\gamma_2<0$. We have just seen that the closure of an energy level set $\Sig_E$ meets neither $\partial D \times D$, nor $D \times \partial D$, nor the part of $\Delta$ lying over $D$. The only remaining place where the closure of $\Sig_E$ could meet $\partial C_2(D)$ is $(\partial D)^2 \cap \Delta$, i.e. a collision of the two vortices on the boundary.

Again using the  half-plane local model near the boundary diagonal, we see that in this neighborhood
  $$H \sim |\gamma_1 \gamma_2| \log \delta - |\gamma_1 \gamma_2| \log|z_1 - \bar{z}_2| + \frac{1}{2}(\gamma_1^2 + \gamma_2^2)\log(2y);$$
  the first and third terms tend to $-\infty$, while the middle one tends to $+\infty$. Group the regular part $h(z_1, z_2)$ back together with the singular collision part to reconstruct the full Green function
  $$\gamma_1\gamma_2\, g(z_1,z_2) \sim |\gamma_1 \gamma_2| \log \delta - |\gamma_1 \gamma_2| \log|z_1 - \bar{z}_2| = -|\gamma_1 \gamma_2| \log \left( \frac{|z_1 - \bar{z}_2|}{|z_1 - z_2|} \right).$$
  Since 
  $$\frac{|z_1 - \bar{z}_2|}{|z_1 - z_2|} \ge 1,$$
  asymptotically,
  $$H \sim \underbrace{-|\gamma_1 \gamma_2| \log \left( \frac{\sqrt{\delta^2 + 4y^2}}{\delta} \right)}_{\le 0} + \underbrace{\frac{1}{2}(\gamma_1^2 + \gamma_2^2)\log(2y)}_{\to -\infty \text{ as } y \to 0}.$$
This seals the proof: all infinities and boundaries of $C_2(D)$ map to $-\infty$, ensuring that $\Sig_E$ is strictly bounded away from the boundaries, and thus $\Sig_E$ is compact.
\end{proof}

\begin{corollary}
  \label{cor:apriori}
  In the opposite-sign problem, for every energy $E$ there exists $\varepsilon_E > 0$, tending to $0$ exponentially as $E \to -\infty$ (like $e^{\varsigma_0 E}$ with $\varsigma_0 = \max\{ 1/|\gamma_1\gamma_2|, 2/\gamma_1^2, 2/\gamma_2^2 \}$), such that, everywhere on $\Sig_E$,
  \[|z_1-z_2| \geq \varepsilon_E \quad \mbox{and} \quad \dist(z_j,\partial D) \geq \varepsilon_E \quad (j=1, 2).\]
\end{corollary}

\begin{proof}
  The space $C_2(D)$ sits inside the compact manifold with corners $\bar{D} \times \bar{D}$. The boundaries of $C_2(D)$ in this ambient space form a closed set:
  $$\partial C_2(D) = \Delta \cup (\partial D \times \bar{D}) \cup (\bar{D} \times \partial D).$$
  Because $\Sig_E$ is a compact set and is completely disjoint from the closed boundary $\partial C_2(D)$, the continuous distance function from $\Sig_E$ to $\partial C_2(D)$ must attain a strictly positive minimum. This minimum dictates that there is a strict lower bound on both the inter-vortex separation $|z_1 - z_2|$ and the boundary distances $\text{dist}(z_j, \partial D)$ across the entire level set.
\end{proof}

\begin{proposition}
  \label{prop:critcompact}
  Equilibrium points of $X_H$ belong to a compact subset of $C_2(D)$. In particular the set of critical values of $H$ is contained in a compact interval $[E_-,E_+]$.
\end{proposition}

\begin{remark}
  In the same-sign problem, we have seen that energy level sets are non-compact; but, thanks to this proposition, the Palais--Smale condition holds for $H$. Compactness of the critical set is what guarantees the existence of critical points of $H$ via Ljusternik--Schnirelmann or Morse theory, the topological changes in the level sets having to correspond to actual critical points rather than to critical points escaping to the boundary of $C_2(D)$. We will not analyse the topology of level sets further than in the two examples above.
\end{remark}

\begin{proof}[Proof of the proposition]
  -- Domain boundary limit (while the two vortices remain at some finite distance): In the differential of $H$, the terms coming from the Robin function go to infinity, and cannot be balanced. 

  -- Binary collision limit ($z_1=z_2$ away from the domain boundary): The dominating term comes from $-\gamma_1\gamma_2 \log |z_2-z_1|$, which has no critical point.

  -- Quadruple collision limit ($z_1=z_2=z_1^*=z_2^*$): For the sake of clarity, we will use a blow-up analysis, which will provide us with a local model for $H$. If a sequence of critical points were to approach a quadruple collision on the boundary of $D$, we could zoom in on the collision point. Because the boundary is smooth, it locally flattens out in the microscopic limit, and the domain converges to the upper half-plane $\H = \{z \in \mathbb{C} : \Im(z) > 0\}$. The critical points would then rescale to critical points of the 2-vortex Hamiltonian on the upper half-plane, $H_{\H}$. However, a scaling obstruction prevents $H_{\H}$ from having any critical points at all. For the upper half-plane, the hydrodynamic Green function is exactly scale-invariant. Specifically, $g_{\H}(z, w) = \log \frac{|z-\bar{w}|}{|z-w|}$. If we apply a uniform dilation by a factor of $\lambda > 0$ mapping $(z_1, z_2) \mapsto (\lambda z_1, \lambda z_2)$, the mutual interaction remains unchanged:
  \[g_{\H}(\lambda z_1, \lambda z_2) = g_{\H}(z_1, z_2).\]
  However, the Robin function, which governs the interaction of a vortex with its own boundary image, is $r_{\H}(z) = \log(2\Im(z))$. Under the same dilation, it scales logarithmically:
  \[r_{\H}(\lambda z) = r_{\H}(z) + \log \lambda.\]
  The total half-plane Hamiltonian is
  \[H_{\H} = \gamma_1 \gamma_2 g_{\H}(z_1, z_2) + \frac{1}{2}\gamma_1^2 r_{\H}(z_1) + \frac{1}{2}\gamma_2^2 r_{\H}(z_2).\]
  Substituting the scaling relations yields:
  \[H_{\H}(\lambda z_1, \lambda z_2) = H_{\H}(z_1, z_2) + \frac{1}{2}(\gamma_1^2 + \gamma_2^2) \log \lambda\]
  Taking the (real) derivative with respect to $\lambda$ and evaluating at $\lambda = 1$ yields
  \[\sum_{i=1}^2 z_i \cdot d_{z_i} H_{\H} = \frac{1}{2}(\gamma_1^2 + \gamma_2^2).\]
  Because $\gamma_1, \gamma_2 \neq 0$, this sum is strictly positive. This proves that $d H_{\H}$ is non-zero everywhere in the configuration space, thus completing the proof for $H$ in the neighborhood of quadruple collisions.
\end{proof}

We may now harvest the opposite-sign case, where compactness (theorem~\ref{thm:compact}) and symplectic convexity (theorem~\ref{thm:Hpsh}) hold simultaneously.

\begin{theorem}[Periodic orbits in the opposite-sign problem]
  \label{thm:HZ}
  Assume $\gamma_1\gamma_2<0$. Then almost every level set of $H$ carries a periodic orbit.
\end{theorem}

\begin{proof}
  Let $E$ be a regular value. By theorem~\ref{thm:compact} and corollary~\ref{cor:apriori}, $\Sig_E$ is compact and contained in $C_2(D)$ at positive distance from $\partial C_2(D)$; so a whole interval of levels $\Sig_{E'}$, $|E'-E|<\eta$, sweeps a bounded open subset $\mathcal N$ of $C_2(D) \subset \C^2$. The symplectic form $\omega$ is constant in the coordinates $(z_1,z_2)$, hence linearly symplectomorphic to the standard form of $\R^4$; a bounded set has finite Hofer--Zehnder capacity, so $c_{HZ}(\mathcal N,\omega)<\infty$. The existence theorem for the Hofer--Zehnder capacity \cite{Hofer:1994,Zehnder:2010} then provides a closed characteristic on $\Sig_{E'}$ for almost every $E'$ in that interval.
\end{proof}

The compactness has a second, purely geometric consequence, which we record because it is what the contact-topological machinery of section~\ref{sec:quadruple} would need.

\begin{proposition}[Stein fillability in the opposite-sign problem]
  \label{prop:fillable}
  Assume $\gamma_1\gamma_2<0$ and let $E$ be a regular value. Then $\Omega_E:=\{H>E\}$ is a relatively compact Stein domain of $C_2(D)$ with strictly pseudoconvex boundary $\Sig_E$, and the induced $CR$ contact structure $\xi_E = T\Sig_E \cap J\,T\Sig_E$ makes $\Sig_E$ a closed, Stein fillable --- hence tight \cite{Eliashberg:1990,Eliashberg:1992} --- contact $3$-manifold.
\end{proposition}

\begin{proof}
  By theorem~\ref{thm:compact} and corollary~\ref{cor:apriori}, $\Omega_E$ is relatively compact in $C_2(D)$; by theorem~\ref{thm:Hpsh}, it is a sublevel set of the strictly plurisubharmonic function $-H$ on the Stein manifold $C_2(D)$ (lemma~\ref{lm:Stein}). Tightness of a holomorphically fillable contact structure is the Eliashberg--Gromov theorem.
\end{proof}

\begin{remark}[What is, and what is not, Reeb]
  \label{rem:notreeb}
  It is tempting to conclude from proposition~\ref{prop:fillable} that the vortex flow on $\Sig_E$ is a Reeb flow, and then to invoke the Weinstein conjecture --- proved in dimension $3$ by Taubes \cite{Taubes:2007} --- to get a periodic orbit on \emph{every} regular level rather than on almost every one. This is a trap, and it is worth saying exactly where it springs. The contact form associated with the filling is $\lambda = d^c(-H)$, and its differential $d\lambda$ is the Kähler form of the strictly plurisubharmonic function $-H$, whose coefficient matrix is $\bigl( \gamma_j\gamma_k K(z_j,z_k) \bigr)$. This is \emph{not} proportional to $\omega$, whose coefficient matrix is $\mathrm{diag}(\gamma_1,\gamma_2)$; indeed, when $\gamma_1\gamma_2<0$, $\omega$ is not even tamed by $J$. The Reeb field of $\lambda$ spans $\ker (d\lambda|_{\Sig_E})$, whereas $X_H$ spans $\ker (\omega|_{\Sig_E})$, and these are different line fields on $\Sig_E$.

  What would be needed is that $\Sig_E$ be of contact type \emph{for $\omega$ itself}, i.e. that there exist a Liouville vector field for $\omega$ transverse to $\Sig_E$. Pseudoconvexity does not provide one, precisely because $\omega$ is $J$-indefinite in the regime in which compactness holds. It is however sometimes available: proposition~\ref{prop:disc} exhibits such a Liouville field --- the Euler field of $\C^2$, halved --- for $D$ the disc, in the same-sign problem; remark~\ref{rem:nearlyround} and proposition~\ref{prop:contact-highE} below extend this to all energies in nearly round domains, and to high energies in arbitrary ones.
\end{remark}

\begin{question}
  \label{qu:contacttype}
  For which domains and intensities are the regular level sets of $H$ of contact type for $\omega$? Whenever they are and are compact --- for instance, presumably, in the opposite-sign problem --- Taubes' theorem upgrades the ``almost every level'' of theorem~\ref{thm:HZ} into ``every regular level''.
\end{question}

Partial answers in the same-sign problem are given in remark~\ref{rem:nearlyround} and section~\ref{sec:highE}.

\medskip The same two ingredients --- a closed contact manifold, and a Reeb
flow --- are what a \emph{global surface of section} would require, and it is
worth recording here what is missing, since the question is a global one and
belongs to this section rather than to the local analysis of the sections to
come.

\begin{remark}[The machinery one would like to use]
  \label{rem:HWZ}
  By the theory of Hofer, Wysocki and Zehnder \cite{Hofer:1998} and its refinements by Hryniewicz and Salom\~ao \cite{Hryniewicz:2011}, a non-degenerate Reeb flow on a closed, tight contact $3$-manifold admits, under suitable hypotheses on the linking numbers and Conley--Zehnder indices of a pair of orbits, a global surface of section. This section acts as a page in an open book decomposition whose binding consists of the two orbits. For flows on $S^3$, this forces the binding to be a Hopf link of linking number $1$. Because our energy levels $\Sigma_E$ exhibit different topologies, a more flexible tool is required. Recently, Colin, Dehornoy, Hryniewicz, and Rechtman \cite{Colin:2024} have shown that every non-degenerate Reeb flow on a closed $3$-manifold admits a \emph{broken} book decomposition, in which the pages are allowed to degenerate near the binding, bypassing the strict topological and linking constraints of classical open books.
\end{remark}

Of the hypotheses of remark~\ref{rem:HWZ}, closedness holds exactly in the
opposite-sign problem (theorem~\ref{thm:compact}) and tightness comes with it
(proposition~\ref{prop:fillable}); that the flow be Reeb is
question~\ref{qu:contacttype}. There remains non-degeneracy, and here is why one
may expect it to be harmless.

\begin{lemma}[Genericity of non-degeneracy; scheme]
  \label{lm:generic-nondeg}
  For a generic domain $D$, the periodic orbits of $X_H$ on regular level sets
  of $H$ are non-degenerate.
\end{lemma}

\begin{proof}[Proof scheme]
  The parameter space is the space of domain boundaries, which is infinite
  dimensional, and the relevant derivative is given by Hadamard's variational
  formula: an outward normal variation $\delta n$ of $\partial D$ changes the
  Green function by
  $\delta g(z,\zeta) = \frac{1}{2\pi}\oint_{\partial D}
  \partial_ng(z,\cdot)\,\partial_ng(\zeta,\cdot)\,\delta n\,\ddd s$ (with the
  outward normal; the sign is the one dictated by the monotonicity of the Dirichlet
  Green function with respect to the domain, $g_{0,D}\leq g_{0,D'}$ if $D\subset D'$
  \cite[Section~4.4]{Ransford:1995}, since
  both normal derivatives are negative; we write the formula for $m=1$, where
  $g = g_0$ is the Dirichlet Green function, the variation of the harmonic-measure
  terms having to be added when $m\geq 2$), hence changes
  the Hamiltonian by
  \[
    \delta H(z_1,z_2)= \frac{1}{4\pi}\oint_{\partial D}
    \Bigl( \frac{\partial \psi(\cdot\,;z_1,z_2)}{\partial n} \Bigr)^{2} \delta n \, \ddd s ,
    \qquad \psi = \gamma_1 g_{z_1}+\gamma_2 g_{z_2},
  \]
  the square arising because the cross terms and the two Robin terms recombine.
  Since $\psi$ is constant along $\partial D$, $\partial_n\psi$ is the tangential
  fluid velocity at the wall; by unique continuation for harmonic functions it
  cannot vanish on an arc. So $(\partial_n\psi)^2>0$ almost everywhere, and a
  bump $\delta n$ supported on a small arc injects a perturbation of $H$ which is
  localized in time along a given orbit. One checks that the resulting map, from
  boundary variations to $2$-jets of $H$ along a periodic orbit, is surjective;
  the Sard--Smale theorem then gives a residual set of domains for which no
  linearized return map has $1$ as an eigenvalue.

  We have not written out the verification of surjectivity, which is the only
  non-formal step.
\end{proof}

\begin{remark}
  A round domain is of course not generic: there the rotational symmetry
  produces circles of periodic orbits. The degeneracy is however of
  Morse--Bott type, the degenerate directions coming exclusively from the $S^1$
  families, and Morse--Bott versions of the contact machinery exist.
\end{remark}

In the same-sign problem there is one more, and more serious, obstruction:
$\Sig_E$ is then open, with $m$ ends\footnote{We use \emph{end} in the sense of Freudenthal. Let $X$ be a connected, locally compact, $\sigma$-compact Hausdorff space (here an open manifold). For a compact subset $C\subset X$, let $\pi_0^\infty(X\setminus C)$ be the set of connected components of $X\setminus C$ with non-compact closure; if $C\subset C'$, each element of $\pi_0^\infty(X\setminus C')$ is contained in a unique element of $\pi_0^\infty(X\setminus C)$, which defines a map $\pi_0^\infty(X\setminus C')\to\pi_0^\infty(X\setminus C)$. The set of \emph{ends} of $X$ is the inverse limit $\varprojlim_C\pi_0^\infty(X\setminus C)$ over the compact subsets of $X$ ordered by inclusion; concretely, an end is a choice of an element of $\pi_0^\infty(X\setminus C)$ for every compact $C$, compatible with these maps. For instance $\R$ has two ends, $\R^n$ ($n\geq2$) has one, and a closed surface punctured at $k$ points has $k$ ends. An end is \emph{collared} if it has a neighbourhood diffeomorphic to $N\times[0,\infty)$, with $N$ a closed manifold. By proposition~\ref{prop:ends}, the $m$ ends of $\Sig_E$ are collared, with cross-section $\T^2$: the complement of a compact subset of $\Sig_E$ is a disjoint union of $m$ copies of $\T^2\times(0,1)$, and of each copy only the side $q\to0$, towards the quadruple collisions, is an end of $\Sig_E$.} which we will describe in
section~\ref{sec:quadruple}, so that ``global'' has to be given a meaning.
Moreover the two natural candidates for the binding --- the two extreme modes of
periodic motion at a given energy, namely a widely separated pair rotating about
the centre of the domain, and a tightly bound pair whose centre of vorticity
slowly revolves --- will only be produced perturbatively, in
sections~\ref{sec:collision} and~\ref{sec:near-domain-boundary}, and not at a
common energy. We therefore state the following as a question.

\begin{question}
  \label{qu:Birkhoff}
  For a generic domain, does the flow of $X_H$ admit a global surface of section
  on every regular level $\Sig_E$, and in particular a Birkhoff annulus of
  section $\mathcal A_E$ whose binding consists of the two extreme periodic
  orbits described above? Equivalently, does $\Sig_E$ carry an open book
  decomposition with annular pages adapted to the flow?
\end{question}

Should the answer be positive, the return map to $\mathcal A_E$ would be an
area-preserving twist-like map of the annulus, whose two boundary circles are
the two extreme modes; and, in the same-sign problem, the invariant tori of
theorem~\ref{thm:kam} would appear as invariant circles accumulating on one of
its ends. Much of what follows may be read as an attempt to see, perturbatively,
the two ends of that hypothetical annulus. At high energy, theorem~\ref{thm:annulus} answers the question positively when $D$ is simply connected and its Robin function has a unique critical point, but in a different form: the annulus is bounded by the binary orbit on one side, and ends at the quadruple collisions on the other.

\begin{remark}[The dichotomy]
  \label{rem:dichotomy}
  It is worth recording the contrast. Plurisubharmonicity holds for all intensities but is of little use without compactness; compactness holds only in the opposite-sign problem. So in the opposite-sign problem the level sets are closed, tight, Stein fillable contact $3$-manifolds, on almost all of which a periodic orbit exists. In the same-sign problem, energy levels are open $3$-manifolds with ends, and global existence has to be replaced by the perturbative analysis of the ends carried out in sections~\ref{sec:collision} to~\ref{sec:quadruple}. The local non-degeneracy results, on the contrary, are insensitive to the signs.
\end{remark}

\section{Near interior binary collisions}
\label{sec:collision}

In this section, keeping the conventions of the prior section, we study the regime in which the pair is close to a binary collision, away from the domain boundary. 

% The regime is a high-$|$energy$|$ regime: the natural parameter is $\eps=e^{-E/(\gamma_1\gamma_2)}\to0$. We show that $\theta=\arg(z_1-z_2)$ is a fast angle with frequency $\Gamma/\rho^2$, average it away to all orders (\S\ref{sec:nf}), and obtain that the guiding centre follows, to all orders in $\eps$, a one-degree-of-freedom flow whose leading Hamiltonian is $\tfrac{\Gamma^2}{2}r$, $r$ the Robin function. Two consequences: every non-degenerate critical point of $r$ generates a one-parameter family of periodic orbits accumulating on $\Delta$ (\S\ref{sec:po}), with the stability type of the critical point; and, under the twist condition that the Hamiltonian flow of $q^{2\kappa}$ ($q=e^{r}$ the harmonic radius, $\kappa=\Gamma^2/2\gamma_1\gamma_2$) be non-isochronous, invariant tori accumulate on $\Delta$ along the level curves of $r$ (\S\ref{sec:kam}). Near a non-degenerate maximum of $r$ the twist condition is automatic. The degenerate case $\Gamma=0$ is treated separately in \S\ref{sec:dipole}: there the fast motion is a translation, not a rotation, and none of the above applies.

Define 
\[
  \Gamma=\gamma_1+\gamma_2,\qquad \mu=\frac{\gamma_1\gamma_2}{\Gamma},\qquad
  \kappa=\frac{\Gamma^2}{2\gamma_1\gamma_2}.
\]
We first assume that $\Gamma \neq 0$, so that $\kappa \neq 0$; the degenerate case $\Gamma=0$, that of a dipole, is of a different nature and is investigated separately in section~\ref{sec:dipole}. Note that $\kappa \geq 2$ in the same-sign problem, while $\kappa<0$ in the opposite-sign one; none of what follows depends on which.

The small parameter of the whole section is the following.

\begin{definition}
  \label{def:eps}
  For $E\in\R$ (energy value) put
  \[\eps:=e^{-E/(\gamma_1\gamma_2)}.\]
  The interior binary-collision regime is $\eps\to0^+$, i.e.\ $E\to \sgn(\gamma_1\gamma_2) \infty$.
\end{definition}

This section is devoted to the following two theorems; the motions they provide are described in physical terms in the introduction, after theorem~\ref{thm:B}.

First note that, since $q\to0$ on $\partial D$, $r$ attains a maximum in $D$.

\begin{theorem}[Periodic orbits in the near collision regime]
  \label{thm:periodic}
  \begin{enumerate}
  \item If the maximum of $r$ is non-degenerate, there is a real analytic family of elliptic periodic orbits accumulating on $\Delta$ as $\eps\to0$.
  \item If $r$ is a Morse function, then for $m\ge2$ there are at least $m-1$ corresponding real analytic families of hyperbolic periodic orbits. 
  \end{enumerate}
\end{theorem}

See also \cite{Bartsch:2017} for periodic orbits in $N$-vortex type problems. 

For a round annulus $r$ is radial, so its critical set is a whole circle and theorem~\ref{thm:periodic} does not apply: instead of isolated periodic orbits one expects a single degenerate, normally elliptic invariant $2$-torus, made of the pairs whose centre of vorticity travels along the critical circle.

On the other hand, for $D$ strictly convex (hence simply connected) the Robin function is strictly concave (see \cite{Pommerenke:1975} or, for example, \cite[Chap. 8]{Flucher:1999}, with the opposite sign convention), so $r$ has a unique non-degenerate maximum and the first item of the theorem applies. Also, for a generic domain the Robin function is Morse~\cite{Bartsch:2018}, in which case neither item of the theorem is empty. 

\medskip In order to state the second theorem, we need some additional notations. 
The function that governs everything is the smooth, positive function on $D$
\begin{equation}
  \label{eq:Theta}
  \Theta:=q^{2\kappa}=\exp\Bigl(\frac{\Gamma^2}{\gamma_1\gamma_2}\,r\Bigr).
\end{equation}

Let $\mathcal U\subset D$ be an open annulus foliated by closed regular level
curves of $\Theta$ (equivalently of $r$). Let $A$ be the action of these curves for the symplectic form $\Gamma\,\ddd X\wedge \ddd Y$ of the centre of vorticity, namely
\[
  A:=\frac{\Gamma}{2\pi}\Area\bigl(\text{enclosed region}\bigr),
  \qquad \Gamma\,\ddd X\wedge \ddd Y=\ddd A\wedge \ddd\psi,
\]
$\psi$ being the conjugate angle,
and regard $\Theta$ as a function $\Theta(A)$ on an interval $\mathcal I$. We say $\mathcal U$ is \emph{non-isochronous} if
\[\Theta''(A)\neq0 \qquad\text{for all }A\in\mathcal I,\]
or, equivalently, $r''(A)+\tfrac{\Gamma^2}{\gamma_1\gamma_2}\,r'(A)^2\neq0$, or, equivalently, the Hamiltonian flow of $\Theta$ has a non-constant period map on $\mathcal U$. (In lemma~\ref{lem:iso}, this condition will turn out to be equivalent to the isoenergetic non-degeneracy of some Hamiltonian which will be later introduced.)

\begin{theorem}[Invariant tori accumulating on $\Delta$]
  \label{thm:qp}
  Let $\mathcal U\Subset D$ be any non-isochronous annulus of closed regular level
  curves of $r$. Close enough to $\Delta \cap \mathcal{U}^2$, that is for $\eps$ small enough:
  \begin{enumerate}[label=(\roman*)]
  \item In each energy level, the flow of $X_H$ has invariant real analytic Lagrangian $2$-tori, graphs over the angles, carrying quasi-periodic motions, and of relative $3$-dimensional Lebesgue measure $1-O(e^{-c/\eps^2})$ for some $c>0$.
  \item These tori separate the energy level inside $\mathcal U\times \SS^1$, so that an orbit caught between two of them stays there for all time; in particular its two vortices neither collide nor reach $\partial D$.
  \item Between two consecutive invariant tori, there exist periodic orbits and Mather sets of every intermediate rotation number.
  \item As $\eps \to 0$, the tori lying over a fixed level curve $\mathcal{O}\subset \mathcal U$ converge, for the Hausdorff distance in $\bar D \times \bar D$, to the circle $\{(Z,Z), \; Z \in \mathcal{O}\}$ of the diagonal.
  \end{enumerate}
\end{theorem}

See \cite{Khanin:1981} for an analogous result for $N$ vortices of positive intensities in the plane. Merging the two proofs, one should expect an analogous result for $N$ vortices in a bounded domain, with arbitrary intensities such that $\gamma_1 + \cdots + \gamma_n \neq 0$ for all $n=2,...,N$.

We now prove theorems~\ref{thm:periodic} and \ref{thm:qp}.

Define the center of vorticities and the relative position
\[Z=\frac{\gamma_1z_1+\gamma_2z_2}{\Gamma},\qquad u=z_1-z_2=\rho e^{i\theta}.\]
A direct computation shows that
\begin{equation}
  \label{eq:omega}
  \omega=\frac{\Gamma}{2i}\,\ddd \bar Z\wedge \ddd Z+\ddd I\wedge \ddd\theta
  = \Gamma\,\ddd X \wedge \ddd Y + \ddd I\wedge \ddd\theta,
  \qquad I=\tfrac{\mu}{2}\rho^2 ,
\end{equation}
the cross terms cancelling because $Z$ is the centre of vorticity. Writing
$F(z_1,z_2)=\gamma_1\gamma_2h+\tfrac{\gamma_1^2}{2}r(z_1)+\tfrac{\gamma_2^2}{2}r(z_2)$
and $F(Z,Z)=\tfrac{\Gamma^2}{2}r(Z)$, set $W:=F(z_1,z_2)-F(Z,Z)$, so that
\begin{equation}
  \label{eq:split}
  H=\underbrace{-\frac{\gamma_1\gamma_2}{2}\log\frac{2I}{\mu}}_{=:H_{\mathrm f}(I)}
  +\underbrace{\frac{\Gamma^2}{2}r(Z)}_{=:H_{\mathrm s}(Z)}+\;W .
\end{equation}

\begin{lemma}[Taylor structure of $W$]
  \label{lem:taylor}
  At fixed $Z$ in a compact $K\Subset D$, uniformly in $K$,
  \[
    W=\gamma_1\gamma_2\,\Ree\bigl(P(Z)u^2\bigr)+O(|u|^3),
    \qquad P(Z)=\partial^2_{z_1}h(z_1,Z)\big|_{z_1=Z},
  \]
  the term linear in $u$ being identically zero; moreover
  $\langle W\rangle_\theta=O(|u|^4)$, where $\langle\cdot\rangle_\theta$ is the
  average over $\theta$. Equivalently, in the variables \eqref{eq:omega},
  \begin{equation}
    \label{eq:W2}
    W=2\Gamma I\,\Ree\bigl(P(Z)e^{2i\theta}\bigr)+O(I^{3/2}),
    \qquad \langle W\rangle_\theta=O(I^2).
  \end{equation}
\end{lemma}

\begin{proof}
  With $a=\gamma_2/\Gamma$, $b=-\gamma_1/\Gamma$
  the coefficient of the linear term is
  $\tfrac{Dr}{2\Gamma}(\gamma_1\gamma_2^2-\gamma_1^2\gamma_2+\gamma_1^2\gamma_2-\gamma_1\gamma_2^2)=0$;
  in the quadratic term the coefficient of $D^2_{12}h$ vanishes and that of
  $D^2_{11}h$ is $\gamma_1\gamma_2$; harmonicity of $h(\cdot,Z)$ turns
  $\tfrac12 D^2_{11}h[u,u]$ into $\Ree(Pu^2)$. The quadratic term has zero
  $\theta$-average and the cubic term is odd under $\theta\mapsto\theta+\pi$.
  Finally $|u|^2=2I/\mu$ and $\Ree(Pu^2)=\tfrac{2I}{\mu}\Ree(Pe^{2i\theta})$,
  and $\gamma_1\gamma_2/\mu=\Gamma$.
\end{proof}

\begin{proposition}[The energy level as a graph]
  \label{prop:graph}
  Fix $K\Subset D$ compact. For $\eps$ small, if $\Sigma_E = H^{-1}(E)$, $\Sig_E\cap\{Z\in K\}$ is the graph over $K_Z \times \SS^1_\theta$ of equation
  \begin{equation}
    \label{eq:rho-graph}
    \rho=\eps\,q(Z)^{\kappa}\,e^{W/(\gamma_1\gamma_2)}
    = \eps\,q(Z)^\kappa\Bigl(1+\eps^2q(Z)^{2\kappa}\Ree\bigl(P(Z)e^{2i\theta}\bigr)+O(\eps^3)\Bigr).
  \end{equation}
  In particular $\rho = O(\eps)$ on $K$, and on $\Sig_E$ the separation
  oscillates at twice the rotation frequency, with relative amplitude
  $\eps^2q^{2\kappa}|P(Z)|$.
\end{proposition}

\begin{proof}
  Solve $H=E$ in \eqref{eq:split} for $\log\rho$:
  $\log\rho=(\gamma_1\gamma_2)^{-1}\bigl(\tfrac{\Gamma^2}{2}r+W-E\bigr)$, i.e.\
  \eqref{eq:rho-graph}; $W=O(\rho^2)=O(\eps^2)$ by Lemma~\ref{lem:taylor}, and
  the graph is well defined because
  $\rho\,\partial_\rho H=-\gamma_1\gamma_2+O(\eps^2)\neq0$. Substituting
  \eqref{eq:W2} gives \eqref{eq:rho-graph}.
\end{proof}

% \begin{remark}
%   Formula \eqref{eq:rho-graph} is the same graph as in the quadruple-collision
%   regime, where $\rho=R_Eq^{\kappa}(1+O(q^{2\beta}))$ with
%   $R_E=e^{-E/(\gamma_1\gamma_2)}=\eps$. The two regimes are the two ways of
%   making $\rho$ small on that graph: $q\to0$ (boundary) or $\eps\to0$
%   (interior). Here $q$ stays bounded away from $0$ and $\eps$ is the small
%   parameter. Note that for $\gamma_1\gamma_2<0$ one has $\kappa<0$, so $\rho$
%   \emph{grows} near $\partial D$ along $\Sig_E$; this is why the restriction to
%   a compact $K$ is not cosmetic in that case.
% \end{remark}

\begin{proposition}[Skew-product structure]
  \label{prop:skew}
  On $\Sig_E\cap\{Z\in K\}$,
  \begin{equation}
    \label{eq:frequencies}
    \dot\theta = \frac{\Gamma}{\rho^2}+O(1) = \frac{\Gamma}{\eps^2}q(Z)^{-2\kappa}\bigl(1+O(\eps^2)\bigr),
    \qquad
    \dot Z=X_{H_{\mathrm s}}(Z)+O(\eps^2), \qquad X_{H_{\mathrm s}}=\tfrac12\Gamma\,\nabla^{\perp}r ,
  \end{equation}
  where $X_{H_{\mathrm s}}$ is the Hamiltonian vector field of $H_{\mathrm s}$ for the symplectic form $\Gamma\,\ddd X\wedge \ddd Y$ of the guiding centre, and $\nabla^\perp=(\partial_Y,-\partial_X)$. Thus $\theta$ is a fast angle
  with frequency $O(\eps^{-2})$, and the guiding centre follows, at leading
  order, the Kirchhoff--Routh flow of a single vortex of circulation
  $\Gamma$.
\end{proposition}

\begin{proof}
  $\dot\theta=-\partial_IH=\tfrac{\gamma_1\gamma_2}{2I}-\partial_IW
  =\tfrac{\Gamma}{\rho^2}+O(1)$, using $\gamma_1\gamma_2/\mu=\Gamma$ and
  $\partial_IW=O(1)$ by \eqref{eq:W2}. The slow equations are
  $\Gamma\dot X=\partial_YH$, $\Gamma\dot Y=-\partial_XH$ with
  $H=H_{\mathrm s}(Z)+O(\eps^2)$ and $H_{\mathrm s}=\tfrac{\Gamma^2}{2}r$.
\end{proof}

Everything below rests on \eqref{eq:frequencies}. The hypothesis $\Gamma\neq0$
enters exactly once, in $\dot\theta\sim\Gamma/\rho^2$: it is what makes $\theta$
a fast \emph{angle}. The signs of $\gamma_1,\gamma_2$ never enter.

Since $h$ and $r$ are real analytic on $D^2$ and on $D$ (definition~\ref{def:Robin}), whatever the regularity of $\partial D$, the Hamiltonian is real analytic on the whole regime, and we estimate it in analytic norms adapted to \eqref{eq:omega}. Fix a compact interval $[S_-,S_+]\subset\R\setminus\{0\}$ whose interior contains $\tfrac\mu2\,q(Z)^{2\kappa}$ for every $Z\in K$, so that, by proposition~\ref{prop:graph}, $I/\eps^2\in[S_-,S_+]$ on $\Sig_E\cap\{Z\in K\}$ for $\eps$ small. For $\sigma>0$ let
\[
  \mathcal B_\sigma:=\bigl\{(X,Y,S,\theta)\in\C^3\times\C/2\pi\Z:\
  \dist((X,Y),K)<\sigma,\ \dist(S,[S_-,S_+])<\sigma,\ |\Imm\theta|<\sigma\bigr\}
\]
(with $K\subset\R^2\subset\C^2$), and, for a function $f$ of $(X,Y,I,\theta)$,
\[
  |f|_\sigma:=\sup_{(X,Y,S,\theta)\in\mathcal B_\sigma}\bigl|f(X,Y,\eps^2S,\theta)\bigr| ,
\]
i.e.\ the action $I$ is measured in its natural scale $\eps^2$. All the functions below are real analytic: holomorphic on some $\mathcal B_\sigma$ and real at real points. We denote by $X_g$ the Hamiltonian vector field of $g$ ($\iota_{X_g}\omega=\ddd g$) and by
\[
  \{f,g\}:=\ddd f(X_g)=\Gamma^{-1}\bigl(\partial_Xf\,\partial_Yg-\partial_Yf\,\partial_Xg\bigr)
  +\partial_If\,\partial_\theta g-\partial_\theta f\,\partial_Ig
\]
the Poisson bracket, so that $\frac{\ddd}{\ddd t}\,f\circ\Phi^t_g=\{f,g\}\circ\Phi^t_g$ along the flow $\Phi^t_g$ of $X_g$.

\begin{lemma}[Analytic estimates]
  \label{lem:analytic}
  There are $\sigma_0>0$, $\eps_0>0$ and $C>0$ such that, for $0<\eps<\eps_0$:
  \begin{enumerate}[label=(\roman*)]
  \item $H_{\mathrm s}$ and $W$ are holomorphic on $\mathcal B_{\sigma_0}$, the first derivatives of $H_{\mathrm s}$ are bounded there by $C$, and
    \[
      |W|_{\sigma_0}\le C\eps^2,\qquad |\langle W\rangle_\theta|_{\sigma_0}\le C\eps^4 ;
    \]
  \item if $0<\delta<\sigma\le\sigma_0$, then $|\partial_Xf|_{\sigma-\delta}$, $|\partial_Yf|_{\sigma-\delta}$ and $|\partial_\theta f|_{\sigma-\delta}$ are at most $\delta^{-1}|f|_\sigma$, while $|\partial_If|_{\sigma-\delta}\le\eps^{-2}\delta^{-1}|f|_\sigma$; consequently
    \begin{equation}
      \label{eq:bracket}
      |\{f,g\}|_{\sigma-\delta}\le C\,\eps^{-2}\delta^{-2}\,|f|_{\sigma}\,|g|_{\sigma} ;
    \end{equation}
  \item if $\langle f\rangle_\theta=0$, the primitive $\partial_\theta^{-1}f$ of $f$ with respect to $\theta$ which has zero average satisfies $|\partial_\theta^{-1}f|_\sigma\le\pi\,|f|_\sigma$.
  \end{enumerate}
\end{lemma}

\begin{proof}
  (i) The function $F$ is real analytic on $D\times D$, hence extends holomorphically to a complex neighbourhood in $\C^4$ of $K'\times K'$, where $K'\Subset D$ is a compact neighbourhood of $K$. On $\mathcal B_\sigma$, for $\sigma$ small, $\rho=(2\eps^2S/\mu)^{1/2}$ is holomorphic and $|\rho|\le C\eps$, so that the two real coordinates $\rho\cos\theta$, $\rho\sin\theta$ of $u$ have modulus at most $C\eps\,e^\sigma$. Hence, for $\sigma_0$ and $\eps_0$ small, $W=F(Z+\tfrac{\gamma_2}\Gamma u,Z-\tfrac{\gamma_1}\Gamma u)-F(Z,Z)$ is holomorphic on $\mathcal B_{\sigma_0}$. By lemma~\ref{lem:taylor}, its Taylor expansion with respect to the coordinates of $u$ starts at order $2$, and that of its average at order $4$ (these are identities between power series, which persist on the complex domain); the two bounds follow from Taylor's formula.

  (ii) These are Cauchy's inequalities, the disc in the variable $I$ having radius $\eps^2\delta$.

  (iii) By Cauchy's theorem, $g=\partial_\theta^{-1}f$ has zero average on each circle $\{\Imm\theta=\mathrm{const}\}$, and two points of such a circle are joined by an arc of length at most $\pi$, along which $|\partial_\theta g|\le|f|_\sigma$.
\end{proof}

\begin{theorem}[Normal form with an exponentially small remainder]
  \label{thm:nf}
  Fix $K\Subset D$. There are $\sigma_*>0$, $c>0$, $C>0$ and $\eps_0>0$ such that, for $0<\eps<\eps_0$, there is a real analytic symplectic map $\Phi:\mathcal B_{\sigma_*}\to\mathcal B_{\sigma_0}$, whose image contains $\mathcal B_{\sigma_*/2}$, satisfying
  \[
    \Phi=\mathrm{id}+\bigl(O(\eps^4),\,O(\eps^4),\,O(\eps^4),\,O(\eps^2)\bigr)
    \quad\text{in }(X,Y,I,\theta),
  \]
  and such that
  \begin{equation}
    \label{eq:nf}
    H\circ\Phi=H_{\mathrm f}(I)+H_{\mathrm s}(Z)+\mathcal Z(Z,I)+\Rem(Z,I,\theta),
  \end{equation}
  with
  \[
    |\mathcal Z|_{\sigma_*}\le C\,\eps^4,
    \qquad
    |\Rem|_{\sigma_*}\le C\,e^{-c/\eps^2} .
  \]
\end{theorem}

The scheme of proof is that of Neishtadt \cite{Neishtadt:1984} for systems with a single rapidly rotating phase: a few averaging steps on a fixed domain make the remainder small enough, and are followed by a number of steps of order $\eps^{-2}$, on domains shrinking linearly.

\begin{proof}
  Throughout, $C$ denotes various constants which depend neither on $\eps$ nor on the step.

  \emph{One step.} Suppose that, on some $\mathcal B_\sigma$, $\sigma\le\sigma_0$,
  \[
    H=H_{\mathrm f}(I)+\mathsf h(Z,I)+R(Z,I,\theta),\qquad\langle R\rangle_\theta=0 ,
  \]
  and let $\mathsf a(\mathsf h)_\sigma:=|\partial_X\mathsf h|_\sigma+|\partial_Y\mathsf h|_\sigma+|\partial_I\mathsf h|_\sigma$. Let $\Omega_\theta=\gamma_1\gamma_2/(2I)$ be the fast frequency, and
  \[
    \chi:=\Omega_\theta^{-1}\,\partial_\theta^{-1}R .
  \]
  Since $H_{\mathrm f}'=-\Omega_\theta$, this function solves the homological equation $\{H_{\mathrm f},\chi\}=-\Omega_\theta\,\partial_\theta\chi=-R$. \emph{No small divisors occur}, because only the fast frequency is used. As $|\Omega_\theta^{-1}|_\sigma\le C\eps^2$, lemma~\ref{lem:analytic}(iii) gives $|\chi|_\sigma\le C\eps^2|R|_\sigma$. Let $0<\delta<\sigma/2$. By lemma~\ref{lem:analytic}(ii), on $\mathcal B_{\sigma-\delta/2}$ the vector field $X_\chi$ has components of size at most $C\eps^2|R|_\sigma/\delta$ along $X$ and $Y$, and at most $C|R|_\sigma/\delta$ along $S=I/\eps^2$ and $\theta$. Hence there is $c_1>0$ such that, if
  \begin{equation}
    \label{eq:smallstep}
    |R|_\sigma\le c_1\,\delta^2 ,
  \end{equation}
  then the maps $\Phi^t_\chi$, $0\le t\le1$, send $\mathcal B_{\sigma-\delta}$ into $\mathcal B_{\sigma-\delta/2}$. Using the homological equation,
  \[
    H\circ\Phi^1_\chi=H_{\mathrm f}+\mathsf h+R^+,\qquad
    R^+=\int_0^1\{\mathsf h+tR,\chi\}\circ\Phi^t_\chi\,\ddd t .
  \]
  Since $\mathsf h$ does not depend on $\theta$, its bracket with $\chi$ involves the derivatives of $\chi$ with respect to $X$, $Y$ and $\theta$ only, so $|\{\mathsf h,\chi\}|_{\sigma-\delta/2}\le C\,\mathsf a(\mathsf h)_\sigma\,\eps^2|R|_\sigma/\delta$; and $|\{R,\chi\}|_{\sigma-\delta/2}\le C\delta^{-2}|R|_\sigma^2$ by \eqref{eq:bracket}. Therefore
  \begin{equation}
    \label{eq:step}
    |R^+|_{\sigma-\delta}\le C\Bigl(\mathsf a(\mathsf h)_\sigma\,\frac{\eps^2}{\delta}+\frac{|R|_\sigma}{\delta^2}\Bigr)|R|_\sigma .
  \end{equation}
  We then replace $\mathsf h$ by $\mathsf h^+=\mathsf h+\langle R^+\rangle_\theta$ and $R^+$ by $R^+-\langle R^+\rangle_\theta$, which at most doubles the norm of the latter, and, by lemma~\ref{lem:analytic}(ii),
  \begin{equation}
    \label{eq:astep}
    \mathsf a(\mathsf h^+)_{\sigma-2\delta}\le\mathsf a(\mathsf h)_\sigma+C\,\eps^{-2}\delta^{-1}\,|R^+|_{\sigma-\delta} .
  \end{equation}

  \emph{Preliminary steps.} We start from $\mathsf h_0=H_{\mathrm s}+\langle W\rangle_\theta$ and $R_0=W-\langle W\rangle_\theta$. By lemma~\ref{lem:analytic}, $|R_0|_{\sigma_0}\le C\eps^2$ and, up to replacing $\sigma_0$ by $\sigma_0/2$, $\mathsf a(\mathsf h_0)_{\sigma_0}\le\mathsf a_0$ for some $\mathsf a_0$ independent of $\eps$ (indeed $|\partial_I\langle W\rangle_\theta|\le C\eps^{-2}\eps^4$). Two steps with $\delta=\sigma_0/8$ are licit for $\eps$ small, since \eqref{eq:smallstep} holds, and by \eqref{eq:step} and \eqref{eq:astep} they lead to
  \[
    |R_1|_{3\sigma_0/4}\le C\eps^4,\qquad |R_2|_{\sigma_0/2}\le C\eps^6,\qquad
    \mathsf a(\mathsf h_2)_{\sigma_0/2}\le\mathsf a_0+C\eps^2 .
  \]

  \emph{Main steps.} Let $N$ be the integer part of $c_3/\eps^2$, and perform $N$ more steps with $\delta=\sigma_0/(8N)$, from $\mathcal B_{\sigma_0/2}$ down to $\mathcal B_{\sigma_0/4}$. We claim that, if $c_3$ and then $\eps$ are small enough, then for $0\le k\le N$
  \[
    |R_{2+k}|\le2^{-k}|R_2|\le C\eps^62^{-k}\qquad\text{and}\qquad\mathsf a(\mathsf h_{2+k})\le2\mathsf a_0
  \]
  on the respective domains. Indeed, if this holds up to $k$, then \eqref{eq:smallstep} holds since $\eps^6\ll\delta^2\asymp\eps^4$; the factor in \eqref{eq:step} is at most $C\bigl(2\mathsf a_0\,\eps^2N+\eps^6N^2\bigr)\le C(c_3+c_3^2\eps^2)\le\tfrac14$; and \eqref{eq:astep} gives $\mathsf a(\mathsf h_{3+k})\le\mathsf a(\mathsf h_{2+k})+C\eps^{-2}N\eps^62^{-k}\le\mathsf a(\mathsf h_{2+k})+C\eps^22^{-k}$, whence $\mathsf a(\mathsf h_{3+k})\le\mathsf a_0+C\eps^2\le2\mathsf a_0$.

  \emph{Conclusion.} Let $\sigma_*=\sigma_0/4$, let $\Phi$ be the composition of the $N+2$ time-one maps, $\Rem=R_{2+N}$ and $\mathcal Z=\mathsf h_{2+N}-H_{\mathrm s}$. Then $|\Rem|_{\sigma_*}\le C\eps^62^{-N}\le Ce^{-c/\eps^2}$, and $|\mathcal Z|_{\sigma_*}\le|\langle W\rangle_\theta|_{\sigma_0}+\sum_{k\ge1}|R_k|\le C\eps^4$. The displacements add up: along $X$ and $Y$, to at most $C\eps^2\sum_k|R_k|/\delta_k\le C\eps^4$; along $S$ and $\theta$, to at most $C\sum_k|R_k|/\delta_k\le C\eps^2$, that is, $O(\eps^4)$ for $I$ and $O(\eps^2)$ for $\theta$. Being $O(\eps^2)$-close to the identity in the variables $(X,Y,S,\theta)$, $\Phi$ is a diffeomorphism onto its image, which contains $\mathcal B_{\sigma_*/2}$. Finally each $\chi$ is real at real points, so $\Phi$ is real analytic, and it is symplectic as a composition of Hamiltonian flows.
\end{proof}

\begin{remark}[Form of the first correction]
  $\mathcal Z=\mathcal Z_2+O(\eps^6)$ where $\mathcal Z_2$ is the sum of
  $\langle W_4\rangle_\theta$ (the $\theta$-average of the quartic Taylor term
  of $W$) and of a second-order averaging term built from the strain
  $W_2=2\Gamma I\Ree(Pe^{2i\theta})$, which is quadratic in $P(Z)$ and
  proportional to $\Gamma^2(\gamma_1\gamma_2)^{-1}|P(Z)|^2I^2$ up to an explicit
  rational constant. We do not use the constant, only
  $\mathcal Z_2=O(I^2)$ and its smoothness, so we have not verified it.
\end{remark}

\begin{corollary}[$I$ is an adiabatic invariant over exponentially long times]
  \label{cor:adiabatic}
  There are $c>0$ and $C>0$ such that, along any orbit of the regime with $Z(s)\in K$ for
  $|s|\le|t|$,
  \[
    \bigl|\rho(t)-\rho(0)\bigr|\le C\bigl(\eps^{3}+\eps^{-1}e^{-c/\eps^2}|t|\bigr).
  \]
  In particular the pair does not collide, and its size does not change by more
  than $O(\eps^3)$, on the time scale $|t|\le e^{c/(2\eps^2)}$, for
  \emph{every} initial condition of the regime.
\end{corollary}

\begin{proof}
  Apply theorem~\ref{thm:nf} to a compact neighbourhood of $K$ in $D$. In the normal-form variables, $\dot I=\partial_\theta\Rem=O(e^{-c/\eps^2})$ by lemma~\ref{lem:analytic}(ii),
  and $\Phi$ changes $I$ by $O(\eps^4)$, i.e.\ $\rho$ by $O(\eps^3)$, since
  $\rho^{-1}=O(\eps^{-1})$.
\end{proof}

\begin{remark}
  Corollary~\ref{cor:adiabatic} is an all-orbits statement, over a finite, though exponentially long, time.
  Theorem~\ref{thm:qp} below upgrades it to an all-time statement, but only in
  the region trapped by invariant tori.
\end{remark}

By Proposition~\ref{prop:skew}, $\dot\theta\neq0$ on the whole regime, so
\[
  \mathcal A_E:=\Sig_E\cap\{Z\in K\}\cap\{\theta=0\}\ \cong\ K
\]
is a global section, with (positive) return time
$T(Z)=2\pi/|\dot\theta|=\tfrac{2\pi\eps^2}{|\Gamma|}q(Z)^{2\kappa}(1+O(\eps^2))$.
Since $\ddd\theta=0$ on the section, $\omega$ restricts to $\Gamma\,\ddd X\wedge \ddd Y$
and the return map $\mathcal P_\eps$ preserves it.

\begin{proposition}[The return map is a near-identity twist map]
  \label{prop:return}
  On $K$,
  \begin{equation}
    \label{eq:P}
    \mathcal P_\eps=\mathrm{id}+\frac{\pi\eps^2|\Gamma|}{2\kappa}\,X_{\Theta}+O(\eps^4),
    \qquad X_\Theta=\Gamma^{-1}\nabla^\perp\Theta,
  \end{equation}
  and more precisely, there are $c>0$ and a real analytic diffeomorphism $\Psi_\eps$ from a neighbourhood of $K$ onto its image, which preserves $\Gamma\,\ddd X\wedge\ddd Y$ and is $O(\eps^4)$-close to the identity, such that $\Psi_\eps^{-1}\circ\mathcal P_\eps\circ\Psi_\eps$ differs by $O(e^{-c/\eps^2})$, uniformly on a complex neighbourhood of $K$, from the time-$\pi\eps^2$ map of the Hamiltonian
  \[G_\eps=\frac{|\Gamma|}{2\kappa}\,\Theta\,\bigl(1+O(\eps^4)\bigr)
    \qquad\text{with respect to }\ \Gamma\,\ddd X\wedge \ddd Y.\]
\end{proposition}

\begin{proof}
  To leading order $\mathcal P_\eps(Z)=Z+T(Z)\,X_{H_{\mathrm s}}(Z)+O(\eps^4)$
  by Proposition~\ref{prop:skew}. Now
  \[
    T\cdot\tfrac{\Gamma}{2}\nabla^\perp r
    =\frac{2\pi\eps^2}{|\Gamma|}\Theta\cdot\frac{\Gamma}{2}\nabla^\perp r
    =\pi\eps^2|\Gamma|\,\Theta\,\Gamma^{-1}\nabla^\perp r
    =\frac{\pi\eps^2|\Gamma|}{2\kappa}\,X_\Theta,
  \]
  because $\nabla\Theta=2\kappa\,\Theta\,\nabla r$.

  For the refined statement, consider first the truncated normal form $\bar H=H_{\mathrm f}(I)+H_{\mathrm s}(Z)+\mathcal Z(Z,I)$ of theorem~\ref{thm:nf}. Its energy level $\{\bar H=E\}$ is the graph $I=I_E(Z)$, with ($\mathcal Z$ being evaluated at $I=I_E(Z)$)
  \[
    \frac{2I_E}{\mu}=\eps^2\,\Theta\,e^{2\mathcal Z/(\gamma_1\gamma_2)}=\eps^2\,\Theta\,\bigl(1+O(\eps^4)\bigr),
  \]
  on which $\omega$ restricts to $\Gamma\,\ddd X\wedge\ddd Y+\ddd I_E\wedge\ddd\theta$. The kernel of this $2$-form is spanned by $\partial_\theta+X_{I_E}$, where $X_{I_E}$ is the Hamiltonian vector field of $I_E$ for $\Gamma\,\ddd X\wedge\ddd Y$. Since $\theta$ varies by $2\pi\sgn(\Gamma)$ between two crossings of $\{\theta=0\}$, the return map of $\bar H$ is \emph{exactly} the time-$2\pi\sgn(\Gamma)$ map of $I_E$, that is, the time-$\pi\eps^2$ map of
  \[
    G_\eps:=\frac{2\sgn(\Gamma)}{\eps^2}\,I_E=\frac{|\Gamma|}{2\kappa}\,\Theta\,\bigl(1+O(\eps^4)\bigr),
  \]
  because $\sgn(\Gamma)\,\mu=\gamma_1\gamma_2/|\Gamma|=|\Gamma|/(2\kappa)$. The return map of $\bar H+\Rem=H\circ\Phi$ to $\{\theta=0\}$ differs from it by $O(e^{-c/\eps^2})$ on a complex neighbourhood of $K$, for a smaller $c$, by lemma~\ref{lem:analytic}(ii) applied to $\Rem$. Finally let $\Psi_\eps$ map a point of $\{H\circ\Phi=E\}\cap\{\theta=0\}$ to the intersection with $\mathcal A_E$, closest in time, of the orbit of its image by $\Phi$ (sections being identified with open subsets of the $Z$-plane): it conjugates the two return maps, it preserves $\Gamma\,\ddd X\wedge\ddd Y$ because $\Phi$ is symplectic and Poincaré maps between sections preserve the induced forms, and it is $O(\eps^4)$-close to the identity, because $\Phi$ displaces $Z$ by $O(\eps^4)$ and $\theta$ by $O(\eps^2)$, an angle which the flow sweeps in a time $O(\eps^4)$.
\end{proof}

% \begin{remark}[Why $\Theta$ and not $r$]
%   The level curves of $\Theta$ and of $r$ coincide, so the \emph{orbits} of the
%   guiding centre are the level curves of the Robin function. What distinguishes
%   $\Theta$ is the \emph{timing}: the return time $T\propto\rho^2\propto q^{2\kappa}$
%   varies along $D$, and reparametrizes the slow flow. This is the source of the
%   twist, and it is why $\kappa$ --- the same exponent as in the
%   quadruple-collision blow-up --- reappears here.
% \end{remark}

Critical points of the Robin function have been an object of study for a long time (see \cite{Gladiali:2025} and references therein). The proposition below allows us to translate the existence of critical points of $r$ into the existence of periodic orbits. 

Let $J=\left(\begin{smallmatrix}0&1\\-1&0\end{smallmatrix}\right)$.

\begin{proposition}[Continuation of Robin critical points]
  \label{prop:po}
  Let $Z_*\in D$ be a non-degenerate critical point of $r$. Then there is $\eps_0>0$, for every $0 < \eps < \eps_0$, a periodic orbit $\mathcal O_\eps$ depending real-analytically on $\eps$ and tending to $Z_*^2$ as $\eps \to 0$, such that 
  \[Z(t)=Z_\eps+O(\eps^4),
    \quad
    |z_1(t)-z_2(t)|=\eps\,q(Z_*)^{\kappa}\bigl(1+O(\eps^2)\bigr)\]
  and whose period satisfies
  \[T_\eps=\frac{2\pi\eps^2}{|\Gamma|}q(Z_*)^{2\kappa}\bigl(1+O(\eps^2)\bigr).\]
  Moreover $\mathcal O_\eps$ is elliptic if $Z_*$ is a local extremum of $r$ and
  hyperbolic if $Z_*$ is a saddle, with Floquet multipliers
  $\exp\bigl(\pm T_\eps\nu+O(\eps^4)\bigr)$.
  % , where $\pm\nu$ are the eigenvalues of $\tfrac12\Gamma\,dX_r(Z_*)$.
\end{proposition}

\begin{proof}
  Periodic orbits in $\Sig_E$ that wind once in $\theta$ are the fixed points of
  $\mathcal P_\eps$. By \eqref{eq:P},
  \[
    \frac{2\kappa}{\pi\eps^2|\Gamma|}\bigl(\mathcal P_\eps-\mathrm{id}\bigr)
    =X_\Theta+O(\eps^2)\qquad\text{in }C^1(K).
  \]
  The zeros of $X_\Theta=2\kappa\,\Gamma^{-1}\Theta\,\nabla^\perp r$ are the critical points of $r$, and the differential of $X_\Theta$ at $Z_*$, which is $2\kappa\,\Gamma^{-1}\Theta(Z_*)\,J\,D^2r(Z_*)$, is invertible. The implicit function theorem gives a unique zero $Z_\eps$ of the left-hand side near $Z_*$, with $Z_\eps=Z_*+O(\eps^2)$.

  The displacement is in fact $O(\eps^4)$, for the following reason. By proposition~\ref{prop:return}, $\mathcal P_\eps$ is conjugate, by a diffeomorphism $O(\eps^4)$-close to the identity, to a map which is exponentially close to the time-$\pi\eps^2$ map of $G_\eps=\tfrac{|\Gamma|}{2\kappa}\Theta\,(1+O(\eps^4))$: the would-be $O(\eps^2)$ correction, namely $W_2=2\Gamma I\Ree(Pe^{2i\theta})$, is killed by the very first averaging step because it has zero $\theta$-average. The critical point of $G_\eps$ close to $Z_*$ is $O(\eps^4)$-close to it, and so is the corresponding fixed point of $\mathcal P_\eps$. The value of $\rho$ and the period come from \eqref{eq:rho-graph} and from the formula for $T$. The multipliers are those of $D\mathcal P_\eps(Z_\eps) =\exp\bigl(T_\eps\cdot\tfrac12\Gamma JD^2r(Z_*)\bigr)+O(\eps^4)$.
\end{proof}

\begin{proof}[End of proof of theorem~\ref{thm:periodic}]
  By proposition~\ref{prop:Robin}, $\ddd r\neq0$ near $\partial D$, and $\{r\ge a\}$ is, for $a$ negative enough, a compact deformation retract of $D$. Since $r \to -\infty$ on $\partial D$, it attains a maximum, and proposition~\ref{prop:po} turns it into an elliptic family. If moreover $r$ is Morse, denoting by $\mathsf c_i$ the number of its critical points of index $i$, Morse theory on $\{r \ge a\}$ gives $\mathsf c_2-\mathsf c_1+\mathsf c_0=\chi(D)=2-m$. By lemma~\ref{lm:superharmonic}, $r$ is superharmonic, hence has no local minimum: $\mathsf c_0 = 0$. With $\mathsf c_2\ge1$ this yields $\mathsf c_1 = \mathsf c_2 - 2 + m \ge m-1$ saddle points, hence as many hyperbolic families.
\end{proof}

There remains to prove the existence of quasiperiodic solutions. Let $\mathcal E(A)=\tfrac{\Gamma^2}{2}r(A)$ be the slow energy and $\mathsf b=\tfrac{\gamma_1\gamma_2}{2}$, so that the truncated Hamiltonian of \eqref{eq:nf} reads $\mathcal K(A,I)=\mathcal E(A)-\mathsf b\log(2I/\mu)+O(\eps^4)$.

\begin{lemma}
  \label{lem:iso}
  Non-isochronicity is isoenergetic non-degeneracy of the truncated Hamiltonian of \eqref{eq:nf}. 
\end{lemma}

\begin{proof}
  The isoenergetic non-degeneracy of the Hamiltonian \eqref{eq:nf} corresponds to $\mathcal{K}$ having a non-zero bordered hessian. The latter quantity is 
  \[\det\begin{pmatrix}\partial^2\mathcal K&\nabla\mathcal K\\[2pt]\nabla\mathcal K^{\!\top}&0\end{pmatrix}
    =-\frac{\mathsf b}{I^2}\Bigl(\mathcal E'(A)^2+\mathsf b\,\mathcal E''(A)\Bigr).\]
  But $\Theta=e^{2\kappa r}$ gives $\Theta''=2\kappa\Theta\bigl(r''+2\kappa(r')^2\bigr)$,
  while $\mathcal E'^2+\mathsf b\mathcal E''
  =\tfrac{\Gamma^2\gamma_1\gamma_2}{4}\bigl(r''+\tfrac{\Gamma^2}{\gamma_1\gamma_2}(r')^2\bigr)$
  and $2\kappa=\Gamma^2/(\gamma_1\gamma_2)$, whence the condition    
  \[-\frac{\mathsf b}{I^2}\cdot\frac{\Gamma^2\gamma_1\gamma_2}{4}\,
    \frac{\Theta''(A)}{2\kappa\,\Theta(A)} \neq 0,\]
  which is equivalent to our non-isochronicity condition.
\end{proof}

\begin{lemma}
  \label{lem:auto}
  Every sufficiently small annulus of level curves around a non-degenerate local extremum $Z_*$ of $r$ is non-isochronous.
\end{lemma}

By this lemma, the hypothesis of theorem~\ref{thm:qp} holds for every sufficiently small annulus around a non-degenerate extremum $Z_*$ of $r$; combined with proposition~\ref{prop:po}, the tori of theorem~\ref{thm:qp} thus accumulate on the elliptic periodic orbit family $\mathcal O_\eps$ and on the point $(Z_*,Z_*)\in\Delta$.

\begin{proof}
  Let $Z_*$ be a non-degenerate local extremum. As $A\to0$,
  \[r(A)=r(Z_*)\mp\varsigma A+O(A^2),\quad \varsigma=\frac{\sqrt{\det D^2r(Z_*)}}{|\Gamma|}>0,
  \]
  hence $\Theta(A)=\Theta(Z_*)e^{\mp2\kappa\varsigma A}\bigl(1+O(A^2)\bigr)$ and
  \[\Theta''(A)\longrightarrow (2\kappa\varsigma)^2\,\Theta(Z_*)\neq0.\]

  For $r=r(Z_*)-\tfrac12(aX^2+bY^2)$ in suitable coordinates, the level $r=r(Z_*)-\sigma$ encloses area $2\pi\sigma/\sqrt{ab}$, so $A=\Gamma\sigma/\sqrt{ab}$ and $r$ is affine in $A$ to leading order. The exponential of an affine function has non-vanishing second derivative because $\kappa\neq0$ and $\varsigma\neq0$.
\end{proof}

\begin{remark}
  Lemma~\ref{lem:auto} is worth pausing on: the slow flow near a non-degenerate extremum is \emph{isochronous} to leading order ($\mathcal E''(0)=0$), so the guiding-centre motion alone carries no twist. The twist comes entirely from the $\mathcal E'^2$ term, i.e.\ from the variation of the fast period $T\propto q^{2\kappa}$ along the level curves. Had we blown up with any exponent other than $\kappa$, this term would not be the right one.
\end{remark}

\begin{proof}[End of proof of theorem~\ref{thm:qp}]
  Work with the return map $\mathcal P_\eps$ in the coordinates given by the diffeomorphism $\Psi_\eps$ of proposition~\ref{prop:return}, and in the action--angle coordinates $(A,\psi)$ of $G_\eps$ on $\mathcal U$, which are $O(\eps^4)$-close to those of $\Theta$. By proposition~\ref{prop:return},
  \[
    \mathcal P_\eps(A,\psi)=\bigl(A+f,\ \psi+\alpha_\eps(A)+g\bigr),
    \qquad
    \alpha_\eps(A)=-\pi\eps^2\,G_\eps'(A)=-\frac{\pi\eps^2|\Gamma|}{2\kappa}\,\Theta'(A)+O(\eps^6),
  \]
  where $f$ and $g$ are real analytic and $O(e^{-c/\eps^2})$ on a complex neighbourhood of the annulus.  The twist is $\alpha_\eps'=-\tfrac{\pi\eps^2|\Gamma|}{2\kappa}\Theta''+O(\eps^6)$, bounded below in modulus by $\mathrm{const}\cdot\eps^{2}$ by the non-isochronicity assumption on $\mathcal{U}$.  So $\mathcal P_\eps$ is a real analytic area-preserving twist map (exact, by the argument of lemma~\ref{lem:section}), whose perturbation is exponentially small compared with its twist. The invariant curve theorem in the analytic category \cite{Moser:1962,Russmann:2001} provides real analytic invariant curves, the complement of which has a relative measure of the order of the square root of the size of the perturbation, up to a power of $\eps$, hence $O(e^{-c/\eps^2})$ for a smaller $c$ (compare \cite{Neishtadt:1984}); their saturations by the Hamiltonian flow are the tori of (i). Alternatively, (i) follows from the invariant torus theorem (see \cite{Russmann:2001} or \cite[Section~6]{Fejoz:2016:KAM}) applied to the Hamiltonian~\eqref{eq:nf}, given the isoenergetic non-degeneracy proved in lemma~\ref{lem:iso}. 

  % By proposition~\ref{prop:skew}, $\dot\theta\neq0$ on the whole regime, so the section $\mathcal A_E$ is a global section of the flow restricted to $\Sig_E\cap(\mathcal U\times\SS^1_\theta)$: every orbit of that region meets $\{\theta=0\}$ transversally, in forward and backward time, at times uniformly bounded by $2\pi/\inf|\dot\theta| = O(\eps^2)$. Hence the saturation
  % \[
  %   \mathcal{T}:=\bigcup_{t\in\R}\phi^t_{X_H}\bigl(\mathcal{C}\bigr)
  % \]
  % of an invariant circle $\mathcal C\subset\mathcal A_E$ under the flow is a compact invariant surface, and the map $(p,t)\mapsto\phi^t(p)$ induces a diffeomorphism from the mapping torus of $\mathcal{P}_\eps|_{\mathcal{C}}$ onto $\mathcal{T}$; since $\mathcal P_\eps|_{\mathcal C}$ is conjugate to a Diophantine rotation, $\mathcal T$ is a $2$-torus carrying a quasi-periodic flow with two frequencies. It is Lagrangian in $C_2(D)$ because $\omega$ vanishes on it: $\mathcal T$ is swept by the Hamiltonian flow from a curve on which the restriction $\ddd I\wedge \ddd\theta+\ddd A\wedge \ddd\psi$ of $\omega$ vanishes for dimension reasons.

  For (ii), let $\mathcal{I}=\,]A_1,A_2[$ be the range of the slow action over $\mathcal U$. By proposition~\ref{prop:graph} the energy level is a graph over $\mathcal U\times\SS^1_\theta$, so that
  \[\Sig_E\cap(\mathcal U\times\SS^1)\ \cong\ \mathcal{I}_A \times\T^2_{(\psi,\theta)} .\]
  Each of our tori is a graph $A=a(\psi,\theta)$ over $\T^2$, hence separates this thickened torus into two pieces. Two of them therefore bound a compact region
  \[\mathcal{R}=\{(A,\psi,\theta): a_1(\psi,\theta)<A<a_2(\psi,\theta)\},\]
  which is invariant by the intermediate value theorem. On $\mathcal R$ the separation $\rho$ stays of order $\eps$ and $Z$ stays in $\mathcal U\Subset D$, whence the last assertion of (ii).

  Item (iii) is the standard theory of exact area-preserving twist maps of the annulus, applied to $\mathcal P_\eps$ between two consecutive invariant circles: Poincaré--Birkhoff produces periodic orbits of every rational rotation number in between, and Aubry--Mather theory the Mather sets of every irrational one.

  Finally, for (iv), $|z_1-z_2|=\eps\, q^\kappa(1+O(\eps^2))\to0$ uniformly on $\mathcal U$ by \eqref{eq:rho-graph}, while $Z$ stays $O(\eps^{2})$-close to $\mathcal O$. So, as $\eps\to0$, the tori over a fixed level curve $\mathcal O\subset\mathcal U$ converge, for the Hausdorff distance in $\overline{D}\times\overline{D}$, to $\{(Z,Z):Z\in\mathcal O\}\subset\Delta$.
\end{proof}

\begin{remark}
  The confinement statement can sometimes be improved from two tori to one. If the level curve $\mathcal{O}$ is null-homotopic in $D$, it bounds a disc $\Subset D$, and the region inside a single corresponding KAM torus $T$ is a solid torus with compact closure, containing the elliptic periodic orbit of proposition~\ref{prop:po}: this one torus already confines the solutions inside it.

  On the other hand, if $\mathcal{O}$ encircles one or several holes of $D$ --- for instance if it is a level curve close to a component $\partial_j D$ --- then the inside of $\mathcal{O}$ reaches down to the quadruple-collision end $\T^2 \times \,]0,1[$ studied in section~\ref{sec:quadruple}, and is not relatively compact. There, two tori are genuinely needed.
\end{remark}

\section{Near the domain boundary, away from collisions}
\label{sec:near-domain-boundary}

We now turn to the opposite corner of the configuration space, and analyse the regimes in which one or both vortices are close to $\partial D$, but far from any collision. There are three of them, all of low energy: one vortex may sit in a boundary collar while the other stays in the interior (regime~I); the two vortices may sit in collars of \emph{different} boundary components (regime~III); or both may sit in the collar of the \emph{same} component (regime~II). Together with the quadruple collisions of section~\ref{sec:quadruple}, these are the faces and corners of the conjectural low-energy model $\partial(\overline D\times\overline D)\setminus N(\partial\Delta)$ of remark~\ref{rem:lowE}.

The three regimes are not equally difficult to study, and it is worth saying at once where the difficulty lies. Regimes~I and~III are perturbative, and regime~I is in fact an instance of a mechanism we have already met in section~\ref{sec:collision}: whenever the system has a single fast angle, the return map to a section is generated by the fast action itself, read on the energy surface as a function of the slow position. We isolate that mechanism once and for all in proposition~\ref{prop:universal}; it gives, with no sign hypothesis on $\gamma_1,\gamma_2$, periodic orbits at the non-degenerate critical points of an explicit potential, together with a twist condition which is automatic near the extrema of that potential. It also explains, after the fact, why the function $q^{2\kappa}$ governed the interior binary regime.

Regime~I then follows at once, with effective potential $V_j=\frac{\gamma_2^2}{2}r
+\gamma_1\gamma_2C_j$. Regime~III is a genuine two-frequency
near-integrable system, and it is the most classical of the three: its isoenergetic non-degeneracy determinant is
$\mathsf a_1\mathsf a_2(\mathsf a_1+\mathsf a_2)/(A_1^2A_2^2)$ with
$\mathsf a_i=\gamma_i^2/2>0$, hence \emph{never} degenerate.
Regime~II, on the contrary, is not perturbative: the two vortices drift along the same boundary component at different speeds, so they overtake each other and undergo an $O(1)$ encounter once per relative period. What survives of the perturbative picture is this: the free part of the Hamiltonian is \emph{exactly} in action form, the encounter is governed by a translation-invariant, hence integrable, local model, and there is exactly one resonance --- equal drift speeds, which is possible only when $\gamma_1\gamma_2>0$ --- around which a thin near-integrable island is expected to persist.

As in the previous sections, $D\subset\C$ is a bounded finitely connected domain, with $C^\infty$-smooth boundary
$\partial D=\partial_1D\sqcup\dots\sqcup\partial_mD$, $|\partial_jD|=L_j$;
$g(z,\zeta)=-\log|z-\zeta|+h(z,\zeta)$ the hydrodynamic Green function,
$r(z)=h(z,z)$, $q=e^r$; $\omega=\sum_k\gamma_k\,\ddd x_k\wedge \ddd y_k$ and
$\iota_{X_H}\omega=\ddd H$. It is convenient here to use the form
\begin{equation}
  \label{eq:H}
  H=\gamma_1\gamma_2\,g(z_1,z_2)+\frac{\gamma_1^2}{2}r(z_1)+\frac{\gamma_2^2}{2}r(z_2).
\end{equation}
We write $d_k=\dist(z_k,\partial D)$.

For $\zeta\in D$ let $C_l(\zeta)$ denote the constant value of $g(\cdot,\zeta)$
on $\partial_lD$, and $\ra_{lp}:=C_l|_{\partial_pD}$.

\begin{lemma}
  \label{lem:C}
  $C_l$ is the harmonic function on $D$, smooth up to $\overline D$, with
  boundary value $\ra_{lp}$ on $\partial_pD$; $\ra_{lp}=\ra_{pl}$; and
  \begin{equation}
    \label{eq:gq}
    g(z,\zeta)=C_l(z)\qquad\text{for all }z\in D,\ \zeta\in\partial_lD .
  \end{equation}
  If $m=1$ then $C_1\equiv\ra_{11}$ is constant. If $m\ge2$ then $C_l$ is
  non-constant as soon as the $\ra_{lp}$ are not all equal, and then
  $\partial_nC_l\not\equiv0$ on $\partial D$.
\end{lemma}

\begin{proof}
  Symmetry of $g$ gives $C_l(z)=g(\zeta,z)$ for any fixed $\zeta\in\partial_lD$,
  which is \eqref{eq:gq} after exchanging the arguments; harmonicity and
  boundary values follow. Symmetry of $\ra$ is the same exchange applied twice. We abbreviate $\ra_j := \ra_{jj}$, the constant boundary value of $g(\cdot,\cdot)$ on $\partial_jD$ (it vanishes when $m=1$, since then $g = g_0$).
\end{proof}

Let $E\to-\infty$. Then $H\to-\infty$ forces at least one $d_k\to0$, and the
possibilities are:
\begin{center}
\renewcommand{\arraystretch}{1.3}
\begin{tabular}{@{}c>{\raggedright\arraybackslash}p{5.1cm}>{\raggedright\arraybackslash}p{5.4cm}@{}}
\hline
 & configuration & structure \\
\hline
\textbf{I} & $z_1\in$ collar of $\partial_jD$, $z_2\in K\Subset D$ &
  one fast angle; proposition~\ref{prop:universal} applies verbatim \\
\textbf{II} & $z_1,z_2\in$ collar of the \emph{same} $\partial_jD$ &
  two fast angles \emph{plus} an $O(1)$ encounter; not perturbative \\
\textbf{III} & $z_1\in$ collar of $\partial_jD$, $z_2\in$ collar of
  $\partial_kD$, $k\neq j$ & two fast angles, classical KAM \\
\hline
(IV) & $z_1,z_2\to$ same point of $\partial D$ & quadruple collision; next section \\
\hline
\end{tabular}
\end{center}
In the conjectural low-energy model
$\Sig_E\cong\partial(\overline D\times\overline D)\setminus N(\partial\Delta)$
of remark~\ref{rem:lowE}, regime~I is the open face $\partial_jD\times D$,
regime~III the corner $\partial_jD\times\partial_kD$, and regime~II the corner
$\partial_jD\times\partial_jD$ with the diagonal removed --- the removal being
exactly regime~IV.

Fix $j$. For $z$ in a collar of $\partial_jD$ let
\[
  \alpha(z):=\Area\{w:\ q(w)<q(z)\}
\]
be the collar area. By lemma~\ref{lm:asymptotics}, $\alpha=L_jd\,(1+O(d))$ and
\begin{equation}
  \label{eq:RR}
  r(z)=\RR_j\bigl(\alpha(z)\bigr),\qquad
  \RR_j(\alpha)=\log\alpha+\mathsf r_j+O(\alpha),\qquad
  \mathsf r_j=\log\frac{2e^{\ra_j}}{L_j}.
\end{equation}
For a vortex of strength $\gamma$ put
\begin{equation}
  \label{eq:aa}
  \JJ:=\frac{\gamma}{2\pi}\,\alpha,\qquad \varphi\ \text{conjugate},\qquad
  \ddd\JJ\wedge \ddd\varphi=\gamma\,\ddd x\wedge \ddd y .
\end{equation}
Note that $\sgn\JJ=\sgn\gamma$ while $\alpha>0$ always; we use $\alpha$ as the
working variable.

\begin{proposition}[Exact action form of the self-energy]
  \label{prop:exact}
  $\tfrac{\gamma^2}{2}r(z)=\tfrac{\gamma^2}{2}\RR_j(\alpha)$ depends on the
  action alone, \emph{exactly} --- there is no $\varphi$-dependence and no
  expansion. The corresponding frequency is
  \begin{equation}
    \label{eq:freq}
    \dot\varphi=-\frac{2\pi}{\gamma}\,\partial_\alpha\Bigl(\frac{\gamma^2}{2}\RR_j\Bigr)
      =-\frac{\pi\gamma}{\alpha}\bigl(1+O(\alpha)\bigr),
  \end{equation}
  i.e.\ the vortex circulates around $\partial_jD$ with period
  $2L_jd/|\gamma|\,(1+O(d))$, drifting at speed $|\gamma|/2d$.
\end{proposition}

This is worth emphasising: \emph{all} $\varphi$-dependence in the two-vortex
problem, in every boundary regime, comes from the interaction term
$\gamma_1\gamma_2\,g(z_1,z_2)$.

\begin{lemma}[Boundary expansion of $g$]
  \label{lem:screening}
  Write $z_k=p_k+d_kn(p_k)$ with $p_k\in\partial D$ and $n$ the inward normal.
  \begin{enumerate}[label=(\roman*)]
    \item \emph{(One point on the boundary.)} If $p_1\in\partial_jD$ and
      $p_2\in\partial_lD$ with $p_1\neq p_2$, then uniformly for
      $|p_1-p_2|\ge\mathrm{const}>0$,
      \begin{equation}
        \label{eq:exp2}
        g(z_1,z_2)=\ra_{jl}
          +d_1\,\partial_nC_l(p_1)+d_2\,\partial_nC_j(p_2)+O\bigl(d_1^2+d_2^2+d_1d_2\bigr).
      \end{equation}
      \emph{The two terms of order $d$ are single-particle: each depends on one
      vortex only.}
    \item \emph{(Simply connected case.)} If $m=1$ then $C_1$ is constant, the
      order-$d$ terms vanish, and $g-\ra_{11}=O(d_1d_2)$. In the half-plane with
      Dirichlet data and along-wall separation $s$,
      \[
        g=\tfrac12\log\frac{s^2+(d_1+d_2)^2}{s^2+(d_1-d_2)^2}
         =\frac{2d_1d_2}{s^2}+O\Bigl(\frac{d^4}{s^4}\Bigr)\quad (s\gg d),
      \]
      the interaction of two dipoles of moments $\propto\gamma_kd_k$: this is
      \emph{image screening}.
    \item \emph{(One point in the interior.)} If $p_1\in\partial_jD$ and
      $z_2\in K\Subset D$ then, uniformly on $K$,
      \begin{equation}
        \label{eq:exp1}
        g(z_1,z_2)=C_j(z_2)+d_1\,\Pi_j(p_1,z_2)+O(d_1^2),
        \qquad \Pi_j(p,\cdot):=\partial_{n_1}g(z_1,\cdot)\big|_{z_1=p}.
      \end{equation}
      For each $p$, $\Pi_j(p,\cdot)$ is harmonic on $D$ with boundary value
      $\partial_nC_l(p)$ on $\partial_lD$. \emph{This term does couple the two
      vortices.} For Dirichlet data on the disc, $\Pi(p,\zeta)=(1-|\zeta|^2)/|p-\zeta|^2>0$.
  \end{enumerate}
\end{lemma}

\begin{proof}
  All three are Taylor expansions in $d_1,d_2$ of $g$ at boundary points, using
  \eqref{eq:gq}. For (i): $g(z_1,z_2)$ at $d_1=d_2=0$ is $\ra_{jl}$; by
  \eqref{eq:gq}, $\partial_{d_1}g(z_1,z_2)|_{z_2\in\partial_lD}
  =\partial_n\bigl(g(\cdot,q)\bigr)(p_1)=\partial_nC_l(p_1)$, independent of
  $q\in\partial_lD$, and symmetrically for $d_2$. For (ii): when $m=1$,
  \eqref{eq:gq} says $g(\cdot,q)\equiv\ra_{11}$, so $g-\ra_{11}$ vanishes
  identically when either point is on the boundary, hence is divisible by
  $d_1d_2$. (iii) is the definition of $\Pi_j$ together with \eqref{eq:gq} for
  the boundary values.
\end{proof}

\begin{proposition}[Return map generated by the fast action]
  \label{prop:universal}
  Suppose that on some region of phase space
  \[
    H=\Phi(\JJ)+V(w)+W,
    \qquad \Phi(\JJ)=\mathsf a\log\JJ+\mathrm{const}+O(\JJ),\quad \mathsf a\ne0,
  \]
  where $(\JJ,\varphi)$ is an action--angle pair, $w$ is a slow planar variable
  with symplectic form $\gamma_s\,\ddd x\wedge \ddd y$, and $W$ is a remainder.
  Let $\JJ_E(w)$ denote the fast action determined by $H=E$, i.e.
  \begin{equation}
    \label{eq:Jc}
    \JJ_E(w)=\JJ_0\,e^{-V(w)/\mathsf a}\bigl(1+o(1)\bigr),\qquad \JJ_0=\mathrm{const}\cdot e^{E/\mathsf a}.
  \end{equation}
  Then the return time to $\{\varphi=0\}\cap\Sig_E$ is
  $T=2\pi|\JJ_E/\mathsf a|\,(1+o(1))$, and the return map is
  \begin{equation}
    \label{eq:Puniv}
    \mathcal P_E=\mathrm{id}\ -\ \sgn(\mathsf a \JJ_E)\,\frac{2\pi}{\gamma_s}\,\nabla^{\perp}\JJ_E\ +\ \cdots
    \qquad \nabla^\perp=(\partial_Y,-\partial_X).
  \end{equation}
  (The action $\JJ$ carries the sign of the intensity of the fast vortex, and may well be negative; the sign above takes care of it.) Consequently:
  \begin{enumerate}[label=(\alph*)]
    \item \emph{Fixed points.} The fixed points of $\mathcal P_E$ --- i.e.\ the
      periodic orbits of the flow winding once in $\varphi$ --- are the critical
      points of $V$, which are also those of $\JJ_E$. At a non-degenerate one
      the continuation is unique and the orbit is elliptic if $V$ has a local
      extremum there and hyperbolic if a saddle.
    \item \emph{Twist.} Let $A$ be the action of the level curves of $\JJ_E$ for
      $\gamma_s\,\ddd x\wedge \ddd y$. The twist condition for $\mathcal P_E$ is
      \[
        \frac{\ddd^2\JJ_E}{\ddd A^2}\neq0,
      \]
      i.e.\ the Hamiltonian flow of $\JJ_E$ is non-isochronous; and this is
      exactly isoenergetic non-degeneracy of the truncated system.
    \item \emph{Automatic near extrema.} If $w_*$ is a non-degenerate local
      extremum of $V$ then $V$ is affine in $A$ to leading order, hence
      $\JJ_E=\JJ_E(w_*)e^{\mp\varsigma A/\mathsf a}(1+O(A^2))$ with
      $\varsigma\ne0$, and the twist condition holds for all small $A$.
  \end{enumerate}
\end{proposition}

\begin{proof}
  $T=2\pi/|\Phi'(\JJ_E)|=2\pi|\JJ_E/\mathsf a|(1+o(1))$, and the slow drift over
  one return is $T\cdot\gamma_s^{-1}\nabla^\perp V$. Since \eqref{eq:Jc} gives
  $\nabla\JJ_E=-\mathsf a^{-1}\JJ_E\nabla V$, we get
  $\JJ_E\,\nabla^\perp V=-\mathsf a\,\nabla^\perp\JJ_E$, whence
  \[
    T\,\gamma_s^{-1}\nabla^\perp V
    =\frac{2\pi}{|\mathsf a|\gamma_s}\,|\JJ_E|\,\nabla^\perp V
    =-\frac{2\pi\sgn(\mathsf a \JJ_E)}{\gamma_s}\,\nabla^\perp\JJ_E .
  \]
  (a) follows by the implicit function theorem applied to
  $\frac{\gamma_s}{2\pi}(\mathcal P_E-\mathrm{id})=-\sgn(\mathsf a\JJ_E)\nabla^\perp\JJ_E+o(1)$ in
  $C^1$, the linearisation at a non-degenerate critical point being
  $-\sgn(\mathsf a\JJ_E) J\,D^2\JJ_E$, $J=\left(\begin{smallmatrix}0&1\\-1&0\end{smallmatrix}\right)$.
  (b) and (c) are lemmas~\ref{lem:iso} and \ref{lem:auto} verbatim, with $\Theta$ replaced
  by $\JJ_E$. The overall sign only reverses the direction of the induced flow;
  it affects neither the fixed points nor the twist.
\end{proof}

\begin{remark}[Why $q^{2\kappa}$ appeared in section~\ref{sec:collision}]
  Take $\JJ=I=\frac{\mu}{2}\rho^2$, so $\Phi=-\frac{\gamma_1\gamma_2}{2}\log I
  +\mathrm{const}$, $\mathsf a=-\frac{\gamma_1\gamma_2}{2}$,
  $V=\frac{\Gamma^2}{2}r$ and $\gamma_s=\Gamma$. Then \eqref{eq:Jc} reads
  $I_E\propto e^{\Gamma^2r/(\gamma_1\gamma_2)}=q^{2\kappa}$. So the generator
  $\Theta=q^{2\kappa}$ of section~\ref{sec:collision} was simply \emph{the rescaled squared
  separation}, i.e.\ the fast action; and the exponent $2\kappa$ is nothing but
  $-1/\mathsf a$ times the coefficient $\Gamma^2/2$ of $r$ in $V$. The same will happen in Regime~I, where the
  generator will be the rescaled distance of the wall vortex with the boundary.
\end{remark}

We start with regime~I, in which one vortex circulates at the wall while the other sits inside. Let $z_1$ lie in the collar of $\partial_jD$ and $z_2\in K\Subset D$. With
$(\JJ_1,\varphi_1)$ as in \eqref{eq:aa} for $\gamma_1$, \eqref{eq:H} splits as
\begin{equation}
  \label{eq:splitI}
  H=\underbrace{\frac{\gamma_1^2}{2}\RR_j(\alpha_1)}_{\Phi(\JJ_1)}
   \ +\ \underbrace{\frac{\gamma_2^2}{2}r(z_2)+\gamma_1\gamma_2\,C_j(z_2)}_{V_j(z_2)}
   \ +\ \underbrace{\gamma_1\gamma_2\bigl(g(z_1,z_2)-C_j(z_2)\bigr)}_{W},
\end{equation}
in which the first two terms are \emph{exact}: $\Phi$ depends on $\JJ_1$ alone
by Proposition~\ref{prop:exact}, and $V_j$ on $z_2$ alone by definition of
$C_j$. By Lemma~\ref{lem:screening}(iii), uniformly on $K$,
\begin{equation}
  \label{eq:WI}
  W=\gamma_1\gamma_2\,d_1\,\Pi_j(p_1,z_2)+O(d_1^2)=O(d_1).
\end{equation}
The truncated system $\Phi+V_j$ is integrable with one slow degree of freedom,
and the averaging of proposition~\ref{prop:universal} involves the single frequency
\eqref{eq:freq}, so \emph{no small divisors occur at any order}.

\begin{definition}
  \label{def:delta}
  Put $\delta:=e^{2E/\gamma_1^2}$, so $\delta\to0$ as $E\to-\infty$.
\end{definition}

\begin{proposition}
  \label{prop:graphI}
  On $\Sig_E\cap\{z_2\in K\}$, for $\delta$ small,
  \begin{equation}
    \label{eq:alphaI}
    \alpha_1=\frac{L_j}{2e^{\ra_j}}\ \delta\ \Theta_{\mathrm I}(z_2)\bigl(1+O(\delta)\bigr),
    \qquad
    \Theta_{\mathrm I}:=e^{-2V_j/\gamma_1^2}
  \end{equation}
  so that $d_1=O(\delta)$ and $\delta=O(d_1)$ uniformly on $K$, and
  $\Sig_E\cap\{z_2\in K\}\cong K\times S^1$ with coordinates $(z_2,\varphi_1)$.
\end{proposition}

\begin{proof}
  Solve $\Phi=E-V_j+O(\delta)$ using \eqref{eq:RR}. Note
  $\mathsf a=\gamma_1^2/2>0$ here, so $\Theta_{\mathrm I}=e^{-V_j/\mathsf a}$ is
  \eqref{eq:Jc} up to the constant $\JJ_0$: as announced, the generator is the
  rescaled distance $\alpha_1/\delta$ of the wall vortex.
\end{proof}

Proposition~\ref{prop:universal} applies with $\mathsf a=\gamma_1^2/2$,
$\gamma_s=\gamma_2$, $V=V_j$, $\JJ_E\propto\Theta_{\mathrm I}$.

\begin{theorem}[Periodic orbits]
  \label{thm:poI}
  Let $w_*\in D$ be a non-degenerate critical point of
  \[
    V_j=\frac{\gamma_2^2}{2}\,r+\gamma_1\gamma_2\,C_j .
  \]
  Then for $\delta$ small the two-vortex system has a periodic orbit on which
  vortex~$1$ circulates once around $\partial_jD$ at distance
  $d_1=\frac{\delta}{2e^{\ra_j}}\Theta_{\mathrm I}(w_*)(1+O(\delta))$ in period
  $2L_jd_1/|\gamma_1|\,(1+O(\delta))$, while
  \[
    z_2(t)=w_*+O(\delta).
  \]
  It is elliptic if $w_*$ is a local extremum of $V_j$ and hyperbolic if a saddle. As $\delta\to0$ these orbits converge to $\partial_jD\times\{w_*\}$.
\end{theorem}

\begin{corollary}[Unconditional existence and counting]
  \label{cor:countI}
  $V_j$ attains a maximum in $D$. Moreover, if $V_j$ is a Morse function, there is at least one elliptic family, and at least $m-1$ hyperbolic ones when $m\ge2$.
\end{corollary}

\begin{proof}
  $V_j\to-\infty$ at $\partial D$, because $r\to-\infty$ while $C_j$ is bounded; moreover $\nabla V_j\neq0$ near $\partial D$, because $|\nabla r|^{-1}=O(d)$ while $|\nabla C_j|=O(1)$. Hence $V_j$ attains a maximum in $D$, and if $V_j$ is Morse, then $\mathsf c_2-\mathsf c_1+\mathsf c_0=\chi(D)=2-m$ with $\mathsf c_2\ge1$, giving $\mathsf c_1\ge m-1$. So there is always at least one elliptic family, and at least $m-1$ hyperbolic ones when $m\ge2$.
\end{proof}

\begin{theorem}[Invariant tori]
  \label{thm:kamI}
  Let $\mathcal U\Subset D$ be an annulus of closed regular level curves of $V_j$ on which $\ddd^2\Theta_{\mathrm I}/\ddd A^2\neq0$. Then for $\delta$ small the flow on $\Sig_E$ has invariant $2$-tori over $\mathcal U$, of relative measure $1-O(\delta^{M})$ for every $M$. These tori are graphs over $(\psi,\varphi_1)$, which separate $\Sig_E\cap(\mathcal U\times S^1)$, so that orbits between two of them are confined for all time: vortex~$2$ stays between two level curves of $V_j$ and vortex~$1$ keeps circulating at a distance $d_1$ with $d_1=O(\delta)$ and $\delta=O(d_1)$. Between consecutive tori, the Poincaré--Birkhoff theorem and Aubry--Mather theory give periodic orbits and Mather sets of every intermediate rotation number.
\end{theorem}

By Proposition~\ref{prop:universal}(c) the hypothesis holds on every sufficiently small annulus around a non-degenerate extremum of $V_j$.  The proofs follow those of section~\ref{sec:collision}, with $(\Theta,\varepsilon^2,\Gamma)$ replaced by $(\Theta_{\mathrm I},\delta,\gamma_2)$ and the normal form driven by the single frequency \eqref{eq:freq}. There is one difference: the fast angle now runs along $\partial_jD$, which is only assumed to be $C^\infty$, and the width of analyticity of $H$ in that angle shrinks like the distance $d_1$ to the boundary. So the normal form is carried out to a finite, arbitrary order in $C^\ell$ norms, as in proposition~\ref{prop:normalform}, instead of the analytic norms of theorem~\ref{thm:nf}; whence a remainder $O(\delta^M)$ for every $M$, and not an exponentially small one. If $\partial D$ is real analytic, the estimates should hold in analytic norms, uniformly up to the boundary, and one should obtain real analytic tori, the complement of which has an exponentially small measure.

\begin{remark}[The effective potential]
  For $m=1$ (or for Dirichlet data) $C_j$ is constant and $V_j=\frac{\gamma_2^2}{2}r+\mathrm{const}$: the wall vortex drops out at leading order and the interior vortex simply sits at a critical point of the Robin function, regardless of $\gamma_1$. For $m\ge2$ the term $\gamma_1\gamma_2C_j$ is genuinely present, and the effective potential \emph{does} depend on both strengths and on which component the  wall vortex is circling --- $V_j$ and $V_k$ are different functions, with different critical points.
\end{remark}

\begin{remark}[Order of the correction]
  Unlike the interior binary regime of section~\ref{sec:collision}, where the first correction to the
  guiding-centre Hamiltonian was $O(\varepsilon^4)$ because
  $\langle W\rangle_\theta$ vanished to second order, here
  $\langle W\rangle_{\varphi_1}=\gamma_1\gamma_2d_1\langle\Pi_j(\cdot,z_2)\rangle_{\varphi_1}
  +O(d_1^2)$ has no reason to vanish --- for Dirichlet data on the disc
  $\Pi>0$, so the average is strictly positive. Hence the corrections are
  $O(\delta)$ and $w_\delta=w_*+O(\delta)$, not $O(\delta^2)$.
\end{remark}

Next we consider regime~III, in which the two vortices circulate along different boundary components (as said, regime~II, in which they share a component, will be dynamically more complicated). Assume $m\ge2$, $z_1$ in the collar of $\partial_jD$ and $z_2$ in the collar of $\partial_kD$, $k\neq j$. Then $|z_1-z_2|\ge\dist(\partial_jD,\partial_kD)>0$ for all time, so no collision can occur and there is no encounter.

Write $\alpha_i=\lambda a_i$ and $A_i=\gamma_ia_i/2\pi$, so that the rescaling $\sigma_\lambda$ satisfies
\begin{equation}
  \label{eq:confIII}
  \sigma_\lambda^*\omega=\lambda\,\omega_0,\qquad
  \omega_0=\ddd A_1\wedge \ddd\varphi_1+\ddd A_2\wedge \ddd\varphi_2 .
\end{equation}
This is a \emph{conformal} rescaling with a single factor, so after the time change $\ddd t=\lambda\,\ddd\tau$ we obtain a Hamiltonian system for the fixed form $\omega_0$. By \eqref{eq:RR} and Lemma~\ref{lem:screening}(i),
\begin{equation}
  \label{eq:hIII}
  \mathsf h_\lambda:=H-\frac{\gamma_1^2+\gamma_2^2}{2}\log\lambda
   =\underbrace{\mathsf a_1\log A_1+\mathsf a_2\log A_2}_{\mathsf h_0(A)}
    +\mathrm{const}_{jk}+\lambda\,\mathsf u(A,\varphi)+O(\lambda^2),
\end{equation}
with
\[
  \mathsf a_i=\frac{\gamma_i^2}{2}>0,\qquad
  \mathsf u(A,\varphi)=\mathsf u_1(A_1,\varphi_1)+\mathsf u_2(A_2,\varphi_2),
\]
the order-$\lambda$ term being \emph{single-particle} by
Lemma~\ref{lem:screening}(i). The rescaled energy curve
$\{\mathsf h_0=\mathrm{const}\}$ is $A_2=\mathrm{const}\cdot A_1^{-\mathsf a_1/\mathsf a_2}$;
we work on a compact arc of it, i.e.\ neither vortex much closer to its wall
than the other.

\begin{proposition}[Integrable to first order]
  \label{prop:III-first}
  Since $\mathsf u$ contains only the Fourier modes $(k,0)$ and $(0,k)$, whose divisors $k\omega_1$ and $k\omega_2$ are bounded away from zero, one non-resonant averaging step --- \emph{with no small divisors} --- brings \eqref{eq:hIII} to
  \[
    \mathsf h_0(A)+\lambda\,\overline{\mathsf u}(A)+O(\lambda^2),
  \]
  again in action form. Genuine coupling appears only at order $\lambda^2$.
\end{proposition}

\begin{lemma}
  \label{lem:III-non-deg}
  For $\mathsf h_0=\mathsf a_1\log A_1+\mathsf a_2\log A_2$,
  \[
    \det\partial^2\mathsf h_0=\frac{\mathsf a_1\mathsf a_2}{A_1^2A_2^2},
    \qquad
    \det\begin{pmatrix}\partial^2\mathsf h_0&\nabla\mathsf h_0\\[2pt]
      \nabla\mathsf h_0^{\!\top}&0\end{pmatrix}
    =\frac{\mathsf a_1\mathsf a_2(\mathsf a_1+\mathsf a_2)}{A_1^2A_2^2}.
  \]
  With $\mathsf a_i=\gamma_i^2/2$ both are strictly positive, \emph{whatever the
  signs of $\gamma_1$ and $\gamma_2$}: Kolmogorov and isoenergetic
  non-degeneracy hold with no hypothesis at all. Equivalently, on the energy
  curve the frequency ratio is
  $\omega_1/\omega_2\propto A_1^{-1-\mathsf a_1/\mathsf a_2}
  =A_1^{-(\gamma_1^2+\gamma_2^2)/\gamma_2^2}$, strictly monotone.
\end{lemma}

\begin{proof}
  The Hessian is diagonal with entries $-\mathsf a_i/A_i^2$. Expanding the
  bordered determinant along its last row gives
  $\frac{\mathsf a_1^2\mathsf a_2+\mathsf a_1\mathsf a_2^2}{A_1^2A_2^2}$.
\end{proof}

That the sign question evaporates here is structural: the integrable part of
\eqref{eq:hIII} is built from the \emph{self}-energies, which carry $\gamma_i^2$.

\begin{theorem}[Perpetual stability in Regime III]
  \label{thm:III}
  Fix a compact arc of the rescaled energy curve. For $\lambda$ small the flow on $\Sig_E$ carries invariant Lagrangian $2$-tori of relative measure $1-O(\lambda)$ (and $1-O(\lambda^{M})$ for every $M$ after enough averaging steps), graphs over $(\varphi_1,\varphi_2)$. They separate the $3$-dimensional energy surface, so:
  \begin{enumerate}[label=(\roman*)]
  \item the distances $d_1,d_2$ of the two vortices to their respective boundary components remain, for all $t\in\R$, within $O(\lambda)$ of their initial values, for every orbit trapped between two tori;
  \item the two vortices circulate around $\partial_jD$ and $\partial_kD$ forever, with incommensurable frequencies $\omega_i$ such that $\omega_i=O(\lambda^{-1})$;
  \item the tori accumulate, as $\lambda\to0$, on the corner torus $\partial_jD\times\partial_kD$;
  \item between consecutive tori one finds periodic orbits and Aubry--Mather sets.
  \end{enumerate}
\end{theorem}

\begin{proof}[Proof scheme]
  \eqref{eq:confIII}, \eqref{eq:hIII} and Proposition~\ref{prop:III-first} put the system in the standard form $\mathsf h_0(A)+\lambda^2\mathsf w(A,\varphi)$ with $\mathsf h_0$ isoenergetically non-degenerate by Lemma~\ref{lem:III-non-deg}. Apply isoenergetic KAM. Item (i) is the separation property of $2$-tori in a $3$-manifold; (iii) follows from $d_i=O(\lambda)$; (iv) follows from  the Poincaré--Birkhoff theorem and Aubry-Mather theory. 
\end{proof}

\begin{remark}
  Regime~III is the only one of the four in which genuine small divisors appear, and also the only one in which the unperturbed system is a two-frequency integrable system rather than a rotation over a one-degree-of-freedom flow.  It is, correspondingly, the most classical.
\end{remark}

\medskip Eventually we look at regime~II, where both vortices circulate along the same boundary component.  Here the perturbative structure breaks down, and it is worth being precise about how.

Both vortices are in the collar of $\partial_jD$, with actions $\alpha_i=\lambda a_i$ and angles $\varphi_i$. By Proposition~\ref{prop:exact},
\begin{equation}
  \label{eq:splitII}
 H= \frac{\gamma_1^2}{2} \RR_j(\alpha_1)+\frac{\gamma_2^2}{2} \RR_j(\alpha_2) +\gamma_1\gamma_2\,g(z_1,z_2),
\end{equation}
and the first two terms are \emph{exactly} in action form, with no
$\varphi$-dependence whatever. The rescaling \eqref{eq:confIII} and
Lemma~\ref{lem:III-non-deg} apply unchanged to them. Everything hinges on the
third term.

Introduce the resonant coordinates
\begin{equation}
  \label{eq:resonant}
  \psi:=\varphi_1-\varphi_2,\quad \Upsilon:=\varphi_2,
  \qquad\text{conjugate to}\qquad \JJ_1,\quad \mathcal S:=\JJ_1+\JJ_2,
\end{equation}
which is symplectic since
$\JJ_1\ddd\varphi_1+\JJ_2\ddd\varphi_2=\JJ_1\ddd\psi+\mathcal S\,\ddd\Upsilon$. Note
$\mathcal S\propto\gamma_1d_1+\gamma_2d_2+O(\lambda^2)$, the linear impulse of
the pair relative to the wall.

Let $s$ be the arclength separation along $\partial_jD$, so
$s=\frac{L_j}{2\pi}\psi+O(\lambda)$.

\begin{proposition}[Two zones]
  \label{prop:II-zones}
  Uniformly in the collar:
  \begin{enumerate}[label=(\roman*)]
    \item \emph{(Drift zone, $s\gg\lambda$.)} By Lemma~\ref{lem:screening}(i),
      \[
        \gamma_1\gamma_2 g=\gamma_1\gamma_2\ra_{jj}
        +\underbrace{O(\lambda)\ \text{single-particle}}_{\text{removable by one averaging step}}
        +\ O\!\left(\frac{\lambda^2}{s^2}\right),
      \]
      the coupling proper being $O(\lambda^2/s^2)$. For $m=1$ the
      order-$\lambda$ term is absent altogether.
    \item \emph{(Encounter zone, $s=O(\lambda)$.)} There
      $g=O(1)$ --- at $s=0$ exactly $g=\log\frac{d_1+d_2}{|d_1-d_2|}+O(\lambda)$.
      The zone is crossed in angular width $O(\lambda)$ and, since the inverse
      of the relative drift speed
      $\bigl|\frac{\gamma_1}{2d_1}-\frac{\gamma_2}{2d_2}\bigr|$ is $O(\lambda)$,
      in physical time $O(\lambda^2)$.
      During the crossing $|\dot\JJ_i|=O(\lambda^{-1})$, so the excursion of
      $\JJ_i$ is $O(\lambda)$, i.e.\ $O(1)$ \emph{relative}: the encounter is a
      strong interaction.
  \end{enumerate}
  Consequently neither classical KAM (the perturbation is not small) nor plain
  averaging (the encounter is not removable) applies to the region as a whole.
\end{proposition}

\begin{proposition}[Local model]
  \label{prop:II-local}
  In the encounter zone, replace $\partial_jD$ by its tangent line. The
  resulting system --- two vortices $\gamma_1,\gamma_2$ in a half-plane, i.e.\
  four vortices $\gamma_1,\gamma_2,-\gamma_1,-\gamma_2$ symmetric about the
  wall --- is translation invariant, hence conserves $\gamma_1d_1+\gamma_2d_2$
  in addition to the energy, and is Liouville integrable. Equivalently: in the
  coordinates \eqref{eq:resonant}, $H$ depends on $\Upsilon$ only through curvature
  and through the closing-up of $\partial_jD$, at relative order $\lambda$, so
  $\mathcal S$ is conserved up to $O(\lambda)$ and the reduced system in
  $(\JJ_1,\psi)$ has one degree of freedom.

  In particular the encounter defines a \emph{scattering map} between the
  asymptotic states, which preserves both $H$ and $\mathcal S$ up to
  $O(\lambda)$. These two relations determine $(\alpha_1,\alpha_2)$ up to
  finitely many branches, so at leading order an encounter produces either no
  net change of the actions (a pure phase shift) or an exchange between the
  branches; the transient excursion during the crossing is $O(1)$ relative in
  both cases.
\end{proposition}

% This is the wall version of Love's leapfrogging problem: two same-signed
% vortices above a wall are two vortex pairs on a common axis, which leapfrog.
We have not carried out the branch analysis of the scattering map --- that is the classical computation and it is what one would need to make the description of regime~II quantitative. We just notice the similarity with a famous celestial mechanical system: the dynamics of the $2$ satellites of Saturn called Janus and Epimetheus (\cite{Niederman:2020} and references therein).

\begin{remark}
  The right picture of Regime~II is therefore \emph{free drift, punctuated once per relative period by an integrable $O(1)$ scattering event}: a structured, impact-system-like problem, closer to a billiard than to a KAM problem. The net action drift per encounter is $O(\lambda)$, so one expects transport on a time scale of $O(\lambda^{-1})$ encounters; we have not estimated it.
\end{remark}

\begin{proposition}[Equal drift speeds]
  \label{prop:II-res}
  By \eqref{eq:freq}, $\dot\psi=-\pi\bigl(\frac{\gamma_1}{\alpha_1}-\frac{\gamma_2}{\alpha_2}\bigr)(1+O(\lambda))$,
  so the vortices never overtake each other exactly when
  \begin{equation}
    \label{eq:res}
    \frac{\gamma_1}{d_1}=\frac{\gamma_2}{d_2}\qquad\text{i.e.}\qquad
    \frac{d_1}{d_2}=\frac{\gamma_1}{\gamma_2},
  \end{equation}
  which requires $\gamma_1\gamma_2>0$. \emph{If $\gamma_1\gamma_2<0$ there is no
  resonance}: the two vortices counter-drift, encounter each other once per
  relative period, and Proposition~\ref{prop:II-zones} applies to every orbit.
\end{proposition}

\begin{proposition}[Expected structure of the island; sketch]
  \label{prop:II-island}
  Assume $\gamma_1\gamma_2>0$ and work at fixed $\mathcal S$. Near the resonance
  \eqref{eq:res} the reduced Hamiltonian is, in rescaled variables,
  \[
    \mathcal H=F(\JJ_1)+U(\psi;\JJ_1),\qquad F'(\JJ_1^*)=0,\quad (F'')^{-1}=O(1),
    \quad U=O(\lambda^2),
  \]
  a pendulum with libration frequency $O(\lambda)$ and half-width $O(\lambda)$ in
  the action. On a straight periodic wall
  $U\propto\gamma_1\gamma_2\,d_1d_2\,\sin^{-2}(\pi s/L_j)>0$, minimal at
  $s=L_j/2$; so the centre of the island is the \emph{antipodal configuration}:
  two vortices diametrically opposite on their level curves, at heights in the
  ratio $\gamma_1:\gamma_2$, drifting around $\partial_jD$ in formation. The
  symmetry-breaking terms are $O(\lambda)$ \emph{relative} to $U$, i.e.\
  $O(\lambda^3)$, so they are small compared with the island's energy scale
  $O(\lambda^2)$ and the island should persist, with KAM curves inside it.
  Orbits in the island never enter the encounter zone; the separatrix of the
  barrier at $s=0$ sits at energies $O(1)$, hence at action distance $O(1)$, far
  outside.
\end{proposition}

We emphasise that Proposition~\ref{prop:II-island} is a scheme: the resonant
normal form, the persistence of the island under the $O(\lambda^3)$ breaking,
and the continuation of the antipodal relative equilibrium to the curved domain
all remain to be written out.

Three features of this section deserve to be recorded, because they will reappear in the next one. The first is that a single mechanism is at work in several places: proposition~\ref{prop:universal} covers the interior binary regime and regime~I at once, and would cover any other regime with a logarithmic fast action. In each case the generator of the return map is the fast action read on the energy surface --- $q^{2\kappa}=\rho^2/\varepsilon^2$ in section~\ref{sec:collision}, $\alpha_1/\delta$ here --- the periodic orbits sit at the critical points of the slow potential, and the twist is automatic near its extrema.

The second is that the signs of the intensities are a boundary phenomenon. Only the quadruple-collision regime needs $\gamma_1\gamma_2>0$, and it needs it only to make $\kappa$ positive, hence \eqref{eq:zeta} a blow-up; and only the island of regime~II needs it, for the resonance to exist at all. Everywhere else the integrable part of the Hamiltonian is built from self-energies, proportional to $\gamma_i^2$, or from $\Gamma^2$, and the question of signs does not arise. This is the local counterpart of remark~\ref{rem:dichotomy}.

The third is that regime~II is the only non-perturbative one, and that its obstruction is specific and identifiable: two vortices on the same boundary component drift at different speeds, so they must meet. What replaces perturbation theory there is an integrable local model together with a scattering map, and a single thin resonant island in which the encounter never happens.

\section{Near quadruple collisions}
\label{sec:quadruple}

In this section we assume $\gamma_1,\gamma_2>0$ and we study the only part of $\partial C_2(D)$ that the closure of an energy level contains, namely the \emph{quadruple collisions} $z_1=z_2=z_1^*=z_2^*$, where both vortices collide on the boundary, hence simultaneously with their images. We set up a quasi-homogeneous blow-up, in which the relative position is measured in the scale $q^{\kappa}$, $q=e^{r}$ being the harmonic radius and $\kappa=\Gamma^2/(2\gamma_1\gamma_2)$. Three logarithmic divergences then cancel identically and the Hamiltonian extends continuously to the blow-up. We determine the ends of the energy levels and, in several regimes, their diffeomorphism type. Finally we prove uniform anisotropic estimates on dyadic shells and reduce the dynamics near a quadruple collision to a perturbed twist map, which yields invariant tori accumulating on the collision set; as a by-product, quadruple collisions turn out to be neither reached nor even approached.

The notation is that of the previous sections; we recall only that
\begin{equation}
  \label{eq:H-quad}
  H=\gamma_1\gamma_2\bigl(-\log|z_1-z_2|+h(z_1,z_2)\bigr)
   +\tfrac12\bigl(\gamma_1^2 r(z_1)+\gamma_2^2 r(z_2)\bigr)
\end{equation}
on $C_2(D)=D\times D\setminus\Delta$, with $\omega=\sum_{k=1}^2\gamma_k\,\ddd x_k\wedge \ddd y_k$, $\iota_{X_H}\omega=\ddd H$, $r(z)=h(z,z)$ and $q=e^r$. Throughout, $\gamma_1,\gamma_2>0$ and we set
\[
  \Gamma=\gamma_1+\gamma_2,\qquad
  \mu=\frac{\gamma_1\gamma_2}{\Gamma},\qquad
  \beta=\frac{\gamma_1^2+\gamma_2^2}{2\gamma_1\gamma_2}\ \ge 1,\qquad
  \kappa=1+\beta=\frac{\Gamma^2}{2\gamma_1\gamma_2}\ \ge 2 .
\]
The identity $\kappa\gamma_1\gamma_2=\Gamma^2/2$ is used constantly. Note that
$\kappa=1+\frac{\gamma_2}{2\gamma_1}+\frac{\gamma_1}{2\gamma_2}$ and that
$\kappa\ge 2$ with equality iff $\gamma_1=\gamma_2$.

Let
\[
  Z=\frac{\gamma_1z_1+\gamma_2z_2}{\Gamma},\qquad u=z_1-z_2,
  \qquad\text{so}\qquad z_1=Z+\tfrac{\gamma_2}{\Gamma}u,\quad z_2=Z-\tfrac{\gamma_1}{\Gamma}u .
\]
A direct computation with $\tfrac{i}{2}\sum\gamma_k\,\ddd z_k\wedge \ddd\bar z_k$ gives
\begin{equation}
  \label{eq:omega-reduced}
  \omega=\Gamma\,\ddd X\wedge \ddd Y+\mu\,\ddd\Ree u\wedge \ddd\Imm u
   =\Gamma\,\ddd X\wedge \ddd Y+\ddd I\wedge \ddd\theta,
   \qquad u=\rho e^{i\theta},\quad I=\tfrac{\mu}{2}\rho^2 .
\end{equation}
The cross terms cancel precisely because $Z$ is the centre of vorticity.

Let $d(Z)=\mathrm{dist}(Z,\partial D)$. Near $\partial_jD$ one has
$q(Z)=2e^{\ra_j}d(Z)\,(1+O(d))$, where
\begin{equation}
  \label{eq:rho-j}
  \ra_j=\lim_{Z\to\partial_jD}\bigl(r(Z)-\log 2d(Z)\bigr)
\end{equation}
is the constant $\ra_{jj}$ of lemma~\ref{lem:C}. Indeed, in a local conformal chart
$f$ mapping an arc of $\partial_jD$ into $\partial\mathbb H$ one has
$r(Z)=\log\bigl(2\Imm f(Z)\bigr)-\log|f'(Z)|+\tilde k(Z,Z)$ with $\tilde k$
smooth up to the arc and $\tilde k=\ra_j$ on it, and
$2\Imm f/|f'|\to 2d$. In particular $q$ extends to a $C^\infty$ function on a
collar of $\partial D$ with $q=0$ and $\ddd q\ne0$ on $\partial D$.

The energy \eqref{eq:H-quad} tends to $+\infty$ on the interior diagonal and to
$-\infty$ when a single vortex approaches $\partial D$ while the other stays
away. The only part of $\partial C_2(D)$ that the closure of a level set can
meet is therefore
\[
  \partial\Delta:=\{(p,p):p\in\partial D\}\ \cong\ \partial D,
\]
i.e.\ $z_1=z_2=\bar z_1=\bar z_2$, in the sense that both vortices collide
\emph{on} the boundary, hence simultaneously with their images. This is
confirmed quantitatively in Proposition~\ref{prop:ends}.

Write $H=-\gamma_1\gamma_2\log|u|+F(z_1,z_2)$ with
\[
  F(z_1,z_2)=\gamma_1\gamma_2h(z_1,z_2)+\tfrac{\gamma_1^2}{2}r(z_1)+\tfrac{\gamma_2^2}{2}r(z_2),
  \qquad F(Z,Z)=\tfrac{\Gamma^2}{2}\,r(Z),
\]
and define the \emph{interaction remainder}
\begin{equation}
  \label{eq:W-def}
  W:=F(z_1,z_2)-F(Z,Z).
\end{equation}

\begin{lemma}[Vanishing of the linear term; normalization of the quadratic one]
  \label{lem:taylor-bdry}
  As $u\to0$ at fixed $Z\in D$,
  \[
    W=\gamma_1\gamma_2\,\Ree\bigl(P(Z)\,u^2\bigr)+O\bigl(|u|^3\,d(Z)^{-3}\bigr),
    \qquad P(Z)=\partial_{z_1}^2h(z_1,Z)\big|_{z_1=Z}.
  \]
  In particular the term linear in $u$ vanishes identically, and
  $\langle W\rangle_{\theta}=O(|u|^4d^{-4})$, where $\langle\cdot\rangle_\theta$
  denotes the average over $\theta=\arg u$.
\end{lemma}

\begin{proof}
  Put $a=\gamma_2/\Gamma$, $b=-\gamma_1/\Gamma$, so $z_1=Z+au$, $z_2=Z+bu$ and
  $a-b=1$. Symmetry of $h$ gives, at the diagonal point $(Z,Z)$,
  $D_1h=D_2h=\tfrac12 Dr$ and $D^2_{11}h=D^2_{22}h=:A$, $D^2_{12}h=:B$, with
  $D^2r=2A+2B$.

  \emph{Linear term.} Its coefficient is
  \[
    \gamma_1\gamma_2(a+b)\tfrac12 Dr+\tfrac{\gamma_1^2}{2}a\,Dr+\tfrac{\gamma_2^2}{2}b\,Dr
    =\frac{Dr}{2\Gamma}\bigl(\gamma_1\gamma_2^2-\gamma_1^2\gamma_2+\gamma_1^2\gamma_2-\gamma_1\gamma_2^2\bigr)=0 .
  \]

  \emph{Quadratic term.} Collecting $\tfrac12[\,\cdot\,][u,u]$, the coefficient of $A$ is
  \[
    (\gamma_1^2a^2+\gamma_2^2b^2)+\gamma_1\gamma_2(a^2+b^2)
    =\frac{\gamma_1\gamma_2}{\Gamma^2}\bigl(2\gamma_1\gamma_2+\gamma_1^2+\gamma_2^2\bigr)=\gamma_1\gamma_2,
  \]
  and that of $B$ is
  $(\gamma_1^2a^2+\gamma_2^2b^2)+2\gamma_1\gamma_2ab
  =\tfrac{2\gamma_1^2\gamma_2^2}{\Gamma^2}-\tfrac{2\gamma_1^2\gamma_2^2}{\Gamma^2}=0$.
  Thus the quadratic term is $\tfrac{\gamma_1\gamma_2}{2}A[u,u]$. Since
  $h(\cdot,Z)$ is harmonic, $A$ is trace free and
  $\tfrac12A[u,u]=\Ree(P u^2)$ with $P=\partial^2_{z_1}h(z_1,Z)|_{z_1=Z}$.

  \emph{Averages.} $\langle\Ree(Pu^2)\rangle_\theta=0$; the cubic term is odd
  under $u\mapsto-u$, i.e.\ $\theta\mapsto\theta+\pi$, so it also averages to
  zero. Derivative bounds $|D^k F|=O(d^{-k})$ for $k\ge2$ follow from the
  local reflection description of $h$.
\end{proof}

\begin{remark}
  The vanishing of the linear term is the structural reason for using the
  centre of vorticity, and the vanishing of $\langle\cdot\rangle_\theta$ up to
  order $4$ is what makes the averaged system of proposition~\ref{prop:normalform} so close to
  integrable.
\end{remark}

\begin{definition}[Quasi-homogeneous blow-up]
  \label{def:blowup}
  Set
  \begin{equation}
    \label{eq:zeta}
    z_1-z_2=q(Z)^{\kappa}\,\zeta,\qquad \zeta=Re^{i\theta},
    \qquad\text{and}\qquad \xi:=q^{\beta}\zeta=\frac{z_1-z_2}{q(Z)} .
  \end{equation}
\end{definition}

Thus $\xi$ is the separation measured in units of the harmonic radius (a
scale-free, smooth variable), while $\zeta$ is the separation measured in the
quasi-homogeneous scale $q^{\kappa}$.

\begin{proposition}[Exact cancellation of the logarithms]
  \label{prop:cancel}
  In the variables \eqref{eq:zeta},
  \begin{equation}
    \label{eq:H-blown}
    H=-\gamma_1\gamma_2\log R+\Wh,\qquad \Wh:=W .
  \end{equation}
\end{proposition}

\begin{proof}
  $-\gamma_1\gamma_2\log|u|=-\gamma_1\gamma_2\kappa\log q-\gamma_1\gamma_2\log R
  =-\tfrac{\Gamma^2}{2}r(Z)-\gamma_1\gamma_2\log R$, which cancels
  $F(Z,Z)=\tfrac{\Gamma^2}{2}r(Z)$ exactly.
\end{proof}

This is the sense in which $\kappa$ is the \emph{only} admissible exponent: it
is the unique weight for which the self-interaction (Robin) term and the mutual
interaction term cancel, leaving a finite limit.

\begin{lemma}[Structure and size of the interaction near $\partial_jD$]
  \label{lem:structure}
  There are $\bar q>0$, $\delta>0$ and a function $\mathcal V_j$, defined and $C^\infty$ on
  $\{(\varphi,q,\xi)\in\T\times[0,\bar q)\times\C:\ |\xi|<\delta\}$ (real analytic
  if $\partial_jD$ is), such that, with $\varphi$ the boundary-parallel
  coordinate introduced in \eqref{eq:JJ} below,
  \[
    W=\mathcal V_j(\varphi,q,\xi),\qquad
    \mathcal V_j=\frac{\gamma_1\gamma_2\,e^{2\ra_j}}{2}\,
      \Ree\bigl(e^{-2i\phi_j(\varphi)}\xi^2\bigr)+O\bigl(|\xi|^2(q+|\xi|)\bigr),
  \]
  where $e^{i\phi_j(\varphi)}$ is the unit tangent to $\partial_jD$ at the foot
  point. Moreover $\langle\mathcal V_j\rangle_{\arg\xi}=O(|\xi|^4)$, and in the
  blown-up variables
  \begin{equation}
    \label{eq:What-size}
    \Wh=\mathcal V_j\bigl(\varphi,q,q^{\beta}Re^{i\theta}\bigr)
       =\frac{\gamma_1\gamma_2e^{2\ra_j}}{2}\,q^{2\beta}R^2\cos2(\theta-\phi_j)
       +O\bigl(q^{2\beta}R^2(q+q^{\beta}R)\bigr).
  \end{equation}
\end{lemma}

\begin{proof}[Sketch]
  Smoothness of $\mathcal V_j$ in $(\varphi,q,\xi)$ up to $q=0$ follows from the
  local reflection description of $h$ used in \eqref{eq:rho-j}: in a conformal
  chart $f$ onto $\mathbb H$, $h(z_1,z_2)=\log|f(z_1)-\overline{f(z_2)}|
  -\log\frac{|f(z_1)-f(z_2)|}{|z_1-z_2|}+\tilde k(z_1,z_2)$, and $\tilde k$
  extends smoothly across the arc by odd reflection about the constant $\ra_j$.
  Writing $u=q\xi$ and using $q=2e^{\ra_j}d\,(1+O(d))$, every occurrence of $u$ is
  accompanied by a factor $q$, which cancels the blow-up of the derivatives of
  $F$; Hadamard's lemma applied to Lemma~\ref{lem:taylor-bdry} gives
  $\mathcal V_j=\sum_{|\gamma|=2}\xi^\gamma\mathcal V_{j,\gamma}(\varphi,q,\xi)$
  with $\mathcal V_{j,\gamma}$ smooth. For the leading coefficient, one checks
  in the half-plane model ($h=\log|z_1-\bar z_2|$) that
  $\partial^2_{z_1}h(z_1,Z)|_{z_1=Z}=\tfrac18 d^{-2}$, so
  $q^2P\to\tfrac12e^{2\ra_j}e^{-2i\phi_j}$, whence the stated constant.
\end{proof}

\begin{corollary}[Regularity of the blow-up]
  \label{cor:regularity}
  $\Wh$ is a smooth function of $(\varphi,q,q^{\beta},R,\theta)$. Hence it is
  $C^\infty$ up to $q=0$ when $\beta\in\N$ (in particular when
  $\gamma_1=\gamma_2$, where $\beta=1$), it is at least $C^{2}$ in general, and
  it becomes $C^\infty$ (analytic, for analytic $\partial D$) after the ramified
  substitution $q=\sigma^{p}$ whenever $\beta=p'/p\in\mathbb Q$.
\end{corollary}

Near $\partial_jD$ the level curves of $q$ are smooth closed curves. Put
\begin{equation}
  \label{eq:JJ}
  \JJ(Z)=\frac{\Gamma}{2\pi}\,\Area\{w\in D:\ q(w)<q(Z)\},
\end{equation}
the area of the collar between $\partial_jD$ and the level curve through $Z$,
and let $\varphi\in\R/2\pi\Z$ be the conjugate angle, i.e.\
$\ddd\JJ\wedge \ddd\varphi=\Gamma\,\ddd X\wedge \ddd Y$. Then $(\JJ,\varphi,I,\theta)$
are symplectic for \eqref{eq:omega-reduced}, they are smooth up to $\JJ=0$, and
\[
  q=Q_j(\JJ),\qquad Q_j\in C^\infty([0,\JJ_0)),\qquad Q_j(0)=0,
  \qquad Q_j'(0)=\frac{4\pi e^{\ra_j}}{\Gamma L_j}>0 .
\]

\begin{proposition}[Exact separation of the unperturbed part]
  \label{prop:H0}
  On the collar,
  \begin{equation}
    \label{eq:H0}
    H=H_0(\JJ,I)+W,\qquad
    H_0(\JJ,I)=\frac{\Gamma^2}{2}\log Q_j(\JJ)-\frac{\gamma_1\gamma_2}{2}\log\frac{2I}{\mu}.
  \end{equation}
  $H_0$ depends on the actions alone. Its slow part is exactly the Hamiltonian of a \emph{single} vortex of circulation $\Gamma$.
\end{proposition}

\begin{proof}
  Immediate from $H=-\gamma_1\gamma_2\log\rho+\tfrac{\Gamma^2}{2}r(Z)+W$,
  $\rho=(2I/\mu)^{1/2}$ and $r(Z)=\log Q_j(\JJ)$.
\end{proof}

This is the structural reason for choosing the level curves of the harmonic
radius (equivalently of the Robin function) as slow coordinates: any other
choice would make $H_0$ depend on $\varphi$.

Let $\widehat U_j=\{0\le q<\bar q,\ R>0,\ q^{\beta}R<\delta\}$ and
\[
  \mathcal E_j:=\widehat U_j\cap\{q=0\}\ \cong\ \partial_jD\times\C^*
\]
be the \emph{collision manifold}. Rescale time by
\begin{equation}
  \label{eq:timechange}
  \ddd t=q^{2\kappa}\,\ddd\tau,\qquad {}'=\frac{\ddd}{\ddd\tau}.
\end{equation}
Writing $Q'=Q_j'(Q_j^{-1}(q))$ and taking all partial derivatives of $\Wh$ in
the coordinates $(\varphi,q,R,\theta)$:

\begin{proposition}[Blown-up equations]
  \label{prop:field}
  \begin{align}
    q'      &= q^{2\kappa}\,Q'\,\partial_\varphi\Wh, \label{eq:qp}\\
    \varphi'&= -q^{2\kappa-1}Q'\Bigl(\tfrac{\Gamma^2}{2}+q\,\partial_q\Wh-\kappa R\,\partial_R\Wh\Bigr),
              \label{eq:phip}\\
    R'      &= \frac{\partial_\theta\Wh}{\mu R}-\frac{\kappa R}{q}\,q', \label{eq:Rp}\\
    \theta' &= \frac{\Gamma}{R^2}-\frac{\partial_R\Wh}{\mu R}. \label{eq:thetap}
  \end{align}
  All right-hand sides extend continuously (indeed with the regularity of
  Corollary~\ref{cor:regularity}) to $\{q=0\}$, and $2\kappa-1=1+2\beta$.
\end{proposition}

\begin{proof}
  From $\dot\JJ=\partial_\varphi H$ and $q=Q_j(\JJ)$ one gets
  $\dot q=Q'\partial_\varphi\Wh$, hence \eqref{eq:qp}. For \eqref{eq:phip},
  $\dot\varphi=-\partial_\JJ H|_{I}$; at fixed $I$ one has
  $R=(2I/\mu)^{1/2}q^{-\kappa}$, so
  $\partial_\JJ(-\gamma_1\gamma_2\log R)=Q'\gamma_1\gamma_2\kappa/q
  =Q'\Gamma^2/(2q)$ and
  $\partial_q|_I\Wh=\partial_q\Wh-\tfrac{\kappa R}{q}\partial_R\Wh$.
  For \eqref{eq:Rp}, differentiate $I=\tfrac{\mu}{2}q^{2\kappa}R^2$ along the
  flow and use $\dot I=\partial_\theta\Wh$. For \eqref{eq:thetap},
  $\dot\theta=-\partial_IH$ with $\partial_I=\bigl(\mu q^{2\kappa}R\bigr)^{-1}\partial_R$
  and $\partial_I(-\gamma_1\gamma_2\log R)=-\gamma_1\gamma_2/(2I)$.
\end{proof}

Substituting \eqref{eq:What-size} and $\gamma_1\gamma_2/\mu=\Gamma$:
\begin{align}
  R' &= -\Gamma e^{2\ra_j}q^{2\beta}R\,\sin2(\theta-\phi_j)+O\bigl(q^{2\beta}R(q+q^\beta R)\bigr),\label{eq:Rlead}\\
  \theta' &= \frac{\Gamma}{R^2}-\Gamma e^{2\ra_j}q^{2\beta}\cos2(\theta-\phi_j)+\cdots,\label{eq:thetalead}\\
  \varphi' &= -\tfrac{\Gamma^2}{2}Q_j'(0)\,q^{1+2\beta}\bigl(1+O(q)\bigr),\qquad
  q'=O\bigl(q^{2+4\beta}R^2\bigr).\label{eq:phiqlead}
\end{align}

\begin{proposition}[The collision manifold]
  \label{prop:collision}
  $\mathcal E_j$ is invariant. The flow it carries reduces to $\theta'=\Gamma/R^2$, $R'=\varphi'=0$: it is the free rotation of the pair. The face $\mathcal E_j$ is foliated by periodic orbits and the linearization of the rescaled field in the normal direction vanishes: the collision manifold is \emph{normally degenerate}.
\end{proposition}

This is in sharp contrast with the collision manifolds of the Newtonian $N$-body problem, which McGehee's blow-up makes normally hyperbolic, and on which the asymptotic behaviour of collision orbits can therefore be read off the stable and unstable manifolds. Here no such thing is available. 

The degeneracy is moreover intrinsic, not an artefact of the chart: by \eqref{eq:thetalead}--\eqref{eq:phiqlead} the internal rotation period is $O(q^{2\kappa})$ in physical time while the drift period along $\partial_jD$ is $O(q)$, and their ratio tends to $0$ like $q^{1+2\beta}$. The nature of this collision manifold is thus more closely related to the manifold at infinity in the (restricted) $3$-body problem \cite{Capinski:2022}.

\begin{theorem}[Quadruple collisions are not reached in finite time]
  \label{thm:complete}
  Fix $E\in\R$. There is $C_E>0$ such that along any orbit in $\Sig_E:=H^{-1}(E)$ entering the collar,
  \[
    |\dot q|\le C_E\,q^{2\beta}\qquad\text{and hence}\qquad
    q(t)\ \ge\ \bigl(q(0)^{1-2\beta}+(2\beta-1)C_E|t|\bigr)^{-1/(2\beta-1)} .
  \]
  In particular the flow of $X_H$ on $C_2(D)$ is complete.
\end{theorem}

\begin{proof}
  On $\Sig_E$ one has $R=R_E\bigl(1+O(q^{2\beta})\bigr)$ with
  $R_E=e^{-E/(\gamma_1\gamma_2)}$, by Proposition~\ref{prop:cancel}; in
  particular $R$ is bounded. By \eqref{eq:qp} in physical time,
  $\dot q=Q'\partial_\varphi\Wh$, and $\partial_\varphi\Wh=O(q^{2\beta}R^2)$ by
  Lemma~\ref{lem:structure} (differentiating in $\varphi$ costs nothing). Since
  $2\beta\ge2>1$, integrate $\bigl|\tfrac{\ddd}{\ddd t}q^{1-2\beta}\bigr|\le(2\beta-1)C_E$.
  Every other way of leaving $C_2(D)$ sends $H$ to $\pm\infty$.
\end{proof}

\begin{remark}
  The bound is only polynomial. Theorem~\ref{thm:kam} upgrades it to
  $\inf_{t\in\R}q(t)>0$.
\end{remark}

\begin{proposition}[The ends]
  \label{prop:ends}
  For every $E\in\R$ there is a compact $K\subset\Sig_E$ with
  \[
    \Sig_E\setminus K\ \cong\ \bigsqcup_{j=1}^{m}\ \T^2\times(0,1),
  \]
  the torus factor being coordinatized by $(\varphi,\theta)$ and the interval by
  $q$. On the $j$-th end,
  \begin{equation}
    \label{eq:level-graph}
    |z_1-z_2|=R_E\,q(Z)^{\kappa}\bigl(1+O(q^{2\beta})\bigr),\qquad R_E=e^{-E/(\gamma_1\gamma_2)} .
  \end{equation}
  The closure of $\Sig_E$ in the blow-up meets $\{q=0\}$ transversally, in the
  disjoint tori $\mathcal C_{j,E}=\{q=0,\ R=R_E\}\cong\partial_jD\times S^1$.
\end{proposition}

\begin{proof}
  $\partial_R H=-\gamma_1\gamma_2/R+O(q^{2\beta}R)\ne0$ on $\widehat U_j$, so
  $\Sig_E$ is the graph $R=R_Ee^{\Wh/(\gamma_1\gamma_2)}$; expand. Since
  $H\to+\infty$ on the interior diagonal and $H\to-\infty$ elsewhere on
  $\partial C_2(D)$, the complement of the collars is compact in $\Sig_E$.
\end{proof}

\begin{proposition}[No change between critical values]
  \label{prop:isotopy}
  If $[E,E']$ contains no critical value of $H$ then $\Sig_E\cong\Sig_{E'}$, and
  $H^{-1}([E,E'])\cong\Sig_E\times[E,E']$.
\end{proposition}

\begin{proof}
  On each $\widehat U_j$ take $V_j:=\bigl(R\,\partial_RH\bigr)^{-1}R\,\partial_R$;
  it satisfies $V_jH=1$, preserves $(q,\varphi,\theta)$ and is therefore
  complete on $H^{-1}([E,E'])\cap\widehat U_j$. On the compact part take
  $\nabla H/|\nabla H|^2$. Glue with a partition of unity depending on $q$
  alone; the resulting $V$ still satisfies $VH=1$ and is complete because it is
  tangent to the level sets of $q$ near the ends.
\end{proof}

\begin{proposition}[High energy, arbitrary $D$]
  \label{prop:highE}
  There is $E_*$ such that for all $E>E_*$ the map $(z_1,z_2)\mapsto(Z,\zeta)$
  identifies $\Sig_E$ with a graph $R=R_E(Z,\theta)$ over $D\times S^1$. Hence
  \[
    \Sig_E\ \cong\ D\times S^1 \qquad (E>E_*),
  \]
  and the blown-up closure is $\overline{D}\times S^1$.
\end{proposition}

\begin{proof}[Sketch]
  Two ingredients. \emph{(i)} $\sup\{H:\ |\xi|\ge\delta_0\}<\infty$ for each
  $\delta_0>0$: the closure of $\{|\xi|\ge\delta_0\}$ in $\overline D\times\overline D$
  meets $\Delta$ only along $\partial\Delta$, $H$ extends continuously with
  values in $[-\infty,\infty)$ off $\partial\Delta$, and near $\partial\Delta$
  one uses $H=-\gamma_1\gamma_2\log|\xi|+\tfrac{\gamma_1^2+\gamma_2^2}{2}\log q+W$
  together with Lemma~\ref{lem:structure}. \emph{(ii)} On $\{|\xi|<\delta_0\}$,
  $R\,\partial_RH=-\gamma_1\gamma_2+O(q^{2\beta}R^2)<0$ uniformly on $D$, again
  by Lemma~\ref{lem:structure}. Choosing $E_*$ as the supremum in (i) forces
  $\{H>E_*\}\subset\{|\xi|<\delta_0\}$, where the level is a radial graph in
  $\zeta$.
\end{proof}

\begin{proposition}[The disc: all levels are open solid tori]
  \label{prop:disc}
  Let $D=\Dd$, so that $h(z,w)=\log|1-z\bar w|$ and $r(z)=\log(1-|z|^2)$. Let
  $\mathsf E=\sum_k(x_k\partial_{x_k}+y_k\partial_{y_k})$ be the Euler field of $\C^2$.
  Then
  \[
    \mathsf EH=-\gamma_1\gamma_2\,\Ree\frac{1+z_1\bar z_2}{1-z_1\bar z_2}
       -\sum_{k}\gamma_k^2\,\frac{|z_k|^2}{1-|z_k|^2}\ <\ 0
    \qquad\text{on }C_2(\Dd).
  \]
  Consequently $H$ has no critical point, $\mathsf E/2$ is a Liouville field transverse
  to every level, and for every $E\in\R$
  \[
    \Sig_E\ \cong\ S^3\setminus(\text{unknot})\ \cong\ \Dd\times S^1 .
  \]
\end{proposition}

\begin{proof}
  $\mathsf E(-\log|z_1-z_2|)=-1$, $\mathsf E\log|1-z_1\bar z_2|=\Ree\frac{-2z_1\bar z_2}{1-z_1\bar z_2}$
  and $\mathsf E\log(1-|z_k|^2)=-2|z_k|^2/(1-|z_k|^2)$; combine the first two using
  $-1-\Ree\frac{2w}{1-w}=-\Ree\frac{1+w}{1-w}$, which is $<0$ for $|w|<1$ by
  the Herglotz representation, applied with $w=z_1\bar z_2$. Along the ray
  $t\mapsto(tu_1,tu_2)$, $t\in(0,t_{\max}(u))$, $t_{\max}=1/\max_k|u_k|$, the
  function $H$ decreases strictly from $+\infty$ to $-\infty$; hence $\Sig_E$
  is a radial graph over the space of rays, which is
  $S^3\setminus K$ with $K=S^3\cap\Delta$ a great circle, i.e.\ an unknot.
\end{proof}

\begin{remark}[Nearly round domains]
  \label{rem:nearlyround}
  Proposition~\ref{prop:disc} extends to an open set of domains. Let $D$ be simply connected, with the origin in $D$, and let $\psi=\gamma_1g_{z_1}+\gamma_2g_{z_2}$. Since $\zeta\mapsto\zeta/t$ maps $D$ conformally onto $D/t$, corollary~\ref{cor:conformal-r} gives $H_D(tz_1,tz_2)=H_{D/t}(z_1,z_2)+\frac{\gamma_1^2+\gamma_2^2}{2}\log t$; differentiating at $t=1$ and applying Hadamard's formula (proof of lemma~\ref{lm:generic-nondeg}; here $g=g_0$ since $D$ is simply connected, and the sign of the formula is that of the monotonicity of $g_0$ with respect to the domain \cite[Section~4.4]{Ransford:1995}) to the normal variation $\delta n=\zeta\cdot n$ of $\partial D$ yields, for arbitrary intensities,
  \[
    \mathsf EH=\frac{\gamma_1^2+\gamma_2^2}{2}-\frac{1}{4\pi}\oint_{\partial D}(\partial_n\psi)^2\,(\zeta\cdot n)\,\ddd s .
  \]
  Now assume $\gamma_1,\gamma_2>0$ and $D$ strictly star-shaped, with boundary $\rho=e^{\ell(\phi)}$ in polar coordinates. The flux of $\psi$ is $\oint\partial_n\psi\,\ddd s=-2\pi\Gamma$, and $\oint\ddd s/(\zeta\cdot n)=2\pi+\int_0^{2\pi}\ell'^2\ddd\phi$, so that the Cauchy--Schwarz inequality gives
  \[
    \mathsf EH\leq\frac{\gamma_1^2+\gamma_2^2}{2}-\frac{\pi\Gamma^2}{2\pi+\int_0^{2\pi}\ell'^2\ddd\phi},
  \]
  which is negative as soon as $\int_0^{2\pi}\ell'^2\ddd\phi<2\pi/\beta=4\pi\gamma_1\gamma_2/(\gamma_1^2+\gamma_2^2)$. Under this condition, which is $C^1$-open and holds for round discs centred at the origin, the proof of proposition~\ref{prop:disc} applies verbatim: $H$ has no critical point, $\frac12\mathsf E$ is a Liouville field transverse to every level, which answers question~\ref{qu:contacttype} for these domains, and $\Sig_E\cong S^3\setminus(\text{unknot})$ for every $E$. The condition is only sufficient, and for $\gamma_1\gamma_2<0$ the argument gives nothing, since then $\Gamma^2<\gamma_1^2+\gamma_2^2$.
\end{remark}

\begin{remark}[Consistency and the shape of the closure]
  Proposition~\ref{prop:disc} is consistent with Proposition~\ref{prop:highE}
  since $\Dd\times S^1\cong S^3\setminus\text{unknot}$. The closure of $\Sig_E$
  in $\C^2$ is a topological $S^3$ obtained by adding $K$ back, but the radial
  function is not Lipschitz there: by \eqref{eq:level-graph} it behaves like
  $t_{\max}(u)-C|u_1-u_2|^{1/\kappa}$, and $1/\kappa\le\tfrac12$. The blow-up of
  definition~\ref{def:blowup} is exactly what removes this cusp.
\end{remark}

\begin{remark}[Low energy and multiple connectivity: open]
  \label{rem:lowE}
  For $E\ll0$ one expects
  $\Sig_E\cong\partial(\overline D\times\overline D)\setminus N(\partial\Delta)$
  (corners smoothed), by pushing with the commuting fields $d_k\partial_{d_k}$
  and $-R\partial_R$, all of which increase $H$. For $D=\Dd$ this gives
  $S^3\setminus N(\text{unknot})$, in agreement with
  Proposition~\ref{prop:disc}. For $m\ge2$ the diffeomorphism type of $\Sig_E$
  changes with $E$ if and only if $H$ has critical points, i.e.\ if and only if
  the pair has equilibria; we do not settle this here.
\end{remark}

From now on we work on the $j$-th collar and drop the index $j$. We
fix $E\in\R$, write $R_E=e^{-E/(\gamma_1\gamma_2)}$, and assume $\partial D$ of
class $C^\infty$.

\begin{definition}
  \label{def:shellnorm}
  For $\lambda>0$ small let
  \[
    \sigma_\lambda(J,\varphi,S,\theta):=(\lambda J,\ \varphi,\ \lambda^{2\kappa}S,\ \theta),
    \qquad
    \mathcal D:=[1,2]\times\T\times[S_-,S_+]\times\T ,
  \]
  where $0<S_-<\tfrac{\mu}{2}Q'(0)^{2\kappa}R_E^2$ and
  $S_+>\tfrac{\mu}{2}Q'(0)^{2\kappa}2^{2\kappa}R_E^2$, so that
  $\sigma_\lambda(\mathcal D)$ is a neighbourhood in phase space of the shell
  $\mathcal S_\lambda:=\Sig_E\cap\{\lambda\le\JJ\le2\lambda\}$. For a function
  $f$ on $\sigma_\lambda(\mathcal D)$ put
  \[
    \|f\|_{\lambda,\ell}:=\bigl\|f\circ\sigma_\lambda\bigr\|_{C^\ell(\mathcal D)}
    =\sum_{|\alpha|\le\ell}\Bigl\|\lambda^{\alpha_1}\lambda^{2\kappa\alpha_3}\,
      \partial_\JJ^{\alpha_1}\partial_\varphi^{\alpha_2}\partial_I^{\alpha_3}\partial_\theta^{\alpha_4}f\Bigr\|_{C^0}.
  \]
\end{definition}

Each $\partial_\JJ$ is thus weighted by $\lambda$ and each $\partial_I$ by
$\lambda^{2\kappa}$, which are the natural scales of the two blocks.

\begin{lemma}[Poisson bracket]
  \label{lem:bracket}
  $\|\{f,g\}\|_{\lambda,\ell}\le C_\ell\,\lambda^{-2\kappa}\|f\|_{\lambda,\ell+1}\|g\|_{\lambda,\ell+1}$.
\end{lemma}

\begin{proof}
  $\{f,g\}\circ\sigma_\lambda=\lambda^{-1}\bigl(\partial_J\tilde f\partial_\varphi\tilde g
  -\partial_\varphi\tilde f\partial_J\tilde g\bigr)
  +\lambda^{-2\kappa}\bigl(\partial_S\tilde f\partial_\theta\tilde g
  -\partial_\theta\tilde f\partial_S\tilde g\bigr)$, and $\lambda^{-1}\le\lambda^{-2\kappa}$.
\end{proof}

\begin{proposition}[The nice limit, quantitatively]
  \label{prop:nicelimit}
  For every $\ell\in\N$ there is $C_\ell$ with
  \begin{equation}
    \label{eq:nicelimit}
    \bigl\|H\circ\sigma_\lambda-\mathsf h_0\bigr\|_{C^\ell(\mathcal D)}\le C_\ell\,\lambda,
    \qquad
    \mathsf h_0(J,S):=\frac{\Gamma^2}{2}\log\bigl(Q'(0)J\bigr)-\frac{\gamma_1\gamma_2}{2}\log\frac{2S}{\mu},
  \end{equation}
  with \emph{no additive renormalization}. Moreover the $\theta$-dependent part
  is much smaller:
  \begin{equation}
    \label{eq:Wsmall}
    \|W\|_{\lambda,\ell}\le C_\ell\,\lambda^{2\beta},\qquad
    \|\langle W\rangle_\theta\|_{\lambda,\ell}\le C_\ell\,\lambda^{4\beta}.
  \end{equation}
\end{proposition}

\begin{proof}
  \emph{(a) Absence of renormalization.} From \eqref{eq:H0},
  \[
    H_0\circ\sigma_\lambda
    =\frac{\Gamma^2}{2}\log Q(\lambda J)-\frac{\gamma_1\gamma_2}{2}\log\frac{2\lambda^{2\kappa}S}{\mu},
  \]
  whose $\log\lambda$ terms have total coefficient
  $\tfrac{\Gamma^2}{2}-\kappa\gamma_1\gamma_2=0$. Since
  $Q(\lambda J)=Q'(0)\lambda J\,(1+O(\lambda))$ in $C^\ell([1,2])$, this gives
  \eqref{eq:nicelimit} for $H_0$.

  \emph{(b) Size of $W$.} By Lemma~\ref{lem:structure},
  $W\circ\sigma_\lambda=\mathcal V\bigl(\varphi,q,q^{\beta}Re^{i\theta}\bigr)$ with
  $q=Q(\lambda J)=Q'(0)\lambda J\,(1+O(\lambda))$ and $R=R(J,S,\lambda)$ given by
  $R^2=2S\bigl(\mu\,(Q(\lambda J)/\lambda)^{2\kappa}\bigr)^{-1}$, a function with
  $C^\ell(\mathcal D)$ norms bounded uniformly in $\lambda$ and bounded away
  from $0$. Since $\mathcal V$ vanishes to order $2$ in its last argument,
  $W\circ\sigma_\lambda=O(q^{2\beta})=O(\lambda^{2\beta})$; each $\partial_J$
  hitting $q^{2\beta}=Q(\lambda J)^{2\beta}$ produces
  $\lambda\cdot O(\lambda^{2\beta-1})=O(\lambda^{2\beta})$, and $\partial_S,
  \partial_\varphi,\partial_\theta$ cost nothing. This is exactly the point of
  the weights in Definition~\ref{def:shellnorm}.

  \emph{(c) The average.} By Lemma~\ref{lem:taylor-bdry} the quadratic and cubic
  terms of $\mathcal V$ have vanishing $\theta$-average, so
  $\langle W\rangle_\theta=O(|\xi|^4)=O(\lambda^{4\beta})$, with the same
  derivative bookkeeping.
\end{proof}

\begin{remark}
  $\mathsf h_0$ is precisely the limiting Hamiltonian
  $\tfrac{\Gamma^2}{2}\log q-\gamma_1\gamma_2\log|z_1-z_2|$ read in the
  quasi-homogeneous scale, and the level $\mathsf h_0=E$ is $R=R_E$, in
  agreement with \eqref{eq:level-graph}.
\end{remark}

Set $\omega_\varphi^0=-\partial_\JJ H_0$, $\omega_\theta^0=-\partial_IH_0$. Then on
$\mathcal S_\lambda$,
\begin{equation}
  \label{eq:freqs}
  \omega_\varphi^0=-\frac{\Gamma^2}{2}\frac{Q'(\JJ)}{Q(\JJ)}
    =-\frac{\Gamma^2}{2\lambda J}\bigl(1+O(\lambda)\bigr),
  \qquad
  \omega_\theta^0=\frac{\gamma_1\gamma_2}{2I}=\frac{\gamma_1\gamma_2}{2\lambda^{2\kappa}S},
  \qquad
  \frac{\omega^0_\varphi}{\omega^0_\theta}
    =-\frac{\Gamma^2S}{\gamma_1\gamma_2J}\,\lambda^{2\kappa-1}\bigl(1+O(\lambda)\bigr),
\end{equation}
with $2\kappa-1=1+2\beta\ge3$. We write
\begin{equation}
  \label{eq:eps}
  \eps:=\lambda^{2\kappa-1}
\end{equation}
for the slow/fast ratio. Note that $\lambda^{2\beta}=\eps^{2\beta/(1+2\beta)}\gg\eps$:
the $\theta$-dependence is small, but \emph{less} small than the frequency ratio.

\begin{lemma}[Isoenergetic non-degeneracy]
  \label{lem:isoenergetic}
  For $\mathcal K(A,I)=\mathsf a\log A-\mathsf b\log I$,
  \[
    \det\begin{pmatrix}\partial^2\mathcal K&\nabla\mathcal K\\[2pt]\nabla\mathcal K^{\!\top}&0\end{pmatrix}
    =-\frac{\mathsf a\,\mathsf b\,(\mathsf a-\mathsf b)}{A^2I^2} .
  \]
  With $\mathsf a=\Gamma^2/2$ and $\mathsf b=\gamma_1\gamma_2/2$ one has
  $\mathsf a-\mathsf b=\tfrac12(\gamma_1^2+\gamma_1\gamma_2+\gamma_2^2)>0$, so
  the determinant is negative. Equivalently, on $\Sig_E$ the ratio
  $\omega_\varphi/\omega_\theta$ is a strictly monotone function of the action.
\end{lemma}

\begin{proof}
  Expand the bordered determinant along the last row. For the equivalent form,
  substitute $I=\tfrac\mu2Q(A)^{2\kappa}R_E^2$ into \eqref{eq:freqs} to get
  \[
    \frac{\omega_\varphi}{\omega_\theta}=-\frac{\Gamma R_E^2}{2}\,Q'(A)\,Q(A)^{2\kappa-1}\bigl(1+o(1)\bigr),
  \]
  whose $A$-derivative is
  $-\tfrac{\Gamma R_E^2}{2}(2\kappa-1)Q'(0)^{2\kappa}A^{2\kappa-2}\bigl(1+O(A)\bigr)\neq0$;
  in particular the ratio is a strictly monotone function of $A$.
\end{proof}

Only the \emph{fast} frequency is used, so there are no small divisors and no
Diophantine conditions at this stage.

\begin{proposition}[Normal form on a shell]
  \label{prop:normalform}
  Fix $N\in\N$ and $\ell\in\N$. For $\lambda$ small there is an exact symplectic
  change of variables $\Phi_N$, defined on $\sigma_\lambda(\mathcal D')$ for a
  slightly smaller $\mathcal D'\Subset\mathcal D$, such that
  \[
    H\circ\Phi_N=H_0(\JJ,I)+Z_N(\JJ,\varphi,I)+\Rem_N ,
  \]
  with
  \[
    \|Z_N\|_{\lambda,\ell}\le C\,\lambda^{4\beta},
    \qquad
    \|\Rem_N\|_{\lambda,\ell}\le C_{N,\ell}\,\lambda^{2(N+1)\beta}.
  \]
  In the rescaled variables $\Phi_N$ is $O(\lambda^{2\beta})$-close to the
  identity; more precisely it displaces $(J,\varphi)$ by
  $O(\lambda^{2\kappa+2\beta-1})=O(\lambda^{1+4\beta})$ and $(S,\theta)$ by
  $O(\lambda^{2\beta})$.
\end{proposition}

\begin{proof}
  Let $\chi$ solve $\omega_\theta^0\,\partial_\theta\chi=-(W-\langle W\rangle_\theta)$,
  i.e.\ $\chi=-\bigl(\omega_\theta^0\bigr)^{-1}\partial_\theta^{-1}\bigl(W-\langle W\rangle_\theta\bigr)$;
  this is an antiderivative in $\theta$ divided by $\omega^0_\theta(I)$, with no
  resonance condition. Since $(\omega_\theta^0)^{-1}=2I/(\gamma_1\gamma_2)$ has
  $\|\cdot\|_{\lambda,\ell}\lesssim\lambda^{2\kappa}$, \eqref{eq:Wsmall} gives
  \[
    \|\chi\|_{\lambda,\ell}\lesssim\lambda^{2\kappa+2\beta}.
  \]
  Let $\Phi^1_\chi$ be the time-one map of $X_\chi$. Then
  \[
    H\circ\Phi^1_\chi=H_0+\langle W\rangle_\theta+E_1+E_2,\qquad
    E_1=\omega^0_\varphi\,\partial_\varphi\chi,
  \]
  $E_2$ collecting the second-order terms of the Lie series. We estimate
  \[
    \|E_1\|_{\lambda,\ell}\lesssim\lambda^{-1}\lambda^{2\kappa+2\beta}=\lambda^{2\beta}\eps,
    \qquad
    \|E_2\|_{\lambda,\ell}\lesssim\lambda^{-2\kappa}\|W\|_{\lambda,\ell+1}\|\chi\|_{\lambda,\ell+1}
    \lesssim\lambda^{4\beta},
  \]
  the second by Lemma~\ref{lem:bracket}. Since $2\kappa-1=1+2\beta>2\beta$, the
  new perturbation is $O(\lambda^{4\beta})$: each step gains a factor
  $\lambda^{2\beta}$. Iterating $N$ times and using
  $\|\langle W\rangle_\theta\|_{\lambda,\ell}\lesssim\lambda^{4\beta}$ from
  \eqref{eq:Wsmall} gives the statement. The displacement estimates follow from
  $\Delta\JJ=\partial_\varphi\chi$, $\Delta\varphi=-\partial_\JJ\chi$,
  $\Delta I=\partial_\theta\chi$, $\Delta\theta=-\partial_I\chi$, rescaled by
  $\lambda^{-1}$ and $\lambda^{-2\kappa}$ respectively.
\end{proof}

\begin{proposition}[Integration of the truncated system]
  \label{prop:integrable}
  $\bar H:=H_0+Z_N$ has $I$ as an integral. There are action--angle variables
  $(A,\psi)$ for the slow block, depending on $I$ as a parameter, with
  $A=\JJ+O(\lambda^{1+4\beta})$, such that
  \[
    \bar H=\mathcal K_N(A,I),\qquad
    \mathcal K_N=\frac{\Gamma^2}{2}\log Q(A)-\frac{\gamma_1\gamma_2}{2}\log\frac{2I}{\mu}
      +O_{\lambda,\ell}(\lambda^{4\beta}) .
  \]
  In particular Lemma~\ref{lem:isoenergetic} applies to $\mathcal K_N$ up to
  $O(\lambda^{4\beta})$, and the frequency ratio on $\Sig_E$ is
  \[
    \frac{\Omega_\psi}{\Omega_\theta}=-\frac{\Gamma R_E^2}{2}\,Q'(A)Q(A)^{2\kappa-1}
      \bigl(1+O(\lambda^{4\beta})+O(\lambda)\bigr).
  \]
\end{proposition}

\begin{proof}
  $\theta$ is cyclic for $\bar H$. At fixed $I$, $|\partial_\JJ\bar H|^{-1}=O(\lambda)$,
  so $\partial_\JJ\bar H$ dominates $\partial_\varphi Z_N=O(\lambda^{4\beta})$, so the level curves of
  $\bar H$ in $(\JJ,\varphi)$ are graphs over $\varphi$ at distance
  $O(\lambda\cdot\lambda^{4\beta})$ from $\{\JJ=\mathrm{const}\}$; apply
  Liouville--Arnold with parameter and set $A=\tfrac{1}{2\pi}\oint\JJ\,\ddd\varphi$.
\end{proof}

\begin{lemma}[A global section of the end]
  \label{lem:section}
  On the collar, $\dot\varphi=-\partial_\JJ H
  =-\tfrac{\Gamma^2}{2}\tfrac{Q'}{Q}\bigl(1+O(q^{2\beta})\bigr)<0$. Hence
  \[
    \mathcal A_E:=\Sig_E\cap\{\varphi=0\}
  \]
  is a global section of the flow restricted to the end, an annulus with
  coordinates $(I,\theta)$, and the return map $\mathcal P$ is exact symplectic for
  $\ddd I\wedge \ddd\theta$.
\end{lemma}

\begin{proof}
  Monotonicity of $\varphi$ is clear from \eqref{eq:H0} and \eqref{eq:Wsmall}.
  On $\{\varphi=0\}\cap\Sig_E$ the restriction of $\omega=\ddd\JJ\wedge \ddd\varphi
  +\ddd I\wedge \ddd\theta$ is $\ddd I\wedge \ddd\theta$, and
  $(\JJ,\theta)\mapsto(I,\theta)$ is a diffeomorphism of the section because
  $\partial_\JJ I|_{H=E}=-\partial_\JJ H/\partial_IH\ne0$. Return maps of
  Hamiltonian flows on a regular level preserve the induced form; exactness
  follows from the usual generating-function (least action) argument, with
  primitive $p\mapsto\int_0^{T(p)}\bigl(\lambda(X_H)-H\bigr)\ddd t+E\,T(p)$ for
  $\lambda=\JJ\,\ddd\varphi+I\,\ddd\theta$.
\end{proof}

\begin{proposition}[Return time and rotation number]
  \label{prop:twist}
  On $\mathcal S_\lambda$ the return time and rotation number of $\mathcal P$ satisfy
  \[
    T=\frac{4\pi\JJ}{\Gamma^2}\bigl(1+O(\lambda)\bigr),
    \qquad
    \varrho(\tilde I)=\eps^{-1}\varsigma(\tilde I)\bigl(1+O(\lambda)+O(\lambda^{4\beta})\bigr),
  \]
  where $\tilde I=I/\lambda^{2\kappa}$, $\eps=\lambda^{2\kappa-1}$, and
  \[
    \varsigma(\tilde I)=\frac{2}{\Gamma R_E^2\,Q'(0)^{2\kappa}}\;a(\tilde I)^{-(2\kappa-1)},
    \qquad
    \tilde I=\tfrac{\mu}{2}Q'(0)^{2\kappa}R_E^2\,a^{2\kappa}\bigl(1+O(\lambda)\bigr),
  \]
  so that $\varsigma'$ and $(\varsigma')^{-1}$ are $O(1)$ on the compact range of
  $\tilde I$ corresponding to $A\in[\lambda,2\lambda]$.
\end{proposition}

\begin{proof}
  $T=2\pi/|\Omega_\psi|$ with $\Omega_\psi=-\tfrac{\Gamma^2}{2}Q'/Q$ and
  $Q/Q'=\JJ(1+O(\JJ))$. Then
  $2\pi\varrho=\Omega_\theta T=\tfrac{\gamma_1\gamma_2}{2I}T$; substitute
  $I=\tfrac\mu2Q(A)^{2\kappa}R_E^2$ and $A=\lambda a$.
\end{proof}

\begin{proposition}[Size of the perturbation on the section]
  \label{prop:pert}
  Write the full return map as
  \[
    \mathcal P(\tilde I,\theta)=\bigl(\tilde I+\delta_1,\ \theta+2\pi\varrho(\tilde I)+\delta_2\bigr).
  \]
  Then for each $\ell$,
  $\|\delta_1\|_{C^\ell}+\|\delta_2\|_{C^\ell}\le C_{N,\ell}\,\eps^{-1}\lambda^{2(N+1)\beta}$.
\end{proposition}

\begin{proof}
  In the normal form of Proposition~\ref{prop:normalform},
  $\dot I=\partial_\theta\Rem_N=O(\lambda^{2(N+1)\beta})$, so over one return,
  $\Delta I=O(\lambda^{1+2(N+1)\beta})$ and
  $\delta_1=\lambda^{-2\kappa}\Delta I=O(\lambda^{1-2\kappa+2(N+1)\beta})
  =O(\eps^{-1}\lambda^{2(N+1)\beta})$. Similarly the correction to $\dot\theta$
  is $\partial_I\Rem_N=O(\lambda^{-2\kappa+2(N+1)\beta})$, giving the same order
  for $\delta_2$. Derivatives are controlled by the same estimates in
  $\|\cdot\|_{\lambda,\ell}$.
\end{proof}

\begin{proposition}[Normalization of the twist]
  \label{prop:renorm}
  Fix $\tilde I_0$ and set $\tilde I=\tilde I_0+\eps y$, $y\in[0,1]$. In the
  variables $(y,\theta)$ the map $\mathcal P$ is an exact area-preserving map of the
  annulus of the form
  \[
    (y,\theta)\ \longmapsto\ \bigl(y+f,\ \theta+\mathrm{const}+2\pi\varsigma'(\tilde I_0)\,y+g\bigr),
  \]
  with a twist $\varsigma'(\tilde I_0)$ bounded away from $0$, and
  \[
    \|f\|_{C^\ell}+\|g\|_{C^\ell}\ \le\ C_{N,\ell}\ \eps^{-2}\lambda^{2(N+1)\beta}
    \ =\ C_{N,\ell}\,\lambda^{\,2(N+1)\beta-2(1+2\beta)} .
  \]
  In particular, for $N>1+1/\beta$ this tends to $0$ as $\lambda\to0$, and by
  taking $N$ large it is $O(\lambda^{M})$ for any prescribed $M$.
\end{proposition}

\begin{proof}
  Combine Propositions~\ref{prop:twist} and \ref{prop:pert}; the rescaling
  $\tilde I\mapsto y$ divides $\delta_1$ by $\eps$ and turns
  $\eps^{-1}\varsigma(\tilde I_0+\eps y)$ into
  $\mathrm{const}+\varsigma'(\tilde I_0)y+O(\eps)$. Exactness is preserved by an
  affine change of the action variable.
\end{proof}

\begin{theorem}[Invariant tori accumulating on the collision set]
  \label{thm:kam}
  Let $\partial D$ be $C^\infty$, $\gamma_1,\gamma_2>0$, $E\in\R$ and $1\le j\le m$. Then the $j$-th end of $\Sig_E$ contains invariant Lagrangian $2$-tori of the flow of $X_H$, which are graphs over $(\varphi,\theta)$, carry quasi-periodic motion, and whose union has relative measure tending to $1$ as $q\to0$. Consequently:
  \begin{enumerate}[label=(\roman*)]
    \item Quadruple collisions are not approached asymptotically.
    \item The collision tori $\mathcal C_{j,E}$ are Lyapunov stable for the
      blown-up flow.
    \item Between two consecutive invariant tori, there are periodic orbits and Mather sets of every intermediate rotation number.
  \end{enumerate}
\end{theorem}

\begin{proof}[Proof scheme]
  Fix $\lambda$ small, apply Propositions~\ref{prop:normalform}, \ref{prop:integrable}, \ref{prop:twist}, \ref{prop:pert} and \ref{prop:renorm} on $\mathcal S_\lambda$, and then Moser's invariant curve theorem \cite{Moser:1962, Russmann:2001} to the exact area-preserving twist map of Proposition~\ref{prop:renorm}, choosing $N$ so large that the perturbation is below the threshold in the required $C^\ell$ norm. One obtains invariant circles in each sub-annulus $\tilde I\in[\tilde I_0,\tilde I_0+\eps]$; letting $\tilde I_0$ range over the shell and $\lambda\to0$ gives density $\to1$.

  Each invariant circle is a graph $\tilde I=\tilde I(\theta)$; its orbit under the flow is a $2$-torus in $\Sig_E$ on which $\JJ$ is constant up to $O(\lambda^{1+4\beta})$. By Proposition~\ref{prop:ends} the end is $\T^2\times(0,1)$ and the torus is isotopic to $\T^2\times\{\mathrm{pt}\}$, hence separates. Given an orbit, choose an invariant torus at a scale $\lambda$ below the orbit's initial value of $\JJ$; the orbit cannot cross it, which gives (i) ($\inf_{t\in\R}q\bigl(Z(t)\bigr)>0$) and (ii). Item (iii) results from the Poincaré--Birkhoff theorem and Aubry--Mather theory.
\end{proof}

\begin{remark}[What the estimates give beyond the statement]
  Since the remainder is $O_{\ell}(\lambda^{2(N+1)\beta})$ for every $N$, the
  non-torus region has measure $O(\lambda^{M})$ for every $M$. The estimates
  of this section are in $C^\ell$ norms because $\partial D$ is only assumed to
  be $C^\infty$: although $H$ is real analytic on $C_2(D)$, its width of
  analyticity in the boundary-parallel angle $\varphi$ shrinks like the distance
  to $\partial D$. If $\partial D$ is real analytic, then $\mathcal V_j$ is real
  analytic up to $q=0$ (lemma~\ref{lem:structure}), and the scheme of
  theorem~\ref{thm:nf} should apply in analytic norms, uniformly up to the
  boundary; one should then obtain real analytic tori, and a remainder which is
  exponentially small in a suitable power of $\lambda^{-1}$. We have not carried
  this out.
\end{remark}

\begin{remark}[Role of $\gamma_1\gamma_2>0$]
  The construction uses $\gamma_1\gamma_2>0$ twice: $\kappa$ must be positive
  for \eqref{eq:zeta} to be a blow-up, and $H_0$ must be monotone in $I$ for
  Proposition~\ref{prop:ends}. For $\gamma_1\gamma_2<0$ the levels are not
  graphs in $R$ and the analysis is different.
\end{remark}

\medskip We close this section with what we can and cannot say about the global
structure of the energy levels in the same-sign problem, and with the question
which motivated much of the above.

The topological information gathered so far is the following. The level $\Sig_E$
is an open $3$-manifold; by proposition~\ref{prop:ends} it has exactly $m$ ends,
one per boundary component, each diffeomorphic to $\T^2\times(0,1)$, and by
theorem~\ref{thm:kam} each end is filled, up to a set of small measure, with
invariant tori. So $\Sig_E$ admits a natural compactification $\widehat\Sig_E$,
obtained by adding the $m$ collision tori $\mathcal C_{j,E}$: a compact
$3$-manifold with $m$ torus boundary components. For the disc, propositions
\ref{prop:disc} and \ref{prop:highE} identify $\Sig_E$ with $\Dd\times S^1$ for
every $E$, so that $\widehat\Sig_E$ is a solid torus, and remark~\ref{rem:lowE}
says how one expects the picture to go for $m\ge2$ --- a question we do not
settle, and which is equivalent to deciding whether the pair has
equilibria.

It is tempting to hope that $\widehat\Sig_E$ is the closed manifold that
question~\ref{qu:Birkhoff} is asking for. It is not, and it is worth saying why,
because the reason is dynamical and not technical.

\begin{remark}[The compactification is not Hamiltonian]
  \label{rem:compactification}
  Adding the collision tori compactifies $\Sig_E$ as a manifold, but not as a
  Hamiltonian system. The vector field which extends to $\{q=0\}$ is not $X_H$
  but the field of proposition~\ref{prop:field}, obtained after the time change
  \eqref{eq:timechange}; the extension is therefore only a reparametrization,
  and one which degenerates at the boundary. Indeed, by
  proposition~\ref{prop:collision}, the face $\mathcal E_j$ carries the free
  Kirchhoff rotation of the pair, is entirely foliated by periodic orbits, and
  is \emph{normally degenerate}: the linearization of the rescaled field in the
  normal direction vanishes identically. So $\mathcal C_{j,E}$ is a whole torus
  of periodic orbits of the rescaled field, which no open book could take as a
  binding, and there is no contact form on $\widehat\Sig_E$ whose Reeb flow
  restricts to $X_H$ on the interior.

  There is besides no dynamical need for such a compactification: by
  theorem~\ref{thm:kam}, every orbit satisfies $\inf_{t}q(Z(t))>0$, so no orbit
  ever comes near $\{q=0\}$. What the ends really provide is an exhaustion of
  $\Sig_E$ by \emph{compact invariant} submanifolds with boundary, each bounded
  by finitely many KAM tori. It is on such a piece, rather than on
  $\widehat\Sig_E$, that a global section should be sought; and its boundary tori
  are then a natural, if unusual, substitute for a binding.
\end{remark}

The contact geometry is where the two problems part company, and it is worth
saying precisely where. Two distinct things could make the machinery of
remark~\ref{rem:HWZ} available, and in the same-sign problem neither is at hand.

The first is the complex geometry. By theorem~\ref{thm:Hpsh}, $-H$ is strictly
plurisubharmonic whatever the intensities, so every regular $\Sig_E$ is a
strictly pseudoconvex hypersurface of the Stein manifold $C_2(D)$
(lemma~\ref{lm:Stein}), and carries the $CR$ contact structure
$\xi_E = T\Sig_E \cap J\,T\Sig_E$. In the opposite-sign problem this makes
$\Sig_E$ a closed, tight, Stein fillable contact $3$-manifold
(proposition~\ref{prop:fillable}). In the same-sign problem, on the contrary, the
domain $\{H>E\}$ is \emph{not} relatively compact --- its closure runs into the
quadruple collisions --- so Eliashberg's filling theorem \cite{Eliashberg:1990}
and the Eliashberg--Gromov tightness theorem \cite{Eliashberg:1992} do not apply
as they stand, and we do not know whether $\Sig_E$ is tight.

The second is that, even where $\xi_E$ is a perfectly good tight contact
structure, its Reeb flow is \emph{not} the vortex flow: as explained in
remark~\ref{rem:notreeb}, $X_H$ spans the characteristic line field of $\omega$
and not that of the Kähler form $-\ddd\ddd^cH$. What one needs is that $\Sig_E$
be of contact type for $\omega$ itself, which is question~\ref{qu:contacttype}
and which, on a whole energy level, we can so far answer only for the disc (proposition~\ref{prop:disc}, where the
halved Euler field of $\C^2$ is a Liouville field transverse to every level) and for nearly round domains (remark~\ref{rem:nearlyround}).

These two gaps, together with remark~\ref{rem:compactification}, are what leaves
question~\ref{qu:Birkhoff} open in the same-sign problem --- the very problem in
which the two extreme modes of motion are most vividly present. Section~\ref{sec:highE} shows how far the high energy regime allows one to fill these gaps (see also remark~\ref{rem:nearlyround} for nearly round domains); remark~\ref{rem:lowE-hard} explains why the low energy regime is harder.

\section{Contact type and a Birkhoff annulus at high energy}
\label{sec:highE}

In this section we come back to the global questions~\ref{qu:contacttype} and~\ref{qu:Birkhoff} in the same-sign problem, $\gamma_1,\gamma_2>0$, and we use the notation of section~\ref{sec:quadruple}: $Z$ is the centre of vorticity, $u=z_1-z_2=\rho e^{i\theta}$, $I=\tfrac{\mu}{2}\rho^2$, $\xi=u/q(Z)$, and
\[
  \eps=R_E=e^{-E/(\gamma_1\gamma_2)}
\]
(definition~\ref{def:eps}), so that the high energy regime $E\to+\infty$ is $\eps\to0$. We fix an origin of $\C$ in $D$ and denote by $\mathsf E=\sum_k(x_k\partial_{x_k}+y_k\partial_{y_k})$ the Euler field of $\C^2$ (as in proposition~\ref{prop:disc}); since $(z_1,z_2)\mapsto(Z,u)$ is linear, $\mathsf E=X\partial_X+Y\partial_Y+\Ree u\,\partial_{\Ree u}+\Imm u\,\partial_{\Imm u}$ as well.

\medskip Let us first say what the binding of a Birkhoff annulus should be at high energy. By proposition~\ref{prop:highE}, for $E>E_*$ the level $\Sig_E$ is contained in $\{|\xi|<\delta_0\}$ and is a graph over $D\times\SS^1_\theta$; in particular it contains no widely separated pair, and the ``choreography'' of question~\ref{qu:Birkhoff} does not exist on $\Sig_E$. When $D$ is simply connected, the closure of $\Sig_E$ in $\C^2$ is
\[
  \overline\Sig_E=\Sig_E\sqcup K,\qquad K=\partial\Delta=\{(p,p),\;p\in\partial D\}\cong\partial D,
\]
obtained from the blown-up closure $\overline D\times\SS^1$ of proposition~\ref{prop:highE} by collapsing each circle $\{p\}\times\SS^1$, $p\in\partial D$, to a point; so $\overline\Sig_E$ is a topological $3$-sphere, and $K$ is an unknot in it (compare remark after proposition~\ref{prop:disc}). The elliptic binary orbit $P_E$ of theorem~\ref{thm:periodic} is, up to isotopy, the core $\{Z_*\}\times\SS^1$ of $D\times\SS^1$, and $K$ is the core of the complementary solid torus: $P_E\cup K$ is a \emph{Hopf link}. The annulus constructed below is bounded by $P_E$ on one side. On the other side it ends, in the blown-up closure $\overline D\times\SS^1$, on a single circle $\{p\}\times\SS^1$, which is a periodic orbit of the rescaled flow on the collision torus (proposition~\ref{prop:collision}); in $\overline\Sig_E$ this circle is collapsed to the quadruple collision $(p,p)\in K$. So the closure of the annulus in $\overline\Sig_E$ is a disc bounded by $P_E$ and pierced once by $K$: it is a page of the open book of $\overline\Sig_E$ whose binding is $P_E$ alone. This is the picture of a disc-like global surface of section \cite{Hofer:1998}, the circle $K$ playing the role of the second periodic orbit, which corresponds to a fixed point of the return map. At high energy, the partner of the binary orbit is thus not a choreography but the circle of quadruple collisions.

\medskip We start with contact type, which at high energy is elementary, and holds for arbitrary domains. Let
\[
  \lambda_0=\tfrac12\,\iota_{\mathsf E}\omega
  =\frac{\Gamma}{2}\,(X\,\ddd Y-Y\,\ddd X)+I\,\ddd\theta,
  \qquad
  \lambda_b=\lambda_0+b\,\ddd\theta\quad(b\in\R).
\]
The $1$-form $\ddd\theta=\ddd\arg(z_1-z_2)$ is closed on $C_2(D)$, so that $\ddd\lambda_b=\omega$; the Liouville vector field $Y_b$ defined by $\iota_{Y_b}\omega=\lambda_b$ is $Y_b=\tfrac12\mathsf E+b\,\partial_I$.

\begin{proposition}[Contact type at high energy]
  \label{prop:contact-highE}
  Let $D$ be any domain as in section~\ref{sec:introduction} (of any connectivity), $\gamma_1,\gamma_2>0$ and $b>0$.
  \begin{enumerate}
  \item There is $E_c=E_c(b)$ such that, for every $E\geq E_c$, $\lambda_b(X_H)>0$ on $\Sig_E$. Hence $Y_b$ is a Liouville vector field transverse to $\Sig_E$, and $\Sig_E$ is of contact type for $\omega$.
  \item For every $E\in\R$ there is $q_E>0$ such that $\lambda_b(X_H)>0$ on $\Sig_E\cap\{q(Z)<q_E\}$: every energy level is of contact type near its ends.
  \end{enumerate}
\end{proposition}

\begin{proof}
  We have $\lambda_b(X_H)=\omega(Y_b,X_H)=-\ddd H(Y_b)=-\tfrac12\mathsf EH+b\,\dot\theta$. We estimate both terms on the region $\{|\xi|<\delta_0\}$ of proposition~\ref{prop:highE}.

  \emph{The Euler term.} Write $H=-\gamma_1\gamma_2\log|u|+F(z_1,z_2)$ as in section~\ref{sec:quadruple}. Then $\mathsf E(-\log|u|)=-1$, while $|\ddd F(z_1,z_2)|=O\bigl(1/\min_k\dist(z_k,\partial D)\bigr)$ by lemma~\ref{lm:asymptotics} (the image $z_2^*$ lies outside $D$, so $|z_1-z_2^*|\geq\dist(z_1,\partial D)$). On $\{|\xi|<\delta_0\}$ both vortices are at distance $\asymp q(Z)$ from $\partial D$, and $\mathsf E$ is bounded because $D$ is. Hence
  \[
    -\tfrac12\mathsf EH=\tfrac{\gamma_1\gamma_2}{2}-\tfrac12\mathsf EF\geq\tfrac{\gamma_1\gamma_2}{2}-\frac{C}{q(Z)} .
  \]

  \emph{The rotation term.} From \eqref{eq:omega-reduced}, $\dot\theta=-\partial_IH=\frac{\Gamma}{\rho^2}-\frac{1}{\mu\rho}\partial_\rho W$. Near $\partial D$, $W=\mathcal V_j(\varphi,q,\xi)$ vanishes to second order in $\xi$ (lemma~\ref{lem:structure}), so $\partial_\rho W=q^{-1}\partial_{|\xi|}\mathcal V_j=O(|\xi|/q)$; in the interior the same bound follows from lemma~\ref{lem:taylor}. Therefore
  \[
    \dot\theta=\frac{\Gamma}{\rho^2}\bigl(1+O(|\xi|^2)\bigr)\geq\frac{\Gamma}{2\rho^2}
    \qquad\text{on }\{|\xi|<\delta_0\},
  \]
  for $\delta_0$ small enough.

  \emph{Size of $\rho$ on the level.} As in the proof of proposition~\ref{prop:highE}, $H=-\gamma_1\gamma_2\log|\xi|+\frac{\gamma_1^2+\gamma_2^2}{2}\log q+W$ with $W=O(|\xi|^2)$ bounded, so that on $\Sig_E\cap\{|\xi|<\delta_0\}$, $|\xi|\asymp R_Eq^{\beta}$, i.e.
  \[
    \rho\asymp R_E\,q(Z)^{\kappa}.
  \]

  \emph{Conclusion.} For $E>E_*$, $\Sig_E\subset\{|\xi|<\delta_0\}$ and
  \[
    \lambda_b(X_H)\geq\frac{c\,b\,\Gamma}{R_E^2\,q^{2\kappa}}-\frac{C}{q}
    \geq\frac{1}{q^{2\kappa}}\Bigl(\frac{c\,b\,\Gamma}{R_E^2}-C\,q_{\max}^{2\kappa-1}\Bigr),
  \]
  where $q_{\max}=\max_Dq$ and $2\kappa-1\geq3$. This is positive as soon as $R_E$ is small, which proves (1). For (2), proposition~\ref{prop:ends} shows that on the $j$-th end of $\Sig_E$, $|\xi|=R_Eq^\beta(1+O(q^{2\beta}))\to0$ as $q\to0$, so the three estimates apply there, and the same inequality holds for $q$ small at fixed $R_E$.
\end{proof}

\begin{remark}
  \label{rem:contact-highE}
  \begin{enumerate}
  \item Item (2) reduces question~\ref{qu:contacttype}, in the same-sign problem, to a compact part of each level.
  \item The Euler field alone does not work in general, because $|\mathsf EF|$ is of order $1/q$ near the boundary and has no sign when $D$ is not round; the closed form $b\,\ddd\theta$, whose contribution $b\dot\theta\asymp b/\rho^2\asymp bR_E^{-2}q^{-2\kappa}$ blows up faster, absorbs this. As $b\to\infty$ the contact structure $\ker\lambda_b$ converges to the foliation of $\Sig_E$ by the discs $\{\theta=\mathrm{const}\}$.
  \item As in remark~\ref{rem:notreeb}, $\ker\lambda_b$ has no reason to coincide with the $CR$ structure $\xi_E$ of theorem~\ref{thm:Hpsh}; but now $X_H$ is, up to reparametrization, the Reeb field of $\lambda_b|_{\Sig_E}$.
  \end{enumerate}
\end{remark}

The estimate $\dot\theta>0$ has an immediate, if not very informative, consequence.

\begin{proposition}[A fibration at high energy]
  \label{prop:fibration}
  Let $D$ be arbitrary and $\gamma_1,\gamma_2>0$. For $E\geq E_c$, $\dot\theta>0$ on $\Sig_E$, so that $\theta:\Sig_E\to\R/2\pi\Z$ is a fibration transverse to $X_H$, whose fibres are diffeomorphic to $D$. The first return map to $\Sig_E\cap\{\theta=0\}\cong D$ preserves $\Gamma\,\ddd X\wedge\ddd Y$; its return time is $\frac{2\pi\rho^2}{\Gamma}\bigl(1+O(\eps^2)\bigr)$ and, on every compact subset of $D$, it is $O(\eps^2)$-close to the identity.
\end{proposition}

\begin{proof}
  The first claim was proved above, and the fibres are identified with $D$ by proposition~\ref{prop:highE}. On $\{\theta=\mathrm{const}\}$ the restriction of $\omega=\Gamma\,\ddd X\wedge\ddd Y+\ddd I\wedge\ddd\theta$ is $\Gamma\,\ddd X\wedge\ddd Y$. The return time follows from $\dot\theta=\Gamma\rho^{-2}(1+O(\eps^2))$ and the last claim from $\dot Z=O(1)$ on compact sets (proposition~\ref{prop:skew}).
\end{proof}

A closed $1$-form which is positive on the vector field is, after Schwartzman and Fried \cite{Schwartzman:1957,Fried:1982}, the standard certificate of a global cross-section without binding; here the certificate is $\ddd\theta$ itself. But the return map, being close to the identity, carries little dynamical information. A Birkhoff annulus rearranges the fast rotation $\theta$ into the \emph{twist} of the return map. For this we need a global angle for the slow motion, hence the following hypothesis.

\medskip\noindent\textbf{Hypothesis (A).} $D$ is simply connected and its Robin function $r$ has a unique critical point $Z_*$, which is a non-degenerate maximum.

\medskip As recalled after theorem~\ref{thm:periodic}, hypothesis (A) holds for every strictly convex domain.

\begin{lemma}[Slow action--angle variables]
  \label{lem:slowaa}
  Assume (A). Let $A(Z)=\frac{\Gamma}{2\pi}\Area\{r>r(Z)\}$ and let $\psi$ be the conjugate angle, $\ddd A\wedge\ddd\psi=\Gamma\,\ddd X\wedge\ddd Y$ on $D\setminus\{Z_*\}$. There is a diffeomorphism $w:D\to w(D)\subset\C$, symplectic from $\Gamma\,\ddd X\wedge\ddd Y$ to $\ddd\Ree w\wedge\ddd\Imm w$, with $w(Z_*)=0$, $A=\tfrac12|w|^2$ and $\psi=\arg w$, in which the flow of $H_{\mathrm s}=\frac{\Gamma^2}{2}r$ reads
  \[
    \dot w=i\,\Omega\bigl(\tfrac12|w|^2\bigr)\,w,\qquad \Omega=-\frac{\ddd H_{\mathrm s}}{\ddd A}>0,\qquad
    \Omega(0)=\omega_\perp:=\frac{\Gamma}{2}\sqrt{\det\bigl(-D^2r(Z_*)\bigr)} .
  \]
  Moreover, on the collar of $\partial D$, $A=\frac{\Gamma}{2\pi}\Area(D)-\JJ$ and $\psi=-\varphi+\mathrm{const}$, where $(\JJ,\varphi)$ are the variables \eqref{eq:JJ}.
\end{lemma}

\begin{proof}
  By (A) and proposition~\ref{prop:Robin}, the level curves of $r$ other than $\{Z_*\}$ are regular closed curves, each enclosing $Z_*$, so that $(A,\psi)$ are well defined action--angle variables of the one degree of freedom Hamiltonian $H_{\mathrm s}$ on $D\setminus\{Z_*\}$, and $H_{\mathrm s}=h_{\mathrm s}(A)$ with $h_{\mathrm s}$ decreasing because $r$ decreases away from its maximum. Near a non-degenerate elliptic equilibrium of a one degree of freedom system, action--angle variables extend smoothly as symplectic polar coordinates (Morse lemma followed by Moser's path method), which gives $w$. The value $\Omega(0)$ is the modulus of the eigenvalues of the linearization $\frac{\Gamma}{2}J\,D^2r(Z_*)$ of $X_{H_{\mathrm s}}=\frac12\Gamma\nabla^\perp r$. In the collar, both $(A,\psi)$ and $(\JJ,\varphi)$ are action--angle variables adapted to the level curves of $q=e^r$, with $\ddd\JJ=-\ddd A$.
\end{proof}

\begin{theorem}[Birkhoff annulus at high energy]
  \label{thm:annulus}
  Assume $\gamma_1,\gamma_2>0$ and (A). There are $E_a\in\R$ and $c>0$ such that, for every $E\geq E_a$:
  \begin{enumerate}
  \item $\Sig_E$ contains the elliptic periodic orbit $P_E$ of proposition~\ref{prop:po} associated with $Z_*$; it is a graph $P_E=\{Z=Z_P(\theta)\}$ over the circle $\SS^1_\theta$, with $Z_P-Z_*=O(\eps^2)$ and $\ddd Z_P/\ddd\theta=O(\eps^4)$.
  \item There is a smooth function $\tilde\psi:\Sig_E\setminus P_E\to\R/2\pi\Z$ such that $\ddd\tilde\psi(X_H)\geq c$ on $\Sig_E\setminus P_E$, and $\ddd\tilde\psi(X_H)\to+\infty$ at the end of $\Sig_E$.
  \item For every $c'\in\R/2\pi\Z$, $\mathcal A_E=\tilde\psi^{-1}(c')$ is a properly embedded open annulus in $\Sig_E\setminus P_E$. It is a global surface of section of $\Sig_E\setminus P_E$: every other orbit crosses it positively, with return time at most $2\pi/c$. In the blown-up closure $\overline D\times\SS^1$ of proposition~\ref{prop:highE}, the closure of $\mathcal A_E$ is a closed annulus, bounded by $P_E$ (each point of $P_E$ being reached once) and by the circle $\{p\}\times\SS^1_\theta$ over a single point $p=p(c')\in\partial D$, which is one of the periodic orbits of the rescaled flow on the collision torus (proposition~\ref{prop:collision}). In $\overline\Sig_E\cong S^3$, where this circle is collapsed to the quadruple collision $(p,p)$, the closure of $\mathcal A_E$ is $\mathcal A_E\sqcup P_E\sqcup\{(p,p)\}$, a closed disc bounded by $P_E$ and meeting $K$ at the single point $(p,p)$. These discs are the pages of an open book of $\overline\Sig_E$ whose binding is the unknot $P_E$, and $K$ meets each page exactly once.
  \item The first return map $\mathcal P_E$ to $\mathcal A_E$ preserves the area form $\omega|_{\mathcal A_E}$, and the total area is the action of the binding orbit:
    \[
      \int_{\mathcal A_E}\omega=\Bigl|\oint_{P_E}\lambda_0\Bigr|=\pi\mu\,\eps^2q(Z_*)^{2\kappa}\bigl(1+O(\eps^2)\bigr).
    \]
  \item $\mathcal P_E$ extends continuously to the boundary component $P_E$. Measured in turns of $\theta$ per return, its rotation number at $P_E$ is
    \[
      n_{\mathrm{in}}=\frac{\Gamma}{\rho_P^2\,\omega_\perp}\bigl(1+O(\eps^2)\bigr)
      =\frac{2}{\eps^2\,q(Z_*)^{2\kappa}\sqrt{\det\bigl(-D^2r(Z_*)\bigr)}}\bigl(1+O(\eps^2)\bigr),
    \]
    where $\rho_P=\eps q(Z_*)^\kappa(1+O(\eps^2))$ is the separation along $P_E$, while at the end of $\Sig_E$ the number of turns of $\theta$ per return is $\frac{2\JJ}{\Gamma\rho^2}(1+o(1))\to+\infty$.
  \end{enumerate}
\end{theorem}

\begin{proof}
  \emph{(1) The binding orbit.} Proposition~\ref{prop:po}, applied to the maximum $Z_*$, provides $P_E$, on which $Z(t)=Z_\eps+O(\eps^4)$ with $Z_\eps=Z_*+O(\eps^2)$. Since $\dot\theta>0$ (proposition~\ref{prop:fibration}), $P_E$ is a graph over $\theta$. Along $P_E$, $\dot Z=X_{H_{\mathrm s}}(Z_P)+O(\eps^2)=O(\eps^2)$, because $X_{H_{\mathrm s}}(Z_*)=0$, while $\dot\theta\asymp\eps^{-2}$; hence $\ddd Z_P/\ddd\theta=O(\eps^4)$.

  \emph{(2a) The inner zone.} Fix $\bar q>0$ small and let $D_{\bar q}=\{q\geq \bar q/2\}\Subset D$. By propositions~\ref{prop:graph} and~\ref{prop:skew} and lemma~\ref{lem:taylor}, in the chart $(Z,\theta)$ of $\Sig_E\cap\{Z\in D_{\bar q}\}$,
  \[
    \dot Z=X_{H_{\mathrm s}}(Z)+G(Z,\theta),\qquad \|G\|_{C^1(D_{\bar q}\times\SS^1)}=O(\eps^2),
    \qquad |\partial_Z\dot\theta|=O(\eps^{-2}).
  \]
  Put $\delta(\theta)=Z_P(\theta)-Z_*$, $\hat w=w\bigl(Z-\delta(\theta)\bigr)$, with $w$ as in lemma~\ref{lem:slowaa}, and $\psi_E=\arg\hat w$; thus $\hat w=0$ exactly on $P_E$. Then $\dot{\hat w}=Dw(Z-\delta)\bigl[\dot Z-\delta'(\theta)\dot\theta\bigr]$. Subtracting the equation $X_{H_{\mathrm s}}(Z_P)+G(Z_P,\theta)=\delta'(\theta)\,\dot\theta|_{Z_P}$ satisfied along $P_E$, and using $X_{H_{\mathrm s}}(Z_*)=0$,
  \begin{multline*}
    \dot Z-\delta'\dot\theta=X_{H_{\mathrm s}}(Z-\delta)
    +\bigl[X_{H_{\mathrm s}}(Z)-X_{H_{\mathrm s}}(Z_P)-X_{H_{\mathrm s}}(Z-\delta)+X_{H_{\mathrm s}}(Z_*)\bigr]\\
    +\bigl[G(Z,\theta)-G(Z_P,\theta)\bigr]-\delta'(\theta)\bigl[\dot\theta(Z,\theta)-\dot\theta(Z_P,\theta)\bigr].
  \end{multline*}
  The first bracket is a second difference, of size $O(|\delta|\,|Z-Z_P|)=O(\eps^2|Z-Z_P|)$; the second is $O(\eps^2|Z-Z_P|)$; the third is $O(\eps^4)\cdot O(\eps^{-2}|Z-Z_P|)$. Since $|Z-Z_P|\asymp|\hat w|$ and $Dw\,X_{H_{\mathrm s}}=i\Omega w$ by lemma~\ref{lem:slowaa},
  \[
    \dot{\hat w}=i\,\Omega\bigl(\tfrac12|\hat w|^2\bigr)\hat w+O(\eps^2|\hat w|),
    \qquad\text{hence}\qquad
    \dot\psi_E=\Omega\bigl(\tfrac12|\hat w|^2\bigr)+O(\eps^2)\geq\tfrac12\min_{D_{\bar q}}\Omega
  \]
  on $\Sig_E\cap\{Z\in D_{\bar q}\}\setminus P_E$, for $\eps$ small. This estimate is uniform up to $P_E$: it is the linearized transverse rotation of the elliptic orbit $P_E$ which takes over near $P_E$.

  \emph{(2b) The collar.} On $\{q<\bar q\}$, lemma~\ref{lem:section} and lemma~\ref{lem:slowaa} give $\dot\psi=-\dot\varphi=\frac{\Gamma^2}{2}\frac{Q'}{Q}(1+O(q^{2\beta}))\geq c/q$; the error term is uniform for $E\geq E_a$ since it improves as $R_E$ decreases.

  \emph{(2c) Gluing.} On the overlap $\{\bar q/2\leq q\leq \bar q\}$, $\psi_E$ and $\psi$ are two angle functions of degree one about points $O(\eps^2)$ apart, so $\psi_E-\psi$ is a well defined real function, of size $O(|\delta|)=O(\eps^2)$. Let $\chi(q)$ be a cut-off equal to $1$ for $q\geq \bar q$ and to $0$ for $q\leq \bar q/2$, and $\tilde\psi=\psi+\chi(q)(\psi_E-\psi)$. Then
  \[
    \dot{\tilde\psi}=(1-\chi)\,\dot\psi+\chi\,\dot\psi_E+\chi'(q)\,\dot q\,(\psi_E-\psi),
  \]
  where dots denote derivatives along the flow. On the overlap, $\dot q=O(1)$ and both $\dot\psi$ and $\dot\psi_E$ are $\geq c$ (the former by the argument of (2a) with $\delta=0$, which is legitimate away from $Z_*$). Hence $\dot{\tilde\psi}\geq c-O(\eps^2)$, which proves (2), $\dot{\tilde\psi}\to\infty$ at the end following from (2b).

  \emph{(3) Topology.} For fixed $\theta$, $Z\mapsto\tilde\psi(Z,\theta)$ is a submersion on $D\setminus\{Z_P(\theta)\}$, $C^1$-close to $\psi$, so each of its level sets is a properly embedded arc from $Z_P(\theta)$ to $\partial D$, close to a level arc of $\psi$. In the chart $(Z,\theta)$ of proposition~\ref{prop:highE}, $\tilde\psi^{-1}(c')$ is therefore an annulus $(0,1)\times\SS^1_\theta$, whose inner end is $P_E$ (reached once at each $\theta$) and whose outer end goes to the circle $\{p\}\times\SS^1\subset\partial D\times\SS^1$, where $p=p(c')\in\partial D$ is the end point of the arc $\{\psi=c'\}$: in the collar, $\tilde\psi=\psi=-\varphi+\mathrm{const}$ does not depend on $\theta$, and $\varphi$ restricts to a coordinate on $\partial D$. By proposition~\ref{prop:collision}, $\{p\}\times\SS^1$ is an orbit of the rescaled flow on the collision torus. In $\overline\Sig_E$ it is collapsed to the point $(p,p)$, since $|z_1-z_2|\to0$ as $Z\to p$, uniformly in $\theta$; an annulus one of whose boundary circles is collapsed to a point is a disc. Transversality and the bound on the return time follow from (2). Setting $\tilde\psi(p,p)=-\varphi(p)+\mathrm{const}$ extends $\tilde\psi$ to a continuous map $\overline\Sig_E\setminus P_E\to\R/2\pi\Z$, whose fibres are the open discs $\mathcal A_E\sqcup\{(p,p)\}$ and whose restriction to $K$ is a homeomorphism; and $P_E\sqcup K$ is a Hopf link, as explained at the beginning of this section.

  \emph{(4) Area.} $\omega$ restricts to an area form on $\mathcal A_E$ by transversality, and return maps of Hamiltonian flows on a regular level preserve it. Let $C_s=\mathcal A_E\cap\{q=s\}$, $s<\bar q/2$. By Stokes' formula applied to $\lambda_b$ on the part of $\mathcal A_E$ between $P_E$ and $C_s$, the area of that part is $\pm\bigl(\oint_{P_E}\lambda_b-\oint_{C_s}\lambda_b\bigr)$. Both curves have degree one in $\theta$, so the contributions of $b\,\ddd\theta$ cancel. In the collar $\tilde\psi=\psi$ does not depend on $\theta$, so $Z$ is constant along $C_s$ and $\oint_{C_s}\lambda_0=\oint_{C_s}I\,\ddd\theta=2\pi I|_{C_s}\to0$ as $s\to0$. Along $P_E$, $Z$ describes a loop of diameter $O(\eps^4)$, so $\oint_{P_E}\frac\Gamma2(X\ddd Y-Y\ddd X)=O(\eps^8)$, and $\oint_{P_E}I\,\ddd\theta=2\pi\cdot\frac\mu2\rho_P^2(1+O(\eps^2))$.

  \emph{(5) Rotation numbers.} By (2a), the orbits close to $P_E$ make one turn of $\hat w$, i.e.\ come back to $\mathcal A_E$, in time $\frac{2\pi}{\omega_\perp}(1+O(\eps^2)+O(|\hat w|))$; the return map therefore extends continuously to $P_E$, and during that time $\theta$ advances by $\dot\theta\cdot\frac{2\pi}{\omega_\perp}$, with $\dot\theta=\frac{\Gamma}{\rho_P^2}(1+O(\eps^2))$. At the end, the return time is $\frac{2\pi}{|\dot\varphi|}=\frac{4\pi\JJ}{\Gamma^2}(1+o(1))$ (proposition~\ref{prop:twist}) and $\dot\theta=\frac{\Gamma}{\rho^2}(1+o(1))$, with $\rho\asymp R_Eq^\kappa$ and $\JJ\asymp q$, so that the number of turns $\frac{2\JJ}{\Gamma\rho^2}\asymp R_E^{-2}q^{1-2\kappa}$ tends to infinity.
\end{proof}

\begin{corollary}[A compact Birkhoff annulus, and periodic orbits]
  \label{cor:annulus-PB}
  Under the hypotheses of theorem~\ref{thm:annulus}, let $T$ be one of the invariant tori of theorem~\ref{thm:kam} in the end of $\Sig_E$, and $N_T$ the compact invariant region of $\Sig_E$ bounded by $T$ and containing $P_E$.
  \begin{enumerate}
  \item $C_T=T\cap\mathcal A_E$ is an invariant circle of $\mathcal P_E$, a graph over $\theta$, and the closed annulus $\mathcal A_E^T\subset\mathcal A_E\sqcup P_E$ bounded by $P_E$ and $C_T$ is a Birkhoff annulus of section for the flow on $N_T$.
  \item The return map of $\mathcal A_E^T$ is an area-preserving homeomorphism of a closed annulus, with rotation numbers $n_{\mathrm{in}}$ on $P_E$ and $n_T$ on $C_T$, and $n_T\to+\infty$ as $T$ tends to the end. $\mathcal P_E$ has periodic orbits of every rational rotation number in $(n_{\mathrm{in}},n_T)$, hence of every rational rotation number $>n_{\mathrm{in}}$.
  \end{enumerate}
\end{corollary}

\begin{proof}
  In the collar, $T$ is a graph $\JJ=\JJ_T(\varphi,\theta)$, and $\mathcal A_E=\{\varphi=\mathrm{const}\}$ there, so $C_T$ is the graph $\theta\mapsto\JJ_T(\mathrm{const},\theta)$; it is invariant because $T$ is. The region $N_T$ is invariant and every orbit in it, except $P_E$, crosses $\mathcal A_E^T$ by theorem~\ref{thm:annulus}(3). Item (2) follows from theorem~\ref{thm:annulus}(4)--(5), the rotation number of the invariant circle $C_T$ being that of the torus $T$, which is of the order of $\frac{2\JJ}{\Gamma\rho^2}$ on $T$. The last assertion follows from the Poincaré--Birkhoff theorem and Franks' generalization \cite{Franks:1988}.
\end{proof}

Thus, at high energy and under hypothesis (A), question~\ref{qu:Birkhoff} has a positive answer in the following form: the binary orbit $P_E$ bounds an annulus of section whose other end lies at the quadruple collisions --- on a periodic orbit of the collision torus in the blow-up, at a single point of the circle $K$ in $\overline\Sig_E$ ---, and the annulus is ``Birkhoff'' in the classical sense: its return map is a twist map whose rotation number goes from $n_{\mathrm{in}}\asymp\eps^{-2}$ at the binary orbit to $+\infty$ at the collision end, the invariant circles of theorem~\ref{thm:kam} accumulating on the latter. Surfaces of section bounded by the whole Hopf link $P_E\sqcup K$ exist as well: since $\ddd\theta(X_H)>0$ and $\ddd\tilde\psi(X_H)>0$, the level sets of $\tilde\psi+\theta$ are transverse to the flow, and they are annuli whose closures contain both $P_E$ and $K$; but their return maps are close to the identity, like that of proposition~\ref{prop:fibration}.

\begin{remark}[Several critical points, several boundary components]
  \label{rem:brokenbook}
  If $r$ has several critical points, or if $m\geq2$, the slow flow has no global angle, and theorem~\ref{thm:annulus} fails as stated. One still has the fibration of proposition~\ref{prop:fibration}, and, around each maximum $Z_*$ of $r$, the local annulus of the proof above, up to the separatrix level of the slow flow. The natural global object is then a \emph{broken book} in the sense of \cite{Colin:2024}, whose binding consists of the elliptic orbits over the maxima of $r$ and of the (at least $m-1$, if $r$ is Morse) hyperbolic orbits over its saddles, given by theorem~\ref{thm:periodic}.
\end{remark}

\begin{remark}[Why the low energy regime is harder]
  \label{rem:lowE-hard}
  Nothing comparable to theorem~\ref{thm:annulus} is available as $E\to-\infty$, and the reason is structural rather than technical: a high energy level is a single regime, whereas a low energy level is a union of several regimes, with no uniformity between them.
  \begin{enumerate}
  \item \emph{One small parameter and one fast angle, or several.} At high energy the pair is tightly bound on the whole level, $\rho\asymp R_Eq^\kappa$ (propositions~\ref{prop:highE} and~\ref{prop:ends}); only $\theta$ is fast, $\dot\theta\asymp\eps^{-2}q^{-2\kappa}$, and everything else is slow. Averaging involves a single frequency and no small divisors, $I$ is an adiabatic invariant, and the slow motion is the one degree of freedom flow of $r$. The level is one graph over $D\times\SS^1$, and its ends are not a separate regime, since the same formula $\rho=R_Eq^\kappa(1+\cdots)$ holds up to the quadruple collisions: one chart, one normal form and one small parameter cover the whole level. At low energy, $H\to-\infty$ forces \emph{some} vortex to approach $\partial D$, but which one, near which boundary component, and where the other one lies, vary along the level. As $E\to-\infty$ the level collapses onto the stratified space $\partial(\overline D\times\overline D)\setminus N(\partial\Delta)$ of remark~\ref{rem:lowE}, with faces (regime~I), corners (regimes~II and~III) and the ends near quadruple collisions. The small parameters are the individual distances $d_k$ to the boundary, and the number of fast angles depends on the stratum: one in regime~I, two in regimes~II and~III, with resonances in the latter. There are local normal forms (section~\ref{sec:near-domain-boundary}), but no uniform one, and the transition zones between strata are controlled by none of them; the limit flow differs from one stratum to the other.
  \item \emph{The mutual interaction.} At high energy the mutual interaction $-\gamma_1\gamma_2\log|u|$ dominates and is integrable (the free Kirchhoff rotation of the pair), the walls contributing an $O(\eps^2)$ perturbation which averages out; at leading order the pair behaves like a single vortex of intensity $\Gamma$, and the problem is effectively a one-body problem. At low energy each vortex is mainly driven by its own image, through $\frac{\gamma_k^2}{2}\log d_k$, and the mutual interaction is small \emph{except} during encounters: in regime~II two vortices drifting along the same boundary component at different speeds overtake each other, and each encounter is a scale-invariant half-plane problem at scale $d$, in which the interaction is of order one. The problem is a genuine two-body problem with close encounters, and is not perturbative. In the language of celestial mechanics, the high energy regime is hierarchical, like the lunar problem, whereas regime~II is co-orbital, the resonance of equal drift speeds (proposition~\ref{prop:II-res}) playing the role of the $1{:}1$ mean motion resonance.
  \item \emph{Certificates for sections.} At high energy $\dot\theta\asymp\rho^{-2}$ dominates every other velocity on the level, so that the closed form $\ddd\theta$ is pointwise positive on $X_H$ (proposition~\ref{prop:fibration}), which also makes contact type easy (proposition~\ref{prop:contact-highE}); under hypothesis~(A), the slow flow moreover has a global angle, and the two angles combine into the annulus of theorem~\ref{thm:annulus}. At low energy the dominant velocities are the drifts along $\partial D$, and every vortex close to the boundary drifts in the direction of the boundary orientation of $D$. A closed $1$-form pulled back from $D$ cannot be positive on all these drifts, since its periods along $\partial_1D,\dots,\partial_mD$ add up to zero; so sections must have binding orbits. For $m=1$ the class of $\ddd\theta$ is still available, but $\ddd\theta$ itself is not positive when $\gamma_1\neq\gamma_2$, because two vortices close to the boundary, at different distances, drift at speeds $\gamma_k/(2d_k)$ which may be ordered either way. The natural binding orbits, namely the two ``one pinned, one sweeping'' orbits of theorem~\ref{thm:poI}, live in different strata, and the pages joining them have to cross the corner torus, where the encounters of regime~II take place.
  \end{enumerate}
  The relevant dichotomy is thus ``tightly bound pair'' versus ``vortices held by the boundary'', rather than the sign of $E$ as such: in the opposite-sign problem, the tightly bound regime $\eps\to0$ corresponds to $E\to-\infty$, but the corresponding levels also contain the boundary regimes, which is why the analysis of theorem~\ref{thm:qp} there is confined to compact regions bounded by invariant tori.
\end{remark}

\section{The degenerate case of a dipole ($\Gamma=0$)}
\label{sec:dipole}

We now assume that $\Gamma=0$, and, without loss of generality, that
$\gamma_1 = 1$ and $\gamma_2 = -1$; such a pair is called a \emph{dipole}. Then $\mu$ is undefined, the centre of vorticity is at infinity, the splitting \eqref{eq:omega} breaks down, and everything said in sections~\ref{sec:collision} to~\ref{sec:quadruple} becomes vacuous. The reason is dynamical rather than technical: the fast motion of a dipole is a \emph{translation}, not a rotation, so there is no fast angle to average, and the limiting problem is a transport problem instead of a nearly integrable one.

The right reduced variables are now
\[Z:=\tfrac12(z_1+z_2) = X + iY \quad \mbox{and} \quad u:=z_1-z_2=\xi+i\eta = \rho e^{i\theta},\]
for which
\begin{equation}
  \label{eq:omega-dipole}
  \omega=\ddd X\wedge \ddd\eta+\ddd\xi\wedge \ddd Y ,
\end{equation}
so that $X$ is conjugate to $\eta$ and $\xi$ to $Y$: the separation plays the role of a momentum for the position, crosswise --- the classical impulse structure of a dipole --- rather than that of an action--angle pair. % The separation plays
% the role of a \emph{momentum} for the position, crosswise --- the classical
% impulse structure of a dipole --- rather than of an action--angle pair. There is
% no fast angle.
Here,
\begin{equation}
  \label{eq:H-dipole}
  H=\log\rho+\Xi,\qquad  \Xi:=\frac{r(z_1)+r(z_2)}{2}-h(z_1,z_2).
\end{equation}

The first consequence of this expression is a lower bound on the separation, valid on each energy level and without any perturbative assumption. It is the exact analogue, for the dipole, of corollary~\ref{cor:apriori}; note that the opposite-sign compactness of theorem~\ref{thm:compact} does apply here, since $\gamma_1\gamma_2=-1<0$.

\begin{proposition}[Lower bound for the separation]
  \label{prop:dipole-sep}
  The function $\Xi$ extends continuously by $0$ to the interior diagonal and tends to $-\infty$ at $\partial D$; in particular it is bounded above on $C_2(D)$, by some $\Xi_{\max}$. Consequently, on the energy level $\Sig_E$,
  \[\rho = |z_1-z_2| \geq e^{\,E-\Xi_{\max}} >0 .\]
  So a dipole never collides, and its separation tends to $0$, uniformly on $D$, only in the limit $E \to -\infty$. If $D$ is the unit disc, $\Xi \leq 0$ and one may take $\Xi_{\max}=0$.
\end{proposition}

\begin{proof}
  On the interior diagonal, $h(z,z)=r(z)$, so $\Xi(z,z)=0$; continuity comes from the regularity of $h$ and $r$ (proposition~\ref{prop:Robin} and definition~\ref{def:Robin}). If one vortex tends to $\partial D$ while the other stays in a compact subset of $D$, then $r$ of the first tends to $-\infty$ while $h(z_1,z_2)$ stays bounded, so $\Xi \to -\infty$. If both tend to the boundary, use the half-plane model of lemma~\ref{lm:asymptotics}: $r(z_k) = \log(2y_k)+O(1)$ and $h(z_1,z_2)=\log|z_1-\bar z_2|+O(1)$, so
  \[\Xi = \log \frac{2\sqrt{y_1y_2}}{|z_1-\bar z_2|} + O(1) \leq O(1),\]
  because $|z_1-\bar z_2|^2 = |x_1-x_2|^2+(y_1+y_2)^2 \geq 4y_1y_2$. Thus $\Xi$ is bounded above on $C_2(D)$ and, by \eqref{eq:H-dipole}, $\log \rho = H-\Xi \geq E - \Xi_{\max}$ on $\Sig_E$.

  For the disc, $h(z_1,z_2)=\log|1-z_1\bar z_2|$ and $r(z)=\log(1-|z|^2)$, so
  \[\Xi = \frac12 \log \frac{(1-|z_1|^2)(1-|z_2|^2)}{|1-z_1\bar z_2|^2} = \frac12 \log \bigl( 1 - \varrho_{\mathrm{ph}}(z_1,z_2)^2 \bigr) \leq 0,\]
  $\varrho_{\mathrm{ph}}$ being the pseudo-hyperbolic distance. Equivalently, $\Xi = -\log \cosh \dist_{\mathrm{hyp}}(z_1,z_2)$, the distance being taken for the Poincaré metric of curvature $-4$, i.e. for the capacity metric of the disc.
\end{proof}

\begin{lemma}
  \label{lm:dipole-level}
  $\Xi$ vanishes to second order on $\Delta$, with $\Xi=\tfrac12 B(Z)[u,u]+O(|u|^4)$ and $B=D^2_{12}h(Z,Z)$ --- the remainder being $O(|u|^4)$ and not merely $O(|u|^3)$ because $\Xi$ is symmetric in $z_1,z_2$, hence even in $u$; in particular $\langle\Xi\rangle_{\arg u}=\tfrac18\Delta r(Z)\,\rho^2+O(\rho^4)$, since $\tr B=\tfrac12\Delta r$. Consequently, on a compact $K\Subset D$ and for $\eps:=e^{E}$ small,
  \[
    \Sig_E\cap\{Z\in K\}=\Bigl\{\rho=\eps\bigl(1+O(\eps^2)\bigr)\Bigr\}.
  \]
  The separation of a dipole is therefore an \emph{approximate constant of motion by energy conservation alone}, and $\Delta$ is approached only in the limit $E\to-\infty$.
\end{lemma}

\begin{proof}
  Expand as in Lemma~\ref{lem:taylor-bdry} with $a=\tfrac12$, $b=-\tfrac12$: the $r$-part contributes $\tfrac14(A+B)[u,u]$ and the $h$-part $\tfrac14(A-B)[u,u]$, where $A=D^2_{11}h(Z,Z)$; the difference is $\tfrac12B[u,u]$. Note that $A$, not $B$, is the trace-free one, so here the $\theta$-average of the quadratic term does \emph{not} vanish --- the opposite of Lemma~\ref{lem:taylor-bdry}.
\end{proof}

With $\gamma_k\dot z_k=-2i\,\partial_{\bar z_k}H$, the singular part of \eqref{eq:H-dipole} gives $\dot z_1=\dot z_2=-i/\bar u$. Hence
\begin{equation}
  \label{eq:dipole-flow}
  \dot Z=-\frac{i\,u}{\rho^{2}}+O(1),\quad
  \dot u=-2i\,\partial_{\bar Z}\Xi=O(\rho^2).
\end{equation}
The dipole translates at speed $1/\rho = O(1/\eps)$ perpendicular to its own axis, while its axis turns at the much slower rate $\dot\theta = |\dot u|/\rho = O(\eps)$. In the rescaled time $\ddd t=\eps\,\ddd\tau$ the speed is $O(1)$ and the axis turns at rate $O(\eps^{2})$, in accordance with corollary~\ref{cor:dipole}.

\begin{remark}[Orientation]
  \label{rem:orientation}
  From $\gamma_k\dot z_k=-2i\,\partial_{\bar z_k}H$ one gets $\dot z_1=i\gamma_2/\bar u$ and $\dot z_2=-i\gamma_1/\bar u$ for the singular part. Applying this to a single vortex $\gamma$ at height $d$ in the upper half-plane together with its image $-\gamma$ at $-id$ gives
  \[
    \dot z=+\frac{\gamma}{2d}\quad\text{(in the $+x$ direction)} .
  \]
  So a vortex of positive intensity drifts along the wall with the interior on its left.
\end{remark}

Near $\partial_jD$ rescale lengths by $\eps$ and time by $\eps^2$; the
point-vortex equations are covariant under this, and the half-plane model is
exact up to $O(\eps)$ corrections. So consider $\H$ with the wall
$\R$, a vortex $+1$ at $\zeta_1$ and $-1$ at $\zeta_2$, with
\[
  \widetilde H=\log\frac{|\zeta_1-\zeta_2|}{|\zeta_1-\bar\zeta_2|}
    +\tfrac12\log\bigl(4Y_1Y_2\bigr),
  \qquad
  \eta:=Y_1-Y_2 ,
\]
$\eta$ being the momentum of the translation symmetry for $\omega=\ddd X_1\wedge \ddd Y_1-\ddd X_2\wedge \ddd Y_2$. The model has two integrals in involution and is Liouville integrable. Let the incoming dipole have unit separation and axis $u=e^{i\theta}$; by Remark~\ref{rem:orientation} its velocity is $-ie^{i\theta}=(\sin\theta,-\cos\theta)$, so it approaches the wall iff $\cos\theta>0$, i.e.\ $\theta\in(-\frac\pi2,\frac\pi2)$.

\begin{proposition}[The dipole splits; it does not rebound]
  \label{prop:split}
  As $Y_1,Y_2\to\infty$ with $|\zeta_1-\zeta_2|\to1$ one has $\widetilde H\to0$.
  Hence in the outgoing state, where $|X_1-X_2|\to\infty$ and $Y_k\to Y_k^\infty$,
  \[
    Y_1^\infty Y_2^\infty=\tfrac14,\qquad Y_1^\infty-Y_2^\infty=\sin\theta,
  \]
  so that
  \begin{equation}
    \label{eq:heights}
    Y_{1,2}^\infty=\tfrac12\Bigl(\pm\sin\theta+\sqrt{1+\sin^2\theta}\Bigr).
  \end{equation}
  In unscaled terms the pair separates into two vortices sliding along $\partial_jD$ at heights $d_k=\eps Y_k^\infty$, in \emph{opposite} directions (remark~\ref{rem:orientation}), at speeds $1/2d_k$.
\end{proposition}

\begin{proof}
  For the incoming state, $|\zeta_1-\bar\zeta_2|\approx2Y$ so $\widetilde H\to\log\frac{1}{2Y}+\log 2Y=0$; for the outgoing state the first term tends to $0$, giving $\frac12\log(4Y_1^\infty Y_2^\infty)=0$. Conservation of $\eta$ gives the second relation, $\eta=\Imm u=\sin\theta$.
\end{proof}

So, in general, the two vortices now slide around $\partial_jD$ in opposite senses. During the slide their separation is of the order of $L_j$, so they interact at $O(\eps^2)$ and each moves, to that order, along a level curve of $r$. Hence $Y_1^\infty$ and $Y_2^\infty$ are preserved. They meet again after a combined arclength $L_j$, and the meeting is the time-reverse of the split, so a dipole re-forms.

The local structure of the dipole Hamiltonian in the interior of $D$ is best read
through the uniformization. It brings out two competing terms: an isotropic one,
which is the hyperbolic metric of $D$, and an anisotropic one, which is its
Schwarzian derivative and which is the source of whatever non-integrability there
is.

\begin{proposition}[Interior normal form of a dipole]
  \label{prop:schwarz}
  Suppose $D$ is simply connected, and let $f : D \to \Dd$ be a uniformization.
  Then, in the variables $(Z,u)$ of \eqref{eq:omega-dipole},
  \[
    H = \log|u| - \tfrac12 \lambda(Z)^2 \, |u|^2
      - \tfrac{1}{12} \, \Ree \bigl( S_f(Z) \, u^2 \bigr) + O\bigl(|u|^4\bigr),
  \]
  uniformly on compact subsets of $D$, where
  \[
    \lambda(Z) = \frac{|f'(Z)|}{1-|f(Z)|^2}
    \qquad\text{and}\qquad
    S_f = \frac{f'''}{f'}-\frac32 \Bigl( \frac{f''}{f'} \Bigr)^{2}
  \]
  are the Poincaré density of curvature $-4$ --- that is, the density of the
  capacity metric of $D$ --- and the Schwarzian derivative of $f$.
\end{proposition}

Observe, as a consistency check, that $\lambda^2 = \pi K(\cdot,\cdot)$ for a
simply connected domain, so that the isotropic coefficient above is
$-\frac{\pi}{2}K$, which is $\frac18 \Delta r$: this is exactly the
$\theta$-average computed in lemma~\ref{lm:dipole-level}. It is also the identity of corollary~\ref{cor:suita}.

\begin{proof}
  Write $z_{1,2}=Z\pm u/2$. By conformal invariance of $g$ and by the
  transformation rule of the Robin function (corollary~\ref{cor:conformal-r}),
  \[
    \Xi(z_1,z_2)= \underbrace{\frac12 \log \frac{(1-|f(z_1)|^2)(1-|f(z_2)|^2)}{|1-f(z_1)\overline{f(z_2)}|^2}}_{T_1}
    + \underbrace{\log \left| \frac{f(z_1)-f(z_2)}{(z_1-z_2)\sqrt{f'(z_1)f'(z_2)}} \right|}_{T_2} .
  \]
  The first term is a metric term, the second a distortion term.

  The argument of the first logarithm equals $1-\varrho_{\mathrm{ph}}^2$, where
  $\varrho_{\mathrm{ph}}$ is the pseudo-hyperbolic distance between $f(z_1)$ and
  $f(z_2)$; and $\varrho_{\mathrm{ph}} = \lambda(Z)|u| + O(|u|^2)$. Hence
  $T_1 = -\tfrac12 \lambda(Z)^2|u|^2+O(|u|^4)$, a term which depends on $|u|$ only
  and not on the orientation of the dipole.

  For the second term, expand numerator and denominator separately:
  \[
    f(Z+\tfrac u2)-f(Z-\tfrac u2) = u\,f'(Z)\Bigl( 1 + \frac{u^2}{24}\frac{f'''}{f'} + O(u^4) \Bigr),
  \]
  \[
    \tfrac12 \Bigl( \log f'(Z+\tfrac u2) + \log f'(Z-\tfrac u2) \Bigr)
    = \log f'(Z) + \frac{u^2}{8} \Bigl( \frac{f'''}{f'}-\Bigl(\frac{f''}{f'}\Bigr)^{2} \Bigr) + O(u^4).
  \]
  Taking the quotient, the terms in $f'''/f'$ combine into
  $\bigl( \tfrac1{24}-\tfrac18 \bigr)\frac{f'''}{f'} = -\tfrac{1}{12}\frac{f'''}{f'}$,
  and there remains $+\tfrac18 (f''/f')^2$; that is exactly $-\tfrac{1}{12}S_f$.
  So $T_2 = -\tfrac1{12}\Ree(S_f(Z)u^2)+O(|u|^4)$, and $H=\log|u|+T_1+T_2$ by
  \eqref{eq:H-dipole}.
\end{proof}

When $D$ is a disc, $f$ is a Möbius transformation, whose Schwarzian derivative
vanishes identically, so the anisotropic term disappears: at this order the
dipole behaves like a particle in a central field for the hyperbolic metric, and
the angular momentum $M$ of the next paragraph is conserved exactly, not just at
this order. For a generic domain $S_f$ does not vanish, and since
$u^2=(\xi^2-\eta^2)+2i\xi\eta$, the term $\Ree(S_fu^2)$ depends on the direction
of the dipole axis: it is a shear which acts as a torque on the dipole, destroys
the conservation of angular momentum, and is the natural suspect for the
destruction of global integrability.
\begin{question}
  For a domain which is not rotationally invariant, is the dipole Hamiltonian non-integrable? A natural strategy is to test meromorphic non-integrability by the Morales--Ramis criterion \cite{MoralesRamis:2001}, examining the identity component of the differential Galois group of the variational equation along collision orbits: by proposition~\ref{prop:schwarz}, this group ``should'' be non-commutative, because of the anisotropic forcing by the Schwarzian derivative. See~\cite{Fejoz:2026} for another strategy, developed for the equal-intensity $2$-vortex problem.
\end{question}

\begin{theorem}[The limit map]
  \label{thm:limitmap}
  Let $\phi\in(0,\pi)$ be the angle of the outgoing dipole ray with the positive
  tangent to $\partial_jD$, so that $\cos\phi=\sin\theta$. Then, to leading order
  in $\eps$:
  \begin{enumerate}[label=(\roman*)]
    \item \emph{(Angle law)} $\theta\mapsto\pi-\theta$; the
      normal component of the velocity is reversed and the tangential component
      preserved.
    \item \emph{(Position law)} The reflection occurs not
      at the point of arrival $s_0$ but at
      \[
        s_0+\mathsf f(\phi)\,L_j,\qquad
        \mathsf f(\phi)=\frac{1}{4(Y_1^\infty)^2+1}
          =\frac12\left(1-\frac{\cos\phi}{\sqrt{1+\cos^2\phi}}\right)
      \]
      in the positive orientation.
    \item \emph{(Return time.)} In the rescaled time $\tau=t/\eps$ the time
      between consecutive wall events is
      $\ell(s,\phi)+\dfrac{L_j}{2\sqrt{1+\cos^2\phi}}$, the chord length plus the
      slide time, both of order one.
  \end{enumerate}
  Thus the limit dynamics is the suspension of the map
  \[
    T_0=\mathcal S_{\mathsf f}\circ T_{\mathrm{billiard}},\qquad
    \mathcal S_{\mathsf f}(s,\phi)=(s+\mathsf f(\phi)L_j,\ \phi),
  \]
  on the annulus $\mathcal A=\partial D\times(0,\pi)$. Since $\mathcal S_{\mathsf f}$
  preserves every circle $\{\phi=\mathrm{const}\}$, $T_0$ is exact
  area-preserving for the billiard measure $\sin\phi\,\ddd s\wedge \ddd\phi$ and
  has the intersection property.
\end{theorem}

\begin{proof}
  (i) $\eta=Y_1-Y_2$ is preserved through the split, the slide and the
  re-formation, so the outgoing axis $\theta'$ satisfies $\sin\theta'=\sin\theta$;
  since the outgoing dipole must move away from the wall, $\cos\theta'<0$, forcing
  $\theta'=\pi-\theta$. The velocity goes from $(\sin\theta,-\cos\theta)$ to
  $(\sin\theta,\cos\theta)$.

  (ii) The vortices meet after a time $T$ with
  $\frac{T}{2\eps Y_1}+\frac{T}{2\eps Y_2}=L_j$, i.e.\
  $T=2\eps L_jY_1Y_2/(Y_1+Y_2)=\eps L_j/\bigl(2(Y_1+Y_2)\bigr)$ using
  $Y_1Y_2=\frac14$. The $+1$ vortex has then covered
  $\frac{T}{2\eps Y_1}=\frac{L_j}{4Y_1(Y_1+Y_2)}=\frac{L_j}{4Y_1^2+1}$ in the
  positive direction. Substituting \eqref{eq:heights},
  $4Y_1^2+1=2\sqrt{1+\sin^2\theta}\bigl(\sqrt{1+\sin^2\theta}+\sin\theta\bigr)$,
  and rationalising gives the stated $\mathsf f$.

  (iii) $Y_1+Y_2=\sqrt{1+\sin^2\theta}$, and $T/\eps$ is the slide time.
\end{proof}

\begin{remark}[Twist of the limit map]
  \label{rem:twist-limit}
  The limit map $T_0$ of theorem~\ref{thm:limitmap} is the composition of the
  billiard map of $D$ with the shift $\mathcal S_{\mathsf f}$, which slides the
  reflection point by $\mathsf f(\phi)L_j$ without changing the angle. This shift
  is a genuine twist: since $\mathsf f$ is a function of $\phi$ alone,
  \[
    \frac{\ddd \mathsf f}{\ddd\phi}(\phi) = \frac{\sin\phi}{2\bigl(1+\cos^2\phi\bigr)^{3/2}}\ >\ 0
    \qquad (0<\phi<\pi),
  \]
  so $\mathcal S_{\mathsf f}$ is an exact area-preserving twist map of the
  annulus $\mathcal A$, with twist bounded below on any $[\phi_0,\pi-\phi_0]$. In
  particular, if $D$ is a disc, the billiard map is itself integrable and
  commutes with rotations, so $T_0$ is an integrable twist map, whose invariant
  circles are the circles $\{\phi=\mathrm{const}\}$. Breaking the rotational
  symmetry of $D$ perturbs $T_0$ within exact area-preserving maps, and Moser's
  invariant curve theorem preserves most of these circles; the dipole in a
  near-circular domain is therefore, in the limit $\eps\to0$, a nearly
  integrable twist map, and there is no reason for the full system to be
  integrable.
\end{remark}

\begin{corollary}
  \label{cor:dipole}
  As $\eps\to0$ the rescaled flow on $\Sig_E$ converges to rectilinear transport on $D$: $Z'=-iu/|u|$ with $u$ frozen. Every orbit reaches a neighbourhood of $\partial D$ in rescaled time $O(\mathrm{diam}\,D)$, having turned by $O(\eps^{2})$.
\end{corollary}

Thus for $\Gamma=0$ there is no fast angle to average, no adiabatic invariant of the kind of corollary~\ref{cor:adiabatic} --- the near-conservation of $\rho$ coming from energy conservation, as in proposition~\ref{prop:dipole-sep}, and not from averaging --- and no interior KAM tori: the limiting problem is a transport problem, and the interesting dynamics is the wall interaction, where a dipole is deflected and slides along the boundary rather than reflecting specularly. That is a question of geometric optics, or of billiard type.

\medskip It is worth noticing that the dipole is a particular case of the opposite-sign problem, so that the global existence result of theorem~\ref{thm:HZ} applies to it: on almost every energy level, and not only for very negative energies, there is a periodic orbit. The elementary compactness estimate behind it reads, in the dipole normalization, as follows, and is the quantitative form of proposition~\ref{prop:dipole-sep}. Using lemma~\ref{lm:asymptotics} near the boundary,
\[H = \frac{1}{2} \log \frac{|z_2-z_1|^2\, \delta(z_1)\, \delta(z_2)}{|z_2-z_1^*|^2} + H_1,\]
with $H_1$ smooth up to $\partial D$; identifying the boundary with the real axis, $|z_2-\bar z_1|^2 = |z_2-z_1|^2+4\delta(z_1)\delta(z_2)$, so that, setting $\mathsf u = |z_2-z_1|^2$ and $\mathsf v = 4\delta(z_1)\delta(z_2)$,
\[e^{2H} \sim \frac{\mathsf u\,\mathsf v}{\mathsf u+\mathsf v}, \qquad \mbox{i.e.} \qquad \frac{1}{\mathsf u}+\frac{1}{\mathsf v} \sim e^{-2H} ,\]
whence $\mathsf u, \mathsf v \gtrsim e^{2E}$ on the level $H=E$: the two vortices stay away from one another \emph{and} from the boundary, quantitatively.
We finally mention two nearly integrable regimes for the dipole, which are the
two natural ways of making the problem perturbative: one may perturb the domain
at fixed energy, or fix the domain and let the energy tend to $-\infty$. The
first route rests on the following exact and pleasant fact.

\begin{proposition}[Conformal form of the perturbation]
  \label{prop:conformal}
  Let $\Psi:D_0\to D_\delta$ be conformal with $\Psi=\mathrm{id}+O(\delta)$ in
  $C^\infty(\overline{D_0})$. Then, exactly,
  \[
    H_{D_\delta}\circ(\Psi\times\Psi)=H_{D_0}
      +\tfrac12\log\bigl|\Psi'(z_1)\Psi'(z_2)\bigr| ,
  \]
  the perturbation being a \emph{sum of single-particle terms} of size
  $O(\delta)$, while $(\Psi\times\Psi)^*\omega
  =\sum_k\gamma_k|\Psi'(z_k)|^2\ddd x_k\wedge \ddd y_k=\omega+O(\delta)$.
\end{proposition}

\begin{proof}
  $g$ is a conformal invariant and, by corollary~\ref{cor:conformal-r}, $r_{D_\delta}(\Psi(z))=r_{D_0}(z)+\log|\Psi'(z)|$;
  the $g$ terms cancel because $\gamma_1\gamma_2 g$ is unchanged, and the
  $r$ terms give $+\frac{\gamma_k^2}{2}$ times the respective
  $\log|\Psi'|$, which is $+\frac12\log|\Psi'(z_1)\Psi'(z_2)|$ when
  $\gamma_k^2=1$. (For instance, if $D_0=\Dd$ and $\Psi(z)=Rz$, then $r_{R\Dd}(Rz)=\log R+r_{\Dd}(z)$.)
\end{proof}

Composing with Moser's trick in each factor to restore $\omega$, one obtains a
symplectomorphism $O(\delta)$-close to the identity, conjugating the dynamics in
$D_\delta$ to an $O(\delta)$ perturbation of the dynamics in $D_0$. Because
$\Sig_E$ is compact here (theorem~\ref{thm:compact}), this is a perturbation on a
compact set, with no end to control --- which is precisely what was impossible in
the same-sign problem.

Take now for $D_0$ a disc. The dipole in a disc is Liouville integrable: besides
the energy, the rotational symmetry of $D_0$ gives the angular momentum
\[
  M = \tfrac12 \bigl( \gamma_1|z_1|^2+\gamma_2|z_2|^2 \bigr) = \tfrac12 \bigl( |z_1|^2-|z_2|^2 \bigr),
\]
which Poisson-commutes with $H_{\Dd}$, and the two are independent away from the
symmetric configurations.

\begin{theorem}[Near-circular domains; scheme]
  \label{thm:kamA}
  Let $D_0$ be a disc, fix $E$, and suppose the integrable system
  $(H_{D_0},M)$ is isoenergetically nondegenerate on an open subset $\mathcal W$
  of $\Sig_E$. Then for every $C^\infty$ domain $D_\delta$ with $\delta$ small,
  the dipole in $D_\delta$ has invariant Lagrangian $2$-tori filling $\mathcal W$
  up to relative measure $O(\sqrt\delta)$. They separate the compact
  $3$-manifold $\Sig_E$, so that the orbits caught between two of them are
  confined for all time; and between consecutive tori Poincaré--Birkhoff and
  Aubry--Mather give periodic orbits and Mather sets.
\end{theorem}

What has to be supplied is the nondegeneracy hypothesis, and this is where the
second route helps. In the regime $E\to-\infty$ the frequency ratio of
$(H_{\Dd},M)$ is, at leading order, the rotation number of the limit map $T_0$ of
theorem~\ref{thm:limitmap}, whose twist is positive by
remark~\ref{rem:twist-limit}. \emph{So the twist of the limit map is the
nondegeneracy of the integrable disc system.} We regard this as strong evidence
for the hypothesis at small $\eps$, and as a computation that should be carried
out directly --- action--angle variables for $(H_{\Dd},M)$ are explicit enough to
attempt --- rather than as a proof.

Alternatively, one may fix a strictly convex $D$ and let $\eps\to0$. By
theorem~\ref{thm:limitmap} the return map is $T_\eps=T_0+o(1)$ on
$\phi\in[\phi_0,\pi-\phi_0]$, with $T_0=\mathcal S_{\mathsf f}\circ T_{\mathrm{billiard}}$
exact area-preserving. For $D$ a $\delta$-perturbation of a disc, $T_0$ is a
$\delta$-perturbation of an integrable twist map, so Moser's invariant curve
theorem gives invariant circles of $T_0$, hence quasi-periodic orbits of the
\emph{limit} dynamics; and $T_\eps$ inherits them provided the error in
theorem~\ref{thm:limitmap} is small in $C^{\ell}$.

Two caveats should be stated. The order in $\eps$ of that error is not
established: we expect $O(\eps)$, coming from the curvature correction to the
wall event, rather than $O(\eps^2)$; numerical experiments are consistent with $O(\eps)$. On
numerical solutions in the disc (with $\eps=2.5\cdot10^{-3}$) and in the domain of
figure~\ref{fig:dipole} (with $\eps=10^{-2}$), for angles $\phi$ between $0.6$ and
$\pi-0.6$, the shift along the boundary agrees with $\mathsf f(\phi)L_j$ within
$2\cdot10^{-3}L_j$, and the angles of arrival and departure agree within
$10^{-3}$; in the disc, the discrepancy on the shift is divided by about four
when $\eps$ goes from $10^{-2}$ to $2.5\cdot10^{-3}$. And the convergence is not uniform as
$\phi_0\to0$, i.e. for dipoles arriving almost tangentially at the wall. The
second route therefore gives the better picture and the weaker theorem, and the
first route the reverse; both are stated here as programmes rather than as
results.

\subsection*{Acknowledgements}
We thank the participants of the Groupe de travail de Systèmes dynamiques in Paris Dauphine, to
whom the first author presented part of this work in 2025, for their remarks.

\subsection*{Funding}
J.~F. was partially supported by the ANR project CoSyDy (ANR-CE40-0014). The
authors were otherwise supported by their home institutions: Université Paris
Dauphine -- PSL and Observatoire de Paris -- PSL (J.~F.), and the University of
Toronto Mississauga (Q.~W.).

\subsection*{Data availability}
No datasets were generated or analysed in this study. The figures were produced
by numerical integration of the equations of motion; the corresponding scripts are
available from the corresponding author upon reasonable request.

\subsection*{Conflict of interest}
The authors declare that they have no conflict of interest.

\bibliographystyle{plain}
\bibliography{2_vortices_in_a_domain}

\end{document}